\documentclass[11pt, reqno]{amsart}

\usepackage[usenames,dvipsnames]{color}

\usepackage[OT2, T1]{fontenc}
\usepackage{url}
\usepackage{amsmath}
\usepackage{array}
\usepackage{graphicx}
\usepackage{amssymb}
\usepackage{amsthm}
\definecolor{link}{RGB}{11,0,128}
\usepackage[colorlinks=true,citecolor=NavyBlue,linkcolor=OliveGreen,urlcolor=link]{hyperref}
\usepackage{hyperref}
\usepackage{paralist}
\usepackage{colonequals}			
\usepackage{color}
\usepackage{mathabx}			
\usepackage{amscd}
\usepackage[all,cmtip]{xy}			
\usepackage{sseq}				
\usepackage{verbatim}
\usepackage{parskip}
\usepackage{microtype}
\usepackage[in]{fullpage}
\usepackage{mathrsfs} 			
\usepackage{cleveref}
\usepackage{enumitem}
\usepackage{moreenum}			
\usepackage[alphabetic,lite]{amsrefs} 	
\usepackage{stmaryrd}			
\usepackage[foot]{amsaddr}
\usepackage{subfiles}
\usepackage{ragged2e}

\newcommand{\gA}{\alpha}

\newcommand{\bA}{\mathbb{A}}

\newcommand{\bC}{\mathbb{C}}

\newcommand{\bF}{\mathbb{F}}
\newcommand{\bG}{\mathbb{G}}

\newcommand{\bN}{\mathbb{N}}
\newcommand{\bP}{\mathbb{P}}
\newcommand{\bQ}{\mathbb{Q}}
\newcommand{\bR}{\mathbb{R}}

\newcommand{\bZ}{\mathbb{Z}}

\newcommand{\bbB}{\mathbf{B}}

\newcommand{\bbG}{\mathbf{G}}
\newcommand{\bbH}{\mathbf{H}}

\newcommand{\bbP}{\mathbf{P}}

\newcommand{\cB}{\mathcal{B}}

\newcommand{\cE}{\mathcal{E}}
\newcommand{\cF}{\mathcal{F}}
\newcommand{\cG}{\mathcal{G}}

\newcommand{\cO}{\mathcal{O}}
\newcommand{\cP}{\mathcal{P}}

\newcommand{\cT}{\mathcal{T}}

\newcommand{\fm}{\mathfrak{m}}

\newcommand{\fp}{\mathfrak{p}}

\newcommand{\sF}{\mathscr{F}}

\newcommand{\sL}{\mathscr{L}}

\newcommand{\sO}{\mathscr{O}}
\newcommand{\sP}{\mathscr{P}}

\newcommand{\sR}{\mathscr{R}}

\newcommand{\sV}{\mathscr{V}}

\newcommand{\sX}{\mathscr{X}}

\DeclareMathOperator{\Aut}{Aut}		
\DeclareMathOperator{\Br}{Br}			
\DeclareMathOperator{\Char}{char}		
\DeclareMathOperator{\depth}{depth}	
\DeclareMathOperator{\Ext}{Ext}		
\DeclareMathOperator{\Frac}{Frac}		
\DeclareMathOperator{\Gal}{Gal}		
\DeclareMathOperator{\GL}{GL}		
\DeclareMathOperator{\Gr}{Gr}			
\DeclareMathOperator{\Hom}{Hom}		
\DeclareMathOperator{\Isom}{Isom}		
\DeclareMathOperator{\Ker}{Ker}		
\DeclareMathOperator{\Lie}{Lie}		
\DeclareMathOperator{\Out}{Out}		
\DeclareMathOperator{\PGL}{PGL}		
\DeclareMathOperator{\Pic}{Pic}		
\DeclareMathOperator{\rad}{rad}		
\DeclareMathOperator{\corad}{corad}		
\DeclareMathOperator{\Res}{Res}		
\DeclareMathOperator{\rk}{rk}			
\DeclareSymbolFont{cyrletters}{OT2}{wncyr}{m}{n}
\DeclareMathSymbol{\Sha}{\mathalpha}{cyrletters}{"58}	
\DeclareMathOperator{\SL}{SL}			
\DeclareMathOperator{\SO}{SO}		
\DeclareMathOperator{\Sp}{Sp}		
\DeclareMathOperator{\Spec}{Spec}		

\newcommand{\ad}{\mathrm{ad}}			

\newcommand{\ce}{\colonequals}

\newcommand{\der}{\mathrm{der}}		

\newcommand{\eps}{\varepsilon}
\newcommand{\et}{\mathrm{\acute{e}t}}	
\newcommand{\fppf}{\mathrm{fppf}}		

\newcommand{\gp}{{\mathrm{gp}}}		
\newcommand{\hra}{\hookrightarrow}
\renewcommand{\i}{^{-1}}
\newcommand{\id}{\mathrm{id}}			
\newcommand{\Inf}{\mathrm{inf}}		

\newcommand{\isomto}{\overset{\sim}{\longrightarrow}}

\newcommand{\llb}{\llbracket}			
\newcommand{\llp}{(\!(}			
	
\newcommand{\ov}{\overline}
\providecommand{\p}[1]{\left(#1\right)}

\newcommand{\ra}{\rightarrow}

\newcommand{\red}{{\mathrm{red}}}		

\newcommand{\rrb}{\rrbracket}			
\newcommand{\rrp}{)\!)}			
\newcommand{\sep}{\mathrm{sep}}		

\newcommand{\sh}{\mathrm{sh}}		
\newcommand{\sm}{\mathrm{sm}}			

\newcommand{\surjects}{\twoheadrightarrow}
\newcommand{\tensor}{\otimes} 			
\newcommand{\tors}{\mathrm{tors}}		

\newcommand{\wt}{\widetilde}
\newcommand{\xra}{\xrightarrow}
\newcommand{\Zar}{\mathrm{Zar}}		

\newcommand{\csub[1]}{\centering\subsection{\text{#1} \nopunct} \hfill \justify}

\providecommand{\up}[1]{{\upshape(}#1{\upshape)}}
\providecommand{\uref}[1]{{\upshape\ref{#1}}}
\providecommand{\uS}{{\upshape\S}}
\providecommand{\f}[2]{\frac{#1}{#2}}

\renewcommand{\b}{\textbf}
\providecommand{\ucolon}{{\upshape:} }
\providecommand{\uscolon}{{\upshape;} }

\newcommand{\brems}{\begin{rems} \hfill \begin{enumerate}[label=\b{\thenumberingbase.},ref=\thenumberingbase]}
\newcommand{\remi}{\addtocounter{numberingbase}{1} \item}
\newcommand{\erems}{\end{enumerate} \end{rems}}
\newcommand{\begs}{\begin{egs} \hfill \begin{enumerate}[label=\b{\thenumberingbase.},ref=\thenumberingbase]}

\newcommand{\eegs}{\end{enumerate} \end{egs}}
\newcommand{\m}{\item}
\newcommand{\bsm}{\begin{smallmatrix}}
\newcommand{\esm}{\end{smallmatrix}}
\newcommand{\blem}{\begin{lemma}}
\newcommand{\elem}{\end{lemma}}

\newcommand{\bconj}{\begin{conj}}
\newcommand{\econj}{\end{conj}}
\newcommand{\bprob}{\begin{Problem}}
\newcommand{\eprob}{\end{Problem}}

\newcommand{\bq}{\begin{Q}}
\newcommand{\eq}{\end{Q}}
\newcommand{\benum}{\begin{enumerate}[label={{\upshape(\alph*)}}]}
\newcommand{\benuma}{\begin{enumerate}[label={{\upshape(\arabic*)}}]}
\newcommand{\benumr}{\begin{enumerate}[label={{\upshape(\roman*)}}]}
\newcommand{\eenum}{\end{enumerate}}
\newcommand{\bitem}{\begin{itemize}}
\newcommand{\eitem}{\end{itemize}}
\newcommand{\bc}{}
\newcommand{\bd}{\begin{defn}}
\newcommand{\ed}{\end{defn}}

\newcommand{\beg}{\begin{eg}}
\newcommand{\eeg}{\end{eg}}

\newcommand{\bcl}{\begin{claim}}
\newcommand{\ecl}{\end{claim}}
\newcommand{\lab}{\label}

\newcommand{\x}{\text}

\newcommand{\q}{\quad}
\providecommand{\qxq}[1]{\quad\text{#1}\quad}
\providecommand{\qx}[1]{\quad\text{#1}}

\newcommand{\qq}{\quad\quad}
\newcommand{\qqq}{\quad\quad\quad}

\newcommand{\tst}{\textstyle}
\newcommand{\ba}{\begin{aligned}}
\newcommand{\ea}{\end{aligned}}
\newcommand{\be}{\begin{equation}}
\newcommand{\ee}{\end{equation}}
\newcommand{\bpf}{\begin{proof}}
\newcommand{\epf}{\end{proof}}
\newcommand{\bthm}{\begin{thm}}
\newcommand{\ethm}{\end{thm}}
\newcommand{\bprop}{\begin{prop}}
\newcommand{\eprop}{\end{prop}}
\newcommand{\bcor}{\begin{cor}}
\newcommand{\ecor}{\end{cor}}
\newcommand{\brem}{\begin{rem}}
\newcommand{\erem}{\end{rem}}
\newcommand*{\QED}{\hfill\ensuremath{\qed}}

\usepackage{aliascnt}
\newaliascnt{numberingbase}{subsubsection}
\numberwithin{equation}{numberingbase}

\newtheoremstyle{thms}{0.2em}{0.2em}{\itshape}{}{\bfseries}{.}{ }{}
\theoremstyle{thms}
\newtheorem{conj}[numberingbase]{Conjecture}

\newtheorem{cor}[numberingbase]{Corollary}

\newtheorem{lemma}[numberingbase]{Lemma}

\newtheorem{prop}[numberingbase]{Proposition}

\newtheorem{Q}[numberingbase]{Question}
\newtheorem{thm}[numberingbase]{Theorem}

\newtheoremstyle{claims}{0.2em}{0.2em}{}{}{\itshape}{.}{ }{}
\theoremstyle{claims}
\newtheorem{claim}[equation]{Claim}

\newtheoremstyle{defs}{0.2em}{0.2em}{}{}{\bfseries}{.}{ }{}
\theoremstyle{defs}
\newtheorem{defn}[numberingbase]{Definition}

\newtheorem{eg}[numberingbase]{Example}
\newtheorem*{egs}{Examples}
\newtheorem{rem}[numberingbase]{Remark}

\newtheorem*{rems}{Remarks}

\Crefname{claim}{Claim}{Claims}
\Crefname{bclaim}{Claim}{Claims}
\Crefname{sublemma}{Lemma}{Lemmas}
\Crefname{conj}{Conjecture}{Conjectures}
\Crefname{cor}{Corollary}{Corollaries}
\Crefname{defn}{Definition}{Definitions}
\Crefname{eg}{Example}{Examples}
\Crefname{prop}{Proposition}{Propositions} 
\Crefname{Q}{Question}{Questions}
\Crefname{rem}{Remark}{Remarks}
\Crefname{thm}{Theorem}{Theorems}
\Crefname{Theorem}{Theorem}{Theorems}
\Crefname{variant}{Variant}{Variants}

\theoremstyle{thms}
\newtheorem{thm-tweak}[subsection]{Theorem}
\Crefname{thm-tweak}{Theorem}{Theorems}
\newtheorem{lemma-tweak}[subsection]{Lemma}
\Crefname{lemma-tweak}{Lemma}{Lemmas}
\newtheorem{cor-tweak}[subsection]{Corollary}
\Crefname{cor-tweak}{Corollary}{Corollaries}
\newtheorem{prop-tweak}[subsection]{Proposition}
\Crefname{prop-tweak}{Proposition}{Propositions} 
\newtheorem{conj-tweak}[subsection]{Conjecture}
\Crefname{conj-tweak}{Conjecture}{Conjectures} 

\theoremstyle{defs}
\newtheorem{defn-tweak}[subsection]{Definition}
\Crefname{defn-tweak}{Definition}{Definitions}
\newtheorem{eg-tweak}[subsection]{Example}
\Crefname{eg-tweak}{Example}{Examples}
\newtheorem*{rems-tweak}{Remarks}
\newtheorem{rem-tweak}[subsection]{Remark}
\Crefname{rem-tweak}{Remark}{Remarks}

\newtheoremstyle{subsection-tweak}
   {2pt}
   {3pt}%
   {}
   {}%
   {\bfseries}
   {}%
   {.5em}
   {\thmnumber{\@{#1}{}\@{#2}.}%
    \thmnote{~{\bfseries#3.}}}    
    
\theoremstyle{subsection-tweak}
\newtheorem{pp}[numberingbase]{}
\newcommand{\bpp}{\begin{pp}}
\newcommand{\epp}{\end{pp}}

\theoremstyle{subsection-tweak}
\newtheorem{pp-tweak}[subsection]{}

\makeatletter
\def\@tocline#1#2#3#4#5#6#7{
    \begingroup 
    \@ifempty{#4}{}{}

    \parindent\z@ \leftskip#3\relax \advance\leftskip\@tempdima\relax
    #5\hskip-\@tempdima
      \ifcase #1
       \or\or \hskip 2em \or \hskip 1em \else \hskip 3em \fi%
      #6\nobreak\relax
    \dotfill\hbox to\@pnumwidth{\@tocpagenum{#7}}\par
    \nobreak
    \endgroup
 }
 \def\l@section{\@tocline{1}{0pt}{1pc}{}{}}

\renewcommand{\tocsection}[3]{%
  \indentlabel{\@ifnotempty{#2}{\makebox[1.3em][l]{%
    \ignorespaces#1 \bfseries{#2}.\hfill}}}\bfseries{#3}
    \vspace{-5pt}}

\renewcommand{\tocsubsection}[3]{%
  \indentlabel{\@ifnotempty{#2}{\hspace*{-0.5em}\makebox[2.1em][l]{%
    \ignorespaces#1#2.\hfill}}}#3
    \vspace{-5pt}}
   
\makeatother 

\makeatletter 
\newcommand\appendix@section[1]{%
  \refstepcounter{section}%
  \orig@section*{Appendix \@Alph\c@section. #1}%
}
\let\orig@section\section
\g@addto@macro\appendix{\let\section\appendix@section}
\makeatother

\makeatletter
\@namedef{subjclassname@2020}{%
  \textup{2020} Mathematics Subject Classification}
\makeatother

\author{K\k{e}stutis \v{C}esnavi\v{c}ius}
\address{CNRS, Universit\'{e} Paris-Saclay,   Laboratoire de math\'{e}matiques d'Orsay, F-91405, Orsay, France}
\email{kestutis@math.u-psud.fr}

\date{\today}

\usepackage{stackengine}
\usepackage{tikz-cd} 
\usepackage{rotating}

\begin{document}

\title{Problems about torsors over regular rings \\ \scriptsize{with an appendix by Yifei Zhao}}

\date{March 29, 2023}

\maketitle

\vspace{-50pt}

\hypersetup{
    linktoc=page,     
}

\renewcommand*\contentsname{}
\q\\
\tableofcontents


\section{Introduction}

The goal of this survey is to discuss a web of conjectures and questions about torsors under reductive groups over regular rings. Some of these are well-known major problems in the field that have stood the tests both of time and of partial results (some recent) by multiple authors. Some others appear to have avoided the spotlight, even though they are close in spirit or even have direct links to the better known of these conjectures. In spite of multiple surveys that some of these problems have already received, we believe that it is worthwhile to discuss them together in the pages that follow with the hope that highlighting their common aspects may eventually lead to further progress. 

Indeed, even though these problems concern torsors, key progress on them involved establishing structural results about regular rings themselves, concrete examples being the Popescu approximation 
(\Cref{thm:Popescu-2} below), the Geometric Presentation Theorem (\Cref{thm:geometric-presentation} below), or the Lindel lemma (\Cref{cor:baby-Lindel} below). These structural results are useful in many contexts, so finding fruitful approaches to torsor problems tends to bring general insights into the geometry of regular rings. This amplifies the significance of problems about torsors, although they are captivating already for the elegance and simplicity of their statements. Heuristically, this geometric approach is suggested by the difficulty of ``enlarging'' the regular rings at hand, for instance, by passing to Henselizations or completions, because \emph{a priori} this may trivialize the torsors one is studying---therefore, one is forced to ``shrink'' the rings instead by studying their fine geometric structure.

The conjectures in question almost exclusively concern regular \emph{local}  rings $R$,\footnote{We expect them to stay true for regular \emph{semilocal} rings, but for the sake of focus we chose to neglect this aspect below. Many of their known special cases are established in this generality in the indicated references.} so they naturally split into three cases of increasingly arithmetic flavor: when $R$ is of equal characteristic, when it is of mixed characteristic but unramified, and when it is of mixed characteristic and ramified (see \S\ref{sec:trichotomy}). Thanks to Popescu approximation, the equal characteristic case essentially concerns local rings of smooth varieties over a field and tends to be the most approachable. Likewise, the mixed characteristic unramified case essentially concerns local rings of smooth schemes over $\bZ$, so it tends to be similar to the equal characteristic case except for new and often rather delicate geometric subtleties of arithmetic flavor. Finally, the mixed characteristic ramified case has so far remained almost entirely out of reach, and our understanding of the geometry of arbitrary ramified regular local rings appears to still be limited. 

As for reductive groups themselves, the simplest are the commutative ones, that is, tori: for them, one typically reduces the problem at hand either to the most basic case of $\bG_m$ or to an abelian question that concerns \'{e}tale or flat cohomology. The next simplest class is that of general linear groups $\GL_n$: for them, studying torsors amounts to studying vector bundles---this case already exhibits general nonabelian phenomena, although it also has some simplifying features, for example, Zariski local triviality of the torsors in question. The case of a general reductive group $G$ is the most delicate: the question at hand may already be highly nontrivial for split $G$ (or even for $\GL_n$), 
and it tends to complicate further once arithmetic structure of a non-split $G$ enters the picture. For instance, in general it may be important to know whether $G$ has nontrivial split subtori, nontrivial parabolic subgroups, or even a Borel subgroup. We segregate problems about $\GL_n$-torsors into \S\ref{sec:vector-bundles} and then discuss problems that concern general $G$ in \S\ref{sec:general-groups}. In subsequent \S\S\ref{sec:ring-structure}--\ref{sec:fix-G}, we discuss some of the techniques that have been used for making progress on the conjectures overviewed in \S\S\ref{sec:vector-bundles}--\ref{sec:general-groups}.

Most of the problems we discuss admit analogues for valuation rings, as we sometimes indicate along the way. For this, the idea is that, by the Zariski local uniformization conjecture (a local form of the resolution of singularities), any valuation ring $V$ ought to be a filtered direct limit of its regular local subrings. Although local uniformization remains open in positive and mixed characteristics, this  heuristic suggests precise formulations for expected valuation ring analogues of problems about torsors. These analogues may often be settled directly: the non-Noetherian nature of general valuation rings brings considerable technicalities but, on the other hand, valuation rings have somewhat straight-forward ring-theoretic properties when compared to regular rings, for instance, their prime ideals are linearly ordered and their arithmetic resembles regular rings of dimension $\le 1$.

Before proceeding to the main body of the text, we review some basic facts about torsors, reductive group schemes, and regular rings that we will use throughout without explicit mention.

\subsection*{Acknowledgements}
I thank Vi\stackon[-9.5pt]{\^{e}}{\d{}}n To\'{a}n H\d{o}c for the invitation to contribute to the special issue of Acta Mathematica Vietnamica. I thank Yifei Zhao for the appendix and for helpful comments on the main body of the text. I thank the referee for helpful comments and suggestions. I thank Alexis Bouthier, Jean-Louis Colliot-Th\'{e}l\`{e}ne, Sean Cotner, Roman Fedorov, Ofer Gabber, Ning Guo, Shang Li, Ivan Panin, Federico Scavia, Yifei Zhao, and many others for useful conversations and correspondence related to the subject of this article. This project received funding from the European Research Council under the European Union's Horizon 2020 research and innovation program (grant agreement No.~851146).

\csub[Notation and conventions] \label{pp:conv}
The rings we consider are commutative and unital. A commutative ring $R$ is \emph{local} (resp.,~\emph{semilocal}) if it has a unique (resp.,~finitely many) maximal ideal(s) $\fm \subset R$ (so the zero ring is semilocal but not local). A scheme is \emph{local} (resp.,~\emph{semilocal}) if it is the spectrum of a local (resp.,~semilocal) ring.  A local ring $(R, \fm)$ is \emph{complete} if it is complete for its $\fm$-adic topology, in other words, if every $\fm$-adic Cauchy sequence in $R$ has a unique limit. We say that a local domain $R$ is \emph{of equicharacteristic} (resp.,~\emph{of mixed characteristic}) if its fraction field and its residue field have the same (resp.,~different) characteristics, equivalently, if $R$ contains either $\bQ$ or some $\bF_p$ (resp.,~if $R$ contains no field). We let $k_\fp$ (resp.,~$k_s$) denote the residue field of a prime ideal $\fp \subset R$ of a ring (resp.,~of a point $s \in S$ of a scheme). For a scheme $S$, and an $S$-scheme $X$, we let $X^\sm$ denote its smooth locus (which is automatically open, see \cite{SP}*{Definition \href{https://stacks.math.columbia.edu/tag/01V5}{01V5}}). A \emph{vector bundle} on a scheme $S$ is an $\sO_S$-module that is finite locally free  of finite rank; a vector bundle $\sV$ on $S$ is \emph{stably free} if $\sV \oplus \sO_S^{\oplus n'} \simeq \sO_S^{\oplus n}$ for some $n, n' \ge 0$. Every vector bundle on a semilocal $S$ scheme is free granted that its rank is constant (an automatic condition if $S$ is connected, for instance, if $S$ is local), see~\cite{SP}*{Lemma~\href{https://stacks.math.columbia.edu/tag/02M9}{02M9}}.

These notations and conventions 
are in place throughout the article, including the appendix.

\csub[Basic properties of torsors] \label{pp:basic}

An exhaustive reference for generalities about torsors is Giraud's book \cite{Gir71}. \emph{Op.~cit.}~may appear difficult to navigate at first, so we now recall some salient~points.

\bpp[The pointed set $H^1(S, G)$] \label{pp:basic-pointed}
For a group sheaf $G$ on a site $S$, a \emph{torsor}, synonymously, a \emph{principal homogeneous space} under $G$ (simply, a \emph{$G$-torsor}) is a sheaf $E$ equipped with a right action of $G$ such that  $E$ becomes trivial locally for the topology in question, that is, becomes isomorphic to the trivial torsor given by $G$ equipped with its right translation action. A simple triviality criterion~is 
\[
E(S) \overset{?}{\neq} \emptyset.
\]
The collection of isomorphism classes of $G$-torsors is denoted by $H^1(S, G)$ and is understood to be pointed by the class of the trivial torsor. For abelian $G$, this agrees with the usual derived functor $H^1$, see \cite{Gir71}*{Chapitre III, Remarque 3.5.4}. Any morphism between $G$-torsors is automatically an isomorphism. The automorphism functor of the trivial $G$-torsor is $G$ itself acting via \emph{left} translation. In particular, for a $G$-torsor $E$, its automorphism functor $\Aut_G(E)$ is a group sheaf that is a pure inner form of $G$ (and every pure inner form of $G$ arises in this way).\footnote{A \emph{form} of $G$ is an $S$-group sheaf isomorphic to $G$ locally on $S$, so it corresponds to an element of $H^1(S, \underline{\Aut}_\gp(G))$. A form is \emph{inner} (resp.,~\emph{pure inner}) if this element lifts to $H^1(S, G/Z_G)$ (resp.,~even to $H^1(S, G)$), where $Z_G \subset G$ is the center and the map $G/Z_G \ra \underline{\Aut}_\gp(G)$ is induced by $G$ acting on itself by conjugation.} The $G$-torsors correspond to the $\Aut_G(E)$-torsors via the ``change of origin'' bijection 
\be \label{eqn:change-origin}
H^1(S, G) \cong H^1(S, \Aut_G(E)) \qxq{given by} E' \mapsto \Isom_G(E, E'),
\ee
where $\Isom_G(E, E')$ is the $\Aut_G(E)$-torsor  that parametrizes $G$-torsor isomorphisms between $E$ and $E'$ (so $E$ gets sent to the class of the trivial $\Aut_G(E)$-torsor), see \cite{Gir71}*{Chapitre III, Remarque~2.6.3}. 
When $G$ is abelian, $\Aut_G(E) \cong G$ and the bijection simply subtracts the~class~of~$E$. 
\epp

\bpp[Subgroups and quotients] \label{pp:basic-subgroups}
For a map of group sheaves $G' \ra G$, every $G'$-torsor $E'$ gives rise to a $G$-torsor defined as the contracted product $E \ce E' \times_{G'} G$; this gives a map of pointed sets 
\[
H^1(S, G') \ra H^1(S, G).
\]
Conversely, for an inclusion of group sheaves $G' \subset G$ and a $G$-torsor $E$, the quotient $E/G'$ parametrizes the reductions of $E$ to a $G'$-torsor $E'$. As in the abelian case, by \cite{Gir71}*{Chapitre III, Proposition~3.3.1; Chapitre IV, Remarque 4.2.10}, a short exact sequence of $S$-sheaves
\[
1 \ra G' \ra G \ra G'' \ra 1
\]
with $G'$ and $G$ group sheaves and $G'' \ce G/G'$ produces a functorial exact sequence of pointed sets
\be \label{eqn:coho-seq}
1 \ra G'(S) \ra G(S) \ra G''(S) \ra H^1(S, G') \ra H^1(S, G) \dashrightarrow H^1(S, G'') \dashrightarrow H^2(S, G'),
\ee
where the first (resp.,~second) dashed arrow exists if $G'$ is normal (resp.,~even central) in $G$, 
and exactness means that the kernel of each arrow is precisely the image of the preceding one. 
\epp

\bpp[The scheme case and the representability properties] \label{pp:basic-representability}
For us, $S$ will be a scheme endowed with its fppf topology and $G$ will be an $S$-group algebraic space (or even an $S$-group scheme), and we will consider torsors for the fppf topology. By \cite{SP}*{Lemma \href{https://stacks.math.columbia.edu/tag/04SK}{04SK}}, the $G$-torsors are then automatically representable by algebraic spaces because they are representable fppf locally on $S$ (by base changes of $G$). In contrast, if $G$ happens to even be an $S$-group scheme, then its torsors \emph{need not be} representable by schemes: nonrepresentable torsors exist already when $G$ is an abelian scheme, see \cite{Ray70b}*{Section XIII 3.2, page 200}. However, this point tends to be moot: modulo more demanding technicalities, working with algebraic spaces is often (but not always!) ``just as good.'' 

On the positive side, if $G$ is an $S$-affine $S$-scheme, or merely an $S$-ind-quasi-affine\footnote{\label{foot:ind-quasi-affine}We recall from \cite{SP}*{Definition \href{https://stacks.math.columbia.edu/tag/0AP6}{0AP6}} that a scheme is \emph{ind-quasi-affine} if all of its quasi-compact opens are quasi-affine, and that a morphism is \emph{ind-quasi-affine} if the preimage of every affine open is ind-quasi-affine. By \cite{SP}*{Lemma \href{https://stacks.math.columbia.edu/tag/0AP8}{0AP8}}, ind-quasi-affineness of a morphism is fpqc local on the target. Useful examples of ind-quasi-affine but not quasi-affine (that is, not quasi-compact) schemes are character groups of tori or automorphism groups of reductive groups, see \S\ref{pp:automorphisms} below. 
} $S$-scheme, then its torsors are representable by schemes: the affine case follows from flat descent for quasi-coherent sheaves, the more general quasi-affine case follows from the effectivity of descent for quasi-affine schemes \cite{SP}*{Lemma \href{https://stacks.math.columbia.edu/tag/0247}{0247}}, and the yet more general ind-quasi-affine case is more subtle and follows from the effectivity of fpqc descent for ind-quasi-affine morphisms due to Gabber \cite{SP}*{Lemma \href{https://stacks.math.columbia.edu/tag/0APK}{0APK}}. 

Of course, by fppf descent, torsors, as well as functors like $\Aut_G(E)$ and $\Isom_G(E, E')$ above, inherit properties of $G$ such as relative ((ind-)quasi-)affineness, or finite presentation, or smoothness, etc. In particular, since a smooth morphism of schemes admits a section \'{e}tale locally on the base (see \cite{EGAIV4}*{Corollaire 17.16.3~(ii)}), torsors under a smooth $G$ trivialize \'{e}tale locally on $S$. Similarly, if $G$ is flat and locally of finite presentation, then so are its torsors even for the fpqc topology, so that they all trivialize already fppf locally; consequently, considering fpqc $G$-torsors is ``no more general.''

As for quotients, it is useful to recall that, for a scheme $S$, the fppf sheaf quotient $X/G$ of an $S$-algebraic space $X$ equipped with a \emph{free} action of a flat, locally finitely presented $S$-group algebraic space $G$ is always representable by an $S$-algebraic space, see \cite{SP}*{Lemma \href{https://stacks.math.columbia.edu/tag/06PH}{06PH}}. The map $X \ra X/G$ is a $G$-torsor over $X/G$, so, by descent, it inherits properties of $G \ra S$, for instance, flatness and local finite presentation; in particular, if $X \ra S$ is smooth, then so is $X/G \ra S$, see \cite{SP}*{Lemma~\href{https://stacks.math.columbia.edu/tag/0AHE}{0AHE}}. Deciding whether $X/G$ is a scheme if $X$ and $G$ are schemes is significantly more delicate, for instance, no general result ensures this even when $X$ and $G$ are both affine (with $S$ general). Some situations in which $X/G$ is a scheme is when $X$ is affine and $G$ is either finite locally free (see \cite{SGA3Inew}*{Th\'{e}or\`{e}me~4.1~(iv)}) or a torus (or even reductive over $S$, see the end of \S\ref{pp:reductive}); in both of these cases the quotient $X/G$ is also affine over $S$.

As far as representability by schemes goes, it may be worth to recall that one does not know any example of a scheme $S$ and a smooth, separated, finitely presented $S$-group algebraic space $G$ with connected fibers that would not be a scheme (compare with \cite{FC90}*{Chapter I, Theorem 1.9}). 
\epp

\csub[Basic properties of reductive group schemes and of their torsors] \label{pp:basic-reductive}
We assume that the reader is familiar with the theory of reductive group schemes (so also with reductive groups over a field) described in \cite{SGA3II}, \cite{SGA3IIInew}, and surveyed by Conrad in \cite{Con14}. Nevertheless, we now review the basic aspects that are particularly relevant for studying torsors---in reality, each individual conjecture discussed in \S\ref{sec:general-groups} only requires a small subset of them. 

\bpp[Reductive groups] \label{pp:reductive}
For a scheme $S$, an $S$-group scheme $G$ is \emph{reductive} (resp.,~\emph{semisimple}) if it is smooth and affine over $S$ and its geometric $S$-fibers are connected reductive (resp.,~semisimple) groups, see \cite{SGA3IIInew}*{Expos\'e XIX, D\'{e}finition 2.7}. Basic examples of reductive $S$-groups are the $S$-tori (see \S\ref{sec-multiplicative-type-closure-properties}) and the split $S$-groups such as $\GL_{n,\, S}$, $\PGL_{n,\, S}$, $\SO_{n,\, S}$, etc. Split (and pinned) reductive group schemes are classified combinatorially by root data and every reductive $S$-group is split (and pinned) \'{e}tale locally on $S$, see \cite{SGA3IIInew}*{Expos\'e XXII, D\'{e}finition~1.13, Corollaire~2.3; Expos\'e XXIII, D\'{e}finition~1.1; Expos\'e  XXV, Th\'{e}or\`{e}me~1.1}. The split type of a general reductive $S$-group $G$ is locally constant on $S$, see \cite{SGA3IIInew}*{Expos\'e XXII, Proposition~2.8}. If this type is constant, then $G$ is a form of the split reductive $S$-group scheme of the same type, see \cite{SGA3IIInew}*{Expos\'e XXII, Corollaire 2.3} and also \S\ref{pp:automorphisms} below.

By \cite{SGA3II}*{Expos\'e XVI, Corollaire 1.5 (a)}, for a reductive $S$-group $H$, any $S$-monomorphism $H \hra G$ to a finitely presented $S$-group scheme $G$ is necessarily an immersion (resp.,~a closed immersion if $G$ is separated).\footnote{Note that the case of \emph{loc.~cit.}~that allows a separated $G$ to be merely locally of finite type over a Noetherian $S$ is false, as is pointed out in \cite{Con14}*{Theorem 5.3.5 and below}: the N\'{e}ron lft model of $\bG_m$ gives a counterexample because its relative identity component is an open but not closed subgroup identified with $\bG_m$.} For a (closed) immersion $H \hra G$ of reductive $S$-group schemes, the quotient $G/H$ is an $S$-affine scheme, more precisely, if $G$ is reductive and $H$ is merely its smooth, closed $S$-subgroup with connected $S$-fibers, then $G/H$ is an $S$-affine scheme if and only if $H$ is reductive, see \cite{Alp14}*{Corollary 9.7.7} in the post-publication arXiv version of \emph{op.~cit.} The affineness of $G/H$ generalizes the Matsushima theorem to arbitrary base schemes, see \cite{Alp14}*{Section 9.4}. Moreover, for any $S$-affine $S$-scheme $X$ equipped with a free action of a reductive $S$-group scheme $H$, the quotient $X/H$ is an $S$-affine $S$-scheme, see \cite{Alp14}*{Proposition 9.7.8}.  
\epp

\bpp[Subtori] \label{pp:subtori}
An $S$-subtorus $T \subset G$ is a \emph{maximal torus} of a reductive $S$-group $G$ if, for every geometric $S$-point $\ov{s}$, the base change $T_{\ov{s}}$ is a maximal subtorus of $G_{\ov{s}}$. By \cite{SGA3II}*{Expos\'e~XIV, Corollaire~3.20}, every reductive $S$-group $G$ admits a maximal torus Zariski locally on $S$, in fact, even Zariski semilocally on $S$: any finite set of points of $S$ contained in a single affine open lies in a smaller affine open over which $G$ has a maximal torus. Any maximal torus $T$ is its own centralizer in a reductive $G$, see \cite{SGA3IIInew}*{Expos\'e XIX, Lemme 1.6.2, Section 2.2} (with \cite{EGAIV4}*{Corollaire~17.9.5}), and any two maximal tori of $G$ are conjugate \'{e}tale locally on $S$, see \cite{SGA3II}*{Expos\'e XI, Corollaire 5.4 bis; Expos\'e XII,~Section 1.0}. The normalizer $N_G(T)$ of any $S$-subtorus $T \subset G$ is a closed, $S$-smooth subgroup of $G$ and, if $T$ is maximal, then the quotient 
\[
W \ce N_G(T)/T
\]
is a finite \'{e}tale $S$-group scheme, the \emph{Weyl group} of $T$ (or simply of $G$), see \cite{SGA3II}*{Expos\'e XI, Corollaire 5.3 bis; Expos\'e XII, Th\'{e}or\`{e}me 2.1}. The functor 
that parametrizes the maximal tori of (base changes of) a reductive $S$-group $G$ is an affine, smooth $S$-scheme (see \cite{SGA3II}*{Expos\'e XII, Corollaire 5.4}) that is isomorphic to $G/N_G(T)$ for any maximal $S$-torus $T \subset G$. The commutative reductive $S$-group schemes are precisely the $S$-tori. The $S$-subtori are closed in $G$, see~\S\ref{pp:reductive}.  
 \epp


\bpp[The reductive center and reductive groups built using it] \label{pp:center}
For a scheme $S$, the \emph{center} of a reductive $S$-group $G$ is the kernel of the conjugation map $G \ra \underline{\Aut}_{\gp}(G)$ and is a closed, finitely presented, $S$-flat group subscheme of multiplicative type $Z_G \subset G$, see \cite{SGA3II}*{Expos\'e~XII, Proposition 4.11} or \cite[Theorem 3.3.4]{Con14}. The self-centralizing property of any maximal $S$-torus $T\subset G$ implies that $Z_G \subset T$, more precisely, \cite{SGA3II}*{Expos\'e XII, Th\'{e}or\`{e}me 4.7 d)} (or \cite[Corollary 3.3.6]{Con14}) shows that $Z_G$ is the kernel of the adjoint action of $T$ on the Lie algebra $\Lie(G)$. In turn, $Z_G$ contains a unique maximal central subtorus of $G$ (see \cite[Expos\'e XII, Proposition 1.12]{SGA3II}), which is of formation compatible with base change, is called the \emph{radical} of $G$, and is denoted by $\rad(G)$ as in \cite[Expos\'e XXII, D\'{e}finition~4.3.6]{SGA3IIInew}. A reductive $S$-group $G$ is semisimple (resp.,~\emph{adjoint}) if and only if $\rad(G)$ (resp.,~$Z_G$) is trivial. By \cite{SGA3IIInew}*{Expos\'e~XXII, Proposition~4.3.5~(ii)}, for a reductive $S$-group $G$, the \emph{adjoint quotient} 
\[
G^\ad \ce G/Z_G
\]
is a semisimple adjoint $S$-group scheme.

For a reductive $S$-group $G$, the fppf sheafification of the group subpresheaf that sends an $S$-scheme $S'$ to the commutator 
\[
[G(S'), G(S')] \subset G(S')
\]
is a semisimple $S$-subgroup $G^{\der} \subset G$, the \emph{derived subgroup} of $G$, see \cite{SGA3IIInew}*{Expos\'e XXII, Th\'{e}or\`{e}me 6.2.1} or \cite[Theorem 5.3.1]{Con14}. By the same references, 
\[
\corad(G) \ce G/G^\der
\]
is a torus, the \emph{coradical} of $G$. By \cite[Expos\'e XXII, Section 6.2.3]{SGA3IIInew}, the multiplication map
\begin{equation}
\label{eq-reductive-group-derived-radical}
G^{\der} \times \rad(G) \rightarrow G
\end{equation}
is a central isogeny (see \S\ref{pp:central-isogeny}) whose kernel is finite and of multiplicative type. When studying torsors, one often combines maps like \eqref{eq-reductive-group-derived-radical} with an analysis of resulting long exact cohomology sequences \eqref{eqn:coho-seq} in attempts to reduce to simpler groups. For a semisimple $S$-group $G$, we denote its simply connected cover by $G^{\mathrm{sc}}$, see \Cref{lemm-simply-connected-cover} for a review. A semisimple $S$-group scheme is \emph{simply connected} if the central isogeny $G^{\mathrm{sc}} \ra G$ is an isomorphism.

The adjoint (resp.,~the simply connected) semisimple $S$-group schemes $G$ canonically decompose as
\be \label{eqn:ad-decomposition}
\tst G \cong \prod_i G_i \qxq{with} G_i \cong \Res_{S_i/S} \wt{G}_i,
\ee
where $i$ runs over the types of nonempty connected Dynkin diagrams, $S_i$ is a finite \'{e}tale $S$-scheme (a clopen in the scheme of Dynkin diagrams of $G$), and $\wt{G}_i$ is an adjoint (resp.,~simply connected) semisimple $S_i$-group scheme with simple geometric $S_i$-fibers of type $i$. This decomposition is one reason why adjoint or simply connected groups tend to be somewhat easier to analyze. 
\epp


\bpp[Parabolic subgroups] \label{pp:parabolic}
For a scheme $S$ and a reductive $S$-group $G$, an $S$-subgroup $P \subset G$ is a \emph{parabolic} (resp.,~a \emph{Borel}) if it is $S$-smooth and each of its geometric $S$-fibers contains (resp.,~is) a maximal solvable subgroup of the corresponding geometric $S$-fiber of $G$. Two parabolic $S$-subgroups of $G$ are \emph{of the same type} if they are conjugate Zariski locally on $S$, see \cite{SGA3IIInew}*{Expos\'e XXVI, D\'{e}finition 3.4, Corollaire 5.2} (also for an equivalent definition). A parabolic $S$-subgroup $P$ is closed in $G$, has connected $S$-fibers, is its own normalizer in $G$, and the quotient $G/P$ is a smooth, projective $S$-scheme, see \cite{SGA3IIInew}*{Expos\'e XXVI, Proposition 1.2}.  The functor that parametrizes parabolic subgroups of (base changes of) a reductive $S$-group $G$ is a smooth, projective $S$-scheme, and the subfunctors that parametrize parabolics of a fixed type are clopen in it---such a subfunctor is simply $G/P$ once a parabolic $P \subset G$ of the type in question is fixed, see \cite{SGA3IIInew}*{Expos\'e XXVI, Corollaires~3.5~et~3.6}. For a $G$-torsor $E$ and a parabolic $S$-subgroup $P \subset G$, the quotient $E/P$ is a smooth, projective $S$-scheme, in fact, it is a scheme that parametrizes parabolics of a fixed type of the inner form $\Aut_G(E)$ of $G$, see \cite{SGA3IIInew}*{Expos\'e XXVI, Lemme~3.20}.\footnote{The representability of $E/P$ is quite remarkable because no general result about quotients ensures it, see \S\ref{pp:basic-representability}.}


A parabolic $S$-subgroup $P \subset G$ has its \emph{unipotent radical}, namely, the largest normal $S$-subgroup $\sR_u(P) \subset P$ that is smooth, finitely presented, and whose geometric $S$-fibers are connected and unipotent, see \cite{SGA3IIInew}*{Expos\'e XXVI, Proposition 1.6~(i)}. The 
quotient $P/\sR_u(P)$ is a reductive $S$-group scheme (\emph{loc.~cit.}),~in particular, $\sR_u(P)$ is closed in $P$. There is a canonical filtration 
\[
\dotsc \subset U_{i + 1} \subset U_i \subset \dotsc \subset U_1 \subset U_0 = \sR_u(P)
\]
by normal, closed $S$-subgroups $U_i$ that are smooth, with connected geometric $S$-fibers, such that the successive quotients $U_i/U_{i + 1}$ are vector groups, that is, are associated to finite locally free $\sO_S$-modules, and $[U_i, U_j] \subset U_{i + j + 1}$ for all $i, j \ge 0$, and such that every automorphism of $P$ preserves $\sR_u(P)$ and the $U_i$ and acts linearly on each $U_i/U_{i + 1}$, see \cite{SGA3IIInew}*{Expos\'e XXVI, Proposition~2.1, Section~2.1.2}. In~particular, whenever $S$ is affine the $S$-scheme $\sR_u(P)$ is isomorphic to the affine space associated to some vector bundle on $S$ and
\be \label{eqn:RuP-H1}
H^1(S, \sR_u(P)) = \{*\},
\ee
see \cite{SGA3IIInew}*{Expos\'e XXVI, Corollaire 2.5}.  A \emph{Levi} $S$-subgroup of a parabolic $S$-subgroup $P \subset G$ is any $S$-subgroup $L \subset P$ that maps isomorphically to the reductive $S$-group $P/\sR_u(P)$, so
\[
P \cong \sR_u(P) \rtimes L.
\]
Any such $L$ is its own normalizer in $P$, any two Levis of $P$ are conjugate by a unique section of $\sR_u(P)$, the functor that parametrizes Levis of (base changes of) $P$ is a $\sR_u(P)$-torsor, and for any maximal $S$-torus $T \subset P$ there exists a unique $S$-Levi containing $T$, see \cite{SGA3IIInew}*{Expos\'e XXVI, Proposition 1.6~(ii), Corollaires 1.8 et 1.9}. A parabolic $S$-subgroup need not have an $S$-Levi but, by \eqref{eqn:RuP-H1}, it does whenever $S$ is affine. If a parabolic $S$-subgroup $P \subset G$ has an $S$-Levi $L \subset P$, then it also has a unique \emph{opposite parabolic relative to $L$}: there is a unique parabolic $S$-subgroup $P' \subset G$ such that $P \cap P' = L$, see \cite{SGA3IIInew}*{Expos\'e XXVI, Th\'{e}or\`{e}me 4.3.2 (a), D\'{e}finition 4.3.3}. 

Due to the uniqueness aspects above, especially, the aspect about any automorphism of $P$ preserving the $U_i$ and the structure of the $U_i/U_{i + 1}$, the claims of the preceding paragraph all continue to hold for any form of $P$, that is, for any $S$-group scheme that becomes isomorphic to $P$ fpqc locally on $S$. 

If $G$ is the restriction of scalars $\Res_{S'/S}(G')$ for a finite \'{e}tale cover $S' \ra S$, then the parabolic $S$-subgroups (resp.,~Borel $S$-subgroups; resp.,~maximal $S$-tori) of $G$ are precisely the restrictions of scalars of the parabolic $S'$-subgroups (resp.,~Borel $S'$-subgroups; resp.,~maximal $S'$-tori) of $G'$, as one checks by reducing to when $S'$ is a disjoint union of copies of $S$, see also \cite{Con14}*{Exercise~6.5.7}. 
\epp


\bpp[Torsors under parabolic subgroups] \label{pp:parabolic-torsors}
The vanishing \eqref{eqn:RuP-H1} is used often when studying torsors. For instance, by also using it in the case of inner forms and exploiting \eqref{eqn:change-origin}--\eqref{eqn:coho-seq}, for any ring $R$ and any $R$-Levi $L$ of a parabolic $R$-subgroup $P$ of a reductive $R$-group $G$ one obtains
\be \label{eqn:M-P}
H^1(R, L) \isomto H^1(R, P),
\ee
see \cite{SGA3IIInew}*{Expos\'e XXVI, Corollaire 2.3}, and similarly when $P$ is merely an $R$-form of a parabolic $R$-group, compare with the end of \S\ref{pp:parabolic}. As we now explain, in the case when $R$ is semilocal we also have
\be \label{eqn:Levi-H1}
H^1(R, L) \hra H^1(R, G), \qxq{equivalently,}  H^1(R, P) \hra H^1(R, G).
\ee
Indeed, by twisting by a variable $P$-torsor as reviewed in \eqref{eqn:change-origin}, at the cost of changing $G$ we reduce to showing that the map on $H^1$ has a trivial kernel. The sequence \eqref{eqn:coho-seq} then reduces us to showing that the map $G(R) \surjects (G/P)(R)$ is surjective. However, by \cite{SGA3IIInew}*{Expos\'e XXVI, Corollaire 5.2}, even the map $\sR_u(P)(R)\sR_u(P')(R) \ra (G/P)(R)$ is surjective for any parabolic $R$-subgroup $P' \subset G$ that is opposite to $P$ (see \S\ref{pp:parabolic}).

The injectivity \eqref{eqn:Levi-H1} deserves to be known more widely, for instance, it is in the same spirit as the Witt cancellation theorem for quadratic forms \cite{Bae78}*{Chapter III, Corollary 4.3}.
\epp


\bpp[Totally isotropic reductive groups]\label{pp:totally-isotropic}
A reductive group $G$ over a semilocal affine scheme $S$ is \emph{quasi-split} (resp.,~\emph{isotropic}; resp.,~\emph{anisotropic}) if it has a Borel $S$-subgroup\footnote{\label{foot:quasi-split}Beyond semilocal $S$, quasi-splitness is a slightly more delicate notion, see \cite{SGA3IIInew}*{Expos\'e XXIV, Section~3.9}.} (resp.,~if it has $\bG_{m,\, S}$ as a subgroup; resp.,~if it has no $\bG_{m,\, S}$ as a subgroup). If $S$ is, in addition, connected, then $G$ is anisotropic if and only if it has no proper parabolic $S$-subgroup and $\rad(G)$ has no $\bG_{m,\, S}$ as an $S$-subgroup, see \cite{SGA3IIInew}*{Expos\'e XXVI, Corollaire 6.12}. 

A reductive group $G$ over a scheme $S$ is \emph{totally isotropic} 
if for every $s \in S$, each $\wt{G}_i$ that appears in the canonical decomposition \eqref{eqn:ad-decomposition} of $G^\ad_{\sO_{S,\, s}}$ (over $\Spec(\sO_{S,\, s})$) is isotropic, 
see  \cite{split-unramified}*{Definition~8.1}, equivalently, if each $\wt{G}_i$ has a parabolic subgroup that contains no fiber of $\wt{G}_i$. This condition on $G^\ad_{\sO_{S,\, s}}$ is stable upon replacing $s$ by a generization, so it suffices to consider those $s$ that exhaust the closed points of the members of some affine open cover of $S$, for instance, for semilocal affine $S$ it suffices to consider the finitely many closed points $s \in S$.  
As an example, if $G$ has a Borel $S$-subgroup, then $G$ is totally isotropic (compare with the end of \S\ref{pp:parabolic}).
\epp


\bpp[Automorphisms of reductive groups] \label{pp:automorphisms}
For a scheme $S$ and a reductive $S$-group $G$, the functor 
that parametrizes group scheme automorphisms of (base changes of) $G$ is an extension
\be \label{eqn:Aut-extn}
1 \ra G^\ad \ra \underline{\Aut}_{\gp}(G) \ra \underline{\Out}_\gp(G) \ra  1
\ee
of  an $S$-group scheme $\underline{\Out}_\gp(G)$ that becomes constant and finitely generated \'{e}tale locally on $S$ by the adjoint group $G^\ad$ that parametrizes the inner automorphisms, see \cite{SGA3IIInew}*{Expos\'e XXIV, Th\'{e}or\`{e}me 1.3}. In particular, the $S$-group $\underline{\Aut}_{\gp}(G)$ is representable, smooth, and ind-quasi-affine (as reviewed in \S\ref{pp:conv}). Torsors under $\underline{\Aut}_{\gp}(G)$ correspond to forms of $G$, that is, to reductive $S$-groups that become  isomorphic to $G$ \'{e}tale locally on $S$. Such a form is inner (resp.,~pure inner) if the corresponding $\underline{\Aut}_{\gp}(G)$-torsor lifts to a $G^\ad$-torsor (resp.,~even to a $G$-torsor). Any $G$ splits \'{e}tale locally on $S$, so, at least for connected $S$, studying forms of $G$ amounts to studying forms of the corresponding split  $S$-group $\bbG$. For split groups, however, the extension \eqref{eqn:Aut-extn} admits a splitting, to the effect that 
\be \label{eqn:semidirect}
\underline{\Aut}_\gp(\bbG) \simeq \bbG^\ad \rtimes \underline{\Out}_\gp(\bbG);
\ee
concretely, a splitting $\underline{\Out}_\gp(\bbG) \le \underline{\Aut}_\gp(\bbG)$ is given by the subgroup of those automorphisms that preserve a fixed pinning of $\bbG$. 
\epp


\bpp[Isotriviality and embeddings into $\GL_n$] \label{pp:isotriviality}
As in the case of groups of multiplicative type reviewed in \S\ref{pp:def-isotrivial}, a reductive group scheme $G$ over a scheme $S$ is said to be \emph{isotrivial} if it becomes split over some finite \'{e}tale cover of $S$. If $S$ is affine and semilocal, then, by \cite{SGA3IIInew}*{Expos\'{e}~XXIV, Th\'{e}or\`{e}me 4.1.5, Corollaire 4.1.6} (see also \cite{Gil21}*{Corollary 4.4}), a reductive $S$-group $G$ of constant type is isotrivial if and only if its maximal central $S$-torus $\rad(G)$ is isotrivial, in which case every $G$-torsor is also \emph{isotrivial} in the sense that it trivializes over some finite \'{e}tale cover of $S$. We recall from \S\S\ref{pp:def-isotrivial}--\ref{pp:geom-uni-split} that $\rad(G)$ is isotrivial if it is $S$-fiberwise of rank $\le 1$, or if $S$ is Noetherian and its local rings are geometrically unibranch (for instance, normal), or if $S$ is local and normal. 

\addtocounter{footnote}{-5}
\renewcommand{\thefootnote}{\fnsymbol{footnote}}

By \cite{Gil21}*{Theorem 1.1} (which refines \cite{Tho87}*{Corollary 3.2}), a reductive group $G$ over a scheme $S$ is a closed subgroup of some $\GL(\sV)$ for a vector bundle $\sV$ on $S$ if and only if $\rad(G)$ is isotrivial.\footnote{\emph{Added after publication.} This sentence is only correct as stated for $G$ of constant type over $S$.} In the case when $S$ is affine, $\sV$ is a direct summand of a finite free $\sO_S$-module, and then one may even choose $\sV$ to be trivial, that is, if $S$ is affine and $\rad(G)$ is isotrivial, one may find a closed embedding
\[
G \hra \GL_{n,\,S} \qxq{for some} n \ge 1.
\]
Another class of $S$-groups $G$ that always admit a closed embedding $G \hra \GL(\sV)$ for some vector bundle $\sV$ on $S$ are the finite, locally free $S$-group schemes $G$: in this case, the translation action of $G$ on itself gives such an embedding by choosing $\sV$ to be the structure sheaf of $G$. Both for reductive and for finite, locally free $G$ and any embedding $G \hra \GL(\sV)$, the quotient $\GL(\sV)/G$ is affine, see \S\ref{pp:basic-representability} and \S\ref{pp:reductive} above.
\epp

\addtocounter{footnote}{4}

\bpp[Extending sections and torsors] \label{pp:extend}
When studying torsors, it is useful to keep in mind extension results that follow from general principles. To recall these, for a scheme $S$, a closed $Z \subset S$, and a $d \ge 0$, we write $\depth_Z(S) \ge d$ to mean that each $\fm_{S,\, z}$ with $z \in Z$ contains an $\sO_{S,\, z}$-regular sequence of length $d$. Then, by \cite{flat-purity}*{Lemma 7.2.7}, 
\be \label{eqn:unique-extn} \ba
E(S) \hra E(S \setminus Z) &\qx{for each separated $S$-scheme $E$, granted that $\depth_Z(S) \ge 1$;}\\
E(S) \isomto E(S \setminus Z) &\qx{for each affine $S$-scheme $E$ granted that $\depth_Z(S) \ge 2$;}
\ea\ee
(when $S$ is locally Noetherian, we may cite \cite{EGAI}*{Corollaire 9.5.6} and \cite{EGAIV2}*{Th\'{e}or\`{e}me~5.10.5}, respectively). In particular, if $\depth_Z(S) \ge 2$, then for any affine $S$-group scheme $G$ we have 
\be \label{eqn:H1-cheap-inj}
\qq H^1(S, G) \hra H^1(S \setminus Z, G),
\ee
in other words, nonisomorphic $G$-torsors over $S$ do not become isomorphic over $S \setminus Z$. 

The surjectivity of \eqref{eqn:H1-cheap-inj} is significantly more delicate, but it does hold if $S$ is regular (see \S\ref{sec:trichotomy}) of dimension $\le 2$, still with $\depth_Z(S) \ge 2$ (so $Z$ of codimension $\ge 2$), and $G$ is either reductive or finite flat: one first reduces to $G = \GL_n$ using embeddings as in \S\ref{pp:isotriviality} and then notes that, thanks to the Auslander--Buchsbaum formula, vector bundles over regular schemes of dimension $2$ extend uniquely over closed points, see \cite{CTS79}*{Corollary 6.14} for a detailed argument. 

Another situation in which \eqref{eqn:H1-cheap-inj} is surjective is if $S$ is regular with $\depth_Z(S) \ge 2$ and $G$ is of multiplicative type (see \S\ref{sec-multiplicative-type-group-schemes}), for instance, an $S$-torus: see \cite{CTS79}*{Corollary 6.9} and note that the key case of $\bG_m$-torsors follows by thinking of line bundles in terms of Weil divisors.
\epp


\csub[Basic properties of regular rings] \label{sec:trichotomy}


A local ring $(R, \fm)$ is \emph{regular} if it is Noetherian and
\[
\dim(R) = \dim_{R/\fm}(\fm/\fm^2),
\]
in which case the same holds for any localization of $R$ at a prime ideal, see \cite{SP}*{Lemma \href{https://stacks.math.columbia.edu/tag/0AFS}{0AFS}}. By Nakayama lemma, this equality amounts to requiring that the maximal ideal $\fm$ be generated by $\dim(R)$ elements.  A scheme is \emph{regular} if it is locally Noetherian and its local rings are regular. A ring is \emph{regular} if it is Noetherian and its localizations at prime (equivalently, maximal) ideals are regular. For a regular local ring $(R, \fm)$, a sequence $r_1, \dotsc, r_d \in \fm$ is a \emph{regular system of parameters} if the elements $r_i$ give a basis of the $R/\fm$-vector space $\fm/\fm^2$; an element $r \in \fm$ is a \emph{regular parameter} if it is a part of a regular system of parameters, in other words, if its image in $\fm/\fm^2$ is nonzero.


A regular ring of dimension $0$ is a product of fields. A regular ring of dimension $\le 1$ is a \emph{Dedekind ring}; each of its local rings is either a field or a discrete valuation ring. As a basic example, any smooth algebra over a Dedekind domain (such as a field or $\bZ$) 
or, more generally, over a regular ring, is regular. The definition of regularity is local, and regular local rings  split into the following classes.

\bd
A regular local ring $(R, \fm)$ of residue characteristic $p \ge 0$ is \emph{unramified} (resp.,~\emph{ramified}) if the ring $R/pR$ is regular (resp.,~is not regular), that is, if $p \in (\fm \setminus \fm^2)\cup \{ 0\}$ (resp.,~if~$p \in \fm^2 \setminus \{ 0\}$). A general regular ring (or a regular scheme) is \emph{unramified} if each of its local rings is unramified. 
\ed

If a regular ring $R$ is unramified, then so is every smooth $R$-algebra.

\brem
In addition to unramified regular local rings, it is useful to consider a larger class consisting of those regular local rings $R$ that are flat over some Dedekind ring $\cO$ and have geometrically regular $\cO$-fibers.\footnote{An algebra $C$ over a field $k$ is \emph{geometrically regular} if $C \tensor_k k'$ is a regular ring for every finite field extension $k'/k$.}  At the expense of more demanding technicalities (for instance, caused by imperfect residue fields of $\cO$), this class tends to be susceptible to the same techniques as the unramified case. The latter is recovered by restricting $\cO$ to be either $\bZ$, or $\bQ$, or some $\bF_p$.
\erem

The ramified regular local rings are necessarily of mixed characteristic (see \S\ref{pp:conv} for this terminology). In turn, the unramified ones split into two further classes: the regular local rings of equicharacteristic and the unramified regular local rings of mixed characteristic.




\beg
The following representative examples illustrate the classes of regular local rings:
\bitem
\m
of equicharacteristic (so also unramified): $k[x_1, \dotsc, x_d]_{(x_1, \dotsc,\, x_d)}$, where $k$ is a field;

\m
unramified of mixed characteristic: $\bZ[x_1, \dotsc, x_d]_{(p,\, x_1, \dotsc,\, x_d)}$, where $p$ is a prime;

\m
ramified: $(\bZ[x_1, \dotsc, x_d]/(p - x_1 \cdots x_r))_{(p,\, x_1, \dotsc,\, x_d)}$, where $1 < r \le d$.
\eitem
\eeg

In the spirit of these examples, the \emph{complete} regular local rings are classified as follows.

\blem[Cohen structure theorem] \lab{lem:Cohen}
\addtocounter{footnote}{-6}
\renewcommand{\thefootnote}{\fnsymbol{footnote}}
Any \emph{complete} regular local ring $(R, \fm)$ is of the form
\[
R \simeq W\llb x_1, \dotsc, x_d\rrb/(p - f) \qxq{with} f \in (p, x_1, \dotsc, x_d),
\]
where $W$ is a complete discrete valuation ring of mixed characteristic $(0, p)$ that has $p$ as a uniformizer.\footnote{\emph{Added after publication.} This sentence is correct as stated only if the residue field of $R$ has characteristic $p > 0$. In the remaining case when $R$ is of equicharacteristic $0$, we have $R \simeq k\llb x_1, \dotsc, x_d\rrb$ as in \ref{m:Cohen-i}.} 
Moreover, letting $k \ce R/\fm$ be the residue field, we have
\benumr
\m \label{m:Cohen-i}
$R$ is of equicharacteristic if and only if we may choose $f = 0$, so that 
\[
R \simeq k \llb x_1, \dotsc, x_d\rrb;
\]

\m
$R$ is unramified of mixed characteristic if and only if we may choose $f = p$, so that 
\[
R \simeq W \llb x_1, \dotsc, x_d\rrb;
\]

\m
$R$ is ramified if and only if we may choose $f \in (p, x_1, \dotsc, x_d)^2\, \setminus\, pR$. 

\eenum
\elem

\bpf
The claim follows from \cite{Mat89}*{Theorem 29.7} (the unramified case) and \cite{Mat89}*{Theorem~29.3 and the proof of Theorem 29.8~(ii)} (the ramified case). In general, the Cohen structure theorem applies beyond regular rings and describes the structure of complete Noetherian local rings, see \cite{EGAIV4}*{Chapitre 0, Th\'{e}or\`{e}me 19.8.8}.
\epf

\addtocounter{footnote}{5}

In principle, \Cref{lem:Cohen} exhaustively describes the structure of \emph{complete} regular local rings. Consequently, it tends to be important in those problems about regular rings that may be reduced to the complete local case. The problems about torsors are typically \emph{not} of this kind: for them, passage to completion may be no less difficult, so other structural results are needed. The central among such is the following highly useful theorem of Popescu that applies in the unramified case.


\bthm[Popescu] \label{thm:Popescu}
For a Noetherian ring $A$, a Noetherian $A$-algebra $B$ is a filtered direct limit of smooth $A$-algebras if and only if it is $A$-flat with geometrically regular $A$-fibers. In particular,
\benum
\m \label{P-a}
a regular ring $R$ that contains a perfect field $F$ \up{such as $\bQ$ or $\bF_p$} is a filtered direct limit of smooth $F$-algebras\uscolon

\m \label{P-c}
a regular ring $R$ that is a flat algebra over a Dedekind ring $\cO$ with geometrically regular $\cO$-fibers is a filtered direct limit of smooth $\cO$-algebras\uscolon

\m \label{P-b}
a regular local ring $R$ of mixed characteristic $(0, p)$ is unramified if and only if it is a filtered direct limit of smooth $\bZ_{(p)}$-algebras, equivalently, of smooth $\bZ$-algebras.
\eenum
\ethm

\bpf
The `only if' claim about $B$ is straight-forward, see \cite{SP}*{Lemma~\href{https://stacks.math.columbia.edu/tag/07DX}{07DX} or Lemma~\href{https://stacks.math.columbia.edu/tag/07EP}{07EP}}. In contrast, the `if' claim is intricate and has been the subject of surveys of its own: it is a result of Popescu \cite{Pop90}, whose proof has been clarified by Swan \cite{Swa98}, who, in turn, built on earlier clarifications due to Andr\'{e} and Ogoma; a modern account of the proof that includes some new simplifications is due to de Jong and is given in \cite{SP}*{Section \href{https://stacks.math.columbia.edu/tag/07BW}{07BW}, especially, Theorem \href{https://stacks.math.columbia.edu/tag/07GC}{07GC}}. 

Part \ref{P-a} (resp.,~\ref{P-c}; resp.,~\ref{P-b}) is a special case of the main claim with $B = R$ and $A = F$ (resp.,~$A = \cO$; resp.,~$A = \bZ_{(p)}$, equivalently, $A = \bZ$): indeed, to check that its conditions are met, we use that every finite field extension $k'$ of a perfect field $k$ is separable and that every \'{e}tale algebra over a regular ring is regular \cite{SP}*{Proposition \href{https://stacks.math.columbia.edu/tag/025N}{025N}}.
\epf

\Cref{thm:Popescu} is useful because smooth algebras are of finite type and may be studied using techniques from algebraic geometry, whereas the geometry of general regular rings is \emph{a priori} difficult to access directly. In effect, this link with algebraic geometry supplied by the Popescu theorem is one of the reasons why unramified regular local rings have been significantly more tractable in problems about torsors. In the ramified case, Popescu has recently established the following version of his theorem. 

\bthm[Popescu, \cite{Pop19}*{Theorem 3.8}] \label{thm:Popescu-2}
Every regular local ring is a filtered direct limit of regular local rings that are essentially of finite type as $\bZ$-algebras. \QED
\ethm

It seems plausible to us that this theorem, or perhaps some version or refinement thereof, could eventually be used for attacking the ramified case of problems about torsors over regular rings. At the moment, however, we are not aware of any application along such lines.


\section{Conjectures about vector bundles over regular rings} \label{sec:vector-bundles}

The most basic nonabelian reductive groups are the general linear groups $\GL_n$. For them, torsors amount to vector bundles, and we discuss the corresponding conjectures in this chapter. These conjectures (and much more) have already been discussed in the survey book \cite{Lam06}, but we hope that our summary would nevertheless be useful to some readers.


\csub[The Bass--Quillen conjecture] \label{sec:Bass-Quillen}

The Bass--Quillen conjecture is the flagship problem about vector bundles over regular rings. Posed in \cite{Bas73}*{Section 4.1} and \cite{Qui76}*{Comment (1) on page 170}, it grew out of Serre's problem solved by Quillen \cite{Qui76} and Suslin \cite{Sus76}: every vector bundle over an affine space over a field is free. Serre's problem and its variants were surveyed in a number of articles and books, notably in \cite{Lam06}.

\bconj[Bass--Quillen] \lab{conj:Bass-Quillen}
For a regular ring $R$, every vector bundle on $\bA^d_R$ descends to $R$.
\econj

\bpp[Basic reductions and known cases of \Cref{conj:Bass-Quillen}] \label{pp:basic-BQ} \hfill
\benuma
\item \label{BBQ-1}
The claim applies with $R[t_1, ..., t_{d - 1}]$ in place of $R$, so induction reduces one to the case~$d = 1$. 

\m \label{BBQ-2}
Once $d = 1$, Quillen patching of \Cref{cor:Quillen-patching} reduces one to local $R$. 
 In effect, it suffices to show that for a regular \emph{local} ring $R$, each finite projective $R[t]$-module is free. 

\item
The case when all the localizations of $R$ at its maximal ideals are unramified follows from results of Quillen, Suslin, Lindel, and Popescu, see \Cref{thm:BQ} below.

\item
Quillen and Suslin (independently) settled the case when $\dim(R) \le 1$, see \Cref{thm:Quillen-Suslin}~below.

\remi
The case when $\dim(R) \le 2$ and $d = 1$ is the Murthy--Horrocks theorem, see \cite{Lam06}*{Chapter~IV, Theorem 6.6}.

\remi
The case when $\dim(R) \le 3$ and $d = 1$ with $6 \in R^\times$ is due to Rao \cite{Rao88}*{Theorem 2}.

\item
It is incorrectly claimed in \cite{Lam06}*{page 330} and \cite{Pop17}*{Theorem 18~(3)} that the case when $R$ is local, Henselian, and excellent (for instance, complete) follows from  results of \cite{Pop89}.

\item
The case of line bundles is known: more generally, by \cite{Swa80}*{Theorem 1}, for any seminormal ring $A$, the map $\Pic(A) \ra \Pic(A[t_1, \ldots, t_d])$ is bijective  (in fact, by \emph{loc.~cit.},~these pullback maps on $\Pic(-)$ are  bijective if and only if $A_\red$ is seminormal).

\m
The analogue for valuation rings was established by Lequain and Simis \cite{LS80}: for a valuation ring $V$, every vector bundle on $\bA^d_V$ is free.
\eenum
\epp

To sum up, the main case in which the conjecture remains (widely) open is when $R$ is a ramified regular local ring. In \S\ref{sec:unramified-BQ} below, we review the proof of \Cref{conj:Bass-Quillen} in the unramified case.

A basic result that is used repeatedly in proving cases of the Bass--Quillen conjecture is the following theorem of Horrocks. We review it here because we will use it in the next section.

\bprop[Horrocks] \label{prop:Horrocks}
For a \up{resp.,~semilocal} ring $A$, a vector bundle on $\bA^1_A$ descends to $A$ \up{resp.,~is free on each connected component of $\bA^1_A$} iff it extends to a vector bundle on $\bP^1_A$. 
\eprop

\bpf
A finite projective module over a semilocal ring with connected spectrum is free (see \S\ref{pp:conv}), so the parenthetical assertion follows from the rest. For the latter, Quillen patching, that is, \Cref{cor:Quillen-patching}, reduces us to the case when $A$ is local, which is Horrocks theorem \cite{Hor64}*{Theorem~1}. Letting $k$ be the residue field of the maximal ideal of $A$, Horrocks used Grothendieck's classification of vector bundles on $\bP^1_k$ as direct sums of the $\sO(n)$ to first analyze the $k$-fiber of a vector bundle on $\bP^1_A$, and then applied the theorem on formal functions from \cite{EGAIII1} to bootstrap to $A$. See \cite{Lam06}*{Chapter IV, Section 2 onwards} for other proofs and also the discussion in \S\ref{sec:Guo} below.
\epf


\csub[The Quillen conjecture] \label{sec:Quillen}

In his resolution of Serre's problem about freeness of vector bundles on an affine space over a field, Quillen proposed an avenue of attack for the general case of the Bass--Quillen conjecture via the Horrocks theorem \ref{prop:Horrocks}. To apply the latter, for a regular local ring $A$ one needs to be able to extend a vector bundle $\sV$ on $\bA^1_{A}$ to $\bP^1_A$. The ring of the formal completion of $\bP^1_{A}$ along the infinity section is $A\llb t \rrb$, where $t$ is the inverse of the standard coordinate of $\bA^1_{A}$, so, by formal glueing (\Cref{lem:patch} below), $\sV$ extends to a vector bundle on $\bP^1_{A}$ if and only if it restricts to a free $A\llp t \rrp$-module. The following conjecture of Quillen \cite{Qui76}*{(2)~on~page 170} predicts that this is always the case.

\bconj[Quillen] \label{conj:Quillen}
For a regular local ring $(R, \fm)$ and an $r \in \fm$ that is a regular parameter \up{so $r \in \fm \setminus \fm^2$, see \uS\uref{sec:trichotomy}}, every finite projective $R[\f1r]$-module is free.
\econj

Equivalently, the conjecture predicts that every finite projective $R[\f1r]$-module extends to a finite projective $R$-module. Thus, formal glueing (see \Cref{lem:patch}) allows one to replace $R$ by its $r$-adic completion to assume without losing generality that $R$ is $r$-adically complete. As we explained above, \Cref{conj:Quillen} implies the Bass--Quillen \cref{conj:Bass-Quillen}.



\bpp[Basic reductions and known cases of \Cref{conj:Quillen}] \label{pp:basic-Quillen} \hfill
\benuma
\item \label{BQ-0}
The case $\dim(R) \le 2$ follows from general extension results reviewed in \S\ref{pp:extend}: by glueing, any finite projective $R[\f1r]$-module extends over the prime ideal $(r) \subset R$, and then the dimension assumption ensures that this extension extends further to a finite projective $R$-module. 

\item
The case when $R$ is of equicharacteristic was settled by Bhatwadekar--Rao \cite{BR83}*{Theorem~2.5} (see also \cite{Rao85}*{Theorem 2.9}) and Popescu (whose input amounts to \Cref{thm:Popescu}~\ref{P-a}). For a statement in mixed characteristic, see \cite{Tei95}*{page 272}.

\item \label{BQ-2}
The case $\dim( R) \le 3$ was settled by Gabber in \cite{Gab81}*{Chapter I, Theorem 1} and later in a simpler way by Swan in \cite{Swa88}. Very crudely,  the dimension assumption helps because any extension of a vector bundle on $R[\f1r]$ over the height one prime $(r) \subset R$ further extends uniquely to the entire ($2$-dimensional) punctured spectrum of $R$ (see \S\ref{pp:extend}), and the reduction of that extension modulo $r$ extends to all of $R/(r)$, so is free. With this premise, both proofs analyze all such extensions to the punctured spectrum of $R$ to show that one of them is free.


\eenum
To sum up, as in the Bass--Quillen conjecture \ref{conj:Bass-Quillen}, the main case in which the Quillen conjecture \ref{conj:Quillen} remains open is when the regular local ring $R$ is of mixed characteristic and, especially, ramified.
\epp

\bpp[Variants and generalizations] \label{pp:Quillen-variants} \hfill
\benuma
\item \label{QV-1}
Rao, following a suggestion of Nori, proposed the following generalization of \Cref{conj:Quillen} in \cite{Rao85}: for a regular local ring $R$ and $r_1, \dotsc, r_t \in \fm$ that form a part of a regular system of parameters, every finite projective $R[\f{1}{r_1\cdots r_t}]$-module is free. Strictly speaking, \emph{loc.~cit.}~only considered those $R$ that are localizations of finite type, regular algebras over some infinite field and, with this restriction, established the generalization when either $t \le 2$ or $\dim(R) \le 5$, see \cite{Rao85}*{Corollaries 2.10 and 2.11}. See also \cite{Gab02}*{Theorem 1.1} (possibly also the earlier \cite{Nis98}) for further results on this generalization of \Cref{conj:Quillen}. 


\item
The assumption that $r \not \in \fm^2$ is critical in \Cref{conj:Quillen}. For instance, as is pointed out in \cite{BR83}*{Example (1) on page 808}, when
\[
\qq R \ce \bR[x, y, z]_{(x,\, y,\, z)} \qxq{and} r \ce x^2 + y^2 + z^2,
\] 
the kernel of the surjection $R[\f1r]^{\oplus 3} \surjects R[\f1r]$ given by $(a, b, c) \mapsto ax + by + cz$ is a nonfree projective $R[\f{1}{r}]$-module of rank $2$, see \cite{Lam06}*{pages 34--35}.

\item
The global version of \Cref{conj:Quillen}, posed in \cite{Qui76}*{(3) on page 170} as a question, is false: for an affine, regular scheme $X$ and a regular divisor $Z \subset X$, a vector bundle on $X \setminus Z$ need not extend to $X$. In \cite{Swa78}*{Section 2}, Swan showed that this happens  for 
\[
\qq X \ce \Spec\p{\bC[x_0, \dotsc, x_4, t]/(x_0^2 + \dotsc + x_4^2 - 1)} \qxq{and} Z \ce \{ t = 0\},
\]
namely, there is a (stably free) vector bundle of rank $2$ on $X \setminus Z$ that does not extend to $X$.
\eenum
\epp




\csub[The Lam conjecture] \label{sec:Lam}

Regularity of $R$ is a critical assumption in the Bass--Quillen conjecture \ref{conj:Bass-Quillen}, see \cite{Lam06}*{page 342} for an overview of counterexamples when this assumption is weakened. In these counterexamples, even the map $K_0(R) \ra K_0(R[t_1, \dotsc, t_d])$ between the Grothendieck rings of vector bundles\footnote{We recall that the Grothendieck ring $K_0(A)$ of a commutative ring $A$ is the quotient of the free abelian group on the set of isomorphism classes of finite projective $A$-modules $P$ by the relations $[P] = [P'] + [P'']$ for finite projective $A$-modules $P$, $P'$, $P''$ with $P \simeq P' \oplus P''$, and that the multiplication in $K_0(A)$ is induced by the tensor product~$\tensor_A $.} fails to be an isomorphism, in other words, the ``obstruction'' is visible already on $K_0$. In contrast, by Grothendieck's theorem \cite{Lam06}*{Chapter II, Theorem 5.8}, for any \emph{regular} ring $R$, we have
\[
K_0(R) \isomto K_0(R[t_1, \dotsc, t_d]).
\]
On the other hand, for a ring $A$, the stably free (see \S\ref{pp:conv}) finite projective $A$-modules $P$  are precisely those whose classes lie in $\bZ \subset K_0(A)$, see \cite{Lam06}*{Chapter I, Corollary 6.2}. Thus, it is conceivable that one could obtain a Bass--Quillen conjecture for arbitrary rings by only considering stably free modules. This is precisely what the following conjecture posed by Lam \cite{Lam06}*{page 180} predicts.

\bconj[Lam] \lab{conj:Lam}
For a ring $A$, every stably free vector bundle on $\bA^d_A$ descends to $A$. Equivalently, for a \emph{local} ring $A$, every stably free $A[t]$-module is free.
\econj

The equivalence of the two formulations, that is, the reduction to local $A$, follows from Quillen patching of \Cref{cor:Quillen-patching}. By the discussion above, for a \emph{regular} local $A$, every finite projective $A[t_1, \dotsc, t_d]$-module is stably free, so \Cref{conj:Lam} implies the Bass--Quillen conjecture \ref{conj:Bass-Quillen}.

Although we attribute the conjecture to Lam, Swan raised it as a question already in \cite{Swa78}*{page~114,~(B)}. We stress that even though \Cref{conj:Lam} offers the advantage of making no assumption on $A$, it has to be regarded as very speculative. Indeed, it is only known in very few cases: the ones in \S\ref{pp:basic-BQ} and, by \cite{Yen08}*{Corollary 5}, also in the case when $\dim (A) \le 1$ and $d = 1$.

The following proposition offers several equivalent versions of \Cref{conj:Lam} that may be useful to keep in mind while contemplating possible arguments or counterexamples. 

\bprop \lab{Lam-conj}
Let $A$ be a local ring and let $A(t)$ be the localization of $A[t]$ with respect to all the monic polynomials. The following statements are equivalent\ucolon
\benum
\m \label{LC-a} 
every stably free $A[t]$-module is free\uscolon

\m \label{LC-b}
every stably free $A(t)$-module is free\uscolon

\m \label{LC-c}
every stably free $A\llp t \rrp$-module is free.

\eenum
\eprop

\bpf
The equivalence between \ref{LC-a} and \ref{LC-b} was established by Bhatwadekar and Rao in \cite{BR83}*{Theorem A}. The equivalence between \ref{LC-b} and \ref{LC-c} was established in \cite{Hitchin-torsors}*{Theorem 2.1.25~(c)--(d)}: roughly, one uses the $A$-isomorphism $A(t)\simeq (A[t]_{1 + tA[t]})[\f1t]$ (induced by ``$t \mapsto t\i$''), notes that the completion of $(A[t]_{1 + tA[t]})[\f1t]$ for the ``$t$-adic'' topology is $A\llp t \rrp$, and, crucially, shows that the functor of isomorphism classes of stably free modules is  invariant under such completion because it is invariant under Zariski pairs. In this last step, Zariski pairs come about via Gabber's technique of considering the ring of ``$t$-adic'' Cauchy sequences valued in $(A[t]_{1 + tA[t]})[\f1t]$: this ring is Zariski along its ideal formed by nil sequences, and the corresponding quotient is the completion $A\llp t \rrp$. 
\epf

Beyond local $A$ and stably free modules, one has the following result of a similar spirit.

\bprop \label{prop:GLn-noniso}
For a ring $A$, with $A(t)$ as in Proposition \uref{Lam-conj}, nonisomorphic finite projective $A$-modules \up{resp.,~$A(t)$-modules} cannot become isomorphic after base change to $A(t)$~or~to~$A\llp t \rrp$.
\eprop

\bpf
The claim is contained in \cite{Hitchin-torsors}*{Theorem 2.1.25 (a)--(b)}. The first main input is the Cauchy sequence technique mentioned in the proof of \Cref{Lam-conj}: it achieves a comparison between finite projective modules over $A(t)$ and $A\llp t\rrp$. The second main input is \cite{Lam06}*{Chapter V, Proposition~2.4}, which is an elementary patching argument due to Bass that shows that nonisomorphic finite projective $A$-modules $P$ and $P'$ cannot become isomorphic after base change to $A(t)$: in more detail, an $A(t)$-isomorphism $P_{A(t)} \simeq P'_{A(t)}$ would permit us to use base changes of $P$ and $P'$ to glue up a vector bundle over $\bP^1_A$ whose restriction along the section $\{t = 1\}$ (resp.,~$\{t = \infty\}$) is $P$ (resp.,~$P'$), and the Horrocks \Cref{prop:Horrocks} would then give the desired $P \simeq P'$.
\epf





\section{Conjectures about torsors under reductive groups over regular rings} \label{sec:general-groups}

We turn our attention to torsors under general reductive group schemes over regular bases. Questions about them tend to be more subtle than the ones about vector bundles discussed in \S\ref{sec:vector-bundles} because group-theoretic properties start playing important roles in the arguments. However, these questions about general reductive groups may be more susceptible to progress, perhaps simply for the reason that some of them do not appear to have been studied as extensively.


\csub[The Grothendieck--Serre conjecture] \label{sec:Grothendieck-Serre}

The following conjecture of Grothendieck and Serre is the flagship problem about torsors under reductive group schemes over regular rings. It originated from its special cases conjectured by Serre \cite{Ser58b}*{Remark on page 31}~and Grothendieck \cite{Gro58}*{Remark on pages 26--27}, \cite{Gro68b}*{Remark~1.11~a)} and was popularized by the article of Colliot-Th\'{e}l\`{e}ne and Ojanguren \cite{CTO92}. It was also the subject of a recent ICM survey of Panin \cite{Pan18}, which we refer to for further discussion.

\bconj[Grothendieck--Serre] \label{conj:Grothendieck-Serre}
For a regular local ring $R$ and a reductive $R$-group scheme $G$, no nontrivial $G$-torsor over $R$ trivializes over $K \ce \Frac(R)$, equivalently,
\be \label{eqn:GS}
H^1(R, G) \hra H^1(K, G).
\ee
\econj

The claimed equivalence of the two statements follows from the twisting bijections \eqref{eqn:change-origin}: more precisely, the injectivity of \eqref{eqn:GS} for $G$ is equivalent to no nontrivial torsor over $R$ trivializing over $K$ \emph{for all} pure inner forms of $G$. 

\beg
For  a reductive group $G$ over an algebraically closed field $k$ (such as $\bC$), the conjecture predicts that every generically trivial $G$-torsor over a smooth algebraic variety $X$ over $k$ is Zariski locally trivial. This was Serre's original formulation and was settled by Colliot-Th\'{e}l\`{e}ne--Ojanguren in \cite{CTO92}, see also \S\ref{pp:basic-GS}~\ref{pp:BGS-3}. By Steinberg theorem \cite{Ser02}*{Chapter~III, Section 2.3, Theorem 1$'$ and Remarks 1) (with Chapter II, Section 3.3, b))}, generic triviality is automatic if $X$ is a curve and, by de Jong--He--Starr theorem \cite{HdJS11}*{Theorem 1.4}, also if both $X$ is a surface and $G$ is semisimple, simply connected; thus, in these cases every $G$-torsor over $X$ is Zariski locally trivial.
\eeg

\beg \label{eg:Brauer-square}
In the cases $G = \GL_n$ or $G = \SL_n$, both the source and the target of \eqref{eqn:GS} vanish, so it is more instructive to consider the case $G = \PGL_n$, in which the conjecture predicts that an Azumaya algebra over $R$ that is isomorphic to a matrix algebra over $K$ is isomorphic to a matrix algebra already over $R$. The central extension 
\[
1 \ra \bG_m \ra \GL_n \ra \PGL_n \ra 1
\]
and its associated long exact sequence \eqref{eqn:coho-seq} give a commutative square
\[
\xymatrix{
H^1(R, \PGL_n) \ar[d] \ar[r] & H^2(R, \bG_m) \ar[d]  \\
H^1(K, \PGL_n) \ar[r] & H^2(K, \bG_m)
}
\]
whose horizontal maps have trivial kernels. Since $R$ is regular, Grothendieck's injectivity result for the Brauer group \cite{Gro68b}*{Corollaire 1.8} implies that the right vertical map is injective. Thus, the left vertical map has trivial kernel, so \Cref{conj:Grothendieck-Serre} holds in the case when $G = \PGL_n$. In fact, this case seems to have been one of the main motivations for the conjecture, see \cite{Gro68b}*{Remarques~1.11~a)}.
\eeg

For relations between the Grothendieck--Serre conjecture \ref{conj:Grothendieck-Serre} and certain group decompositions, see \cite{split-unramified}*{Corollary 1.3}. For consequences for quadratic forms, see \cite{split-unramified}*{Corollary 9.6}

\bpp[Basic reductions and known cases of \Cref{conj:Grothendieck-Serre}] \label{pp:basic-GS}
The known cases of \Cref{conj:Grothendieck-Serre} have already been summarized in \cite{split-unramified}*{Section 1.4} or \cite{Pan18}, and the literature is vast, so here we do not attempt to be exhaustive and focus on overviewing the main known cases.
\benuma
\m \label{pp:BGS-1}
The case when $G$ is a torus was settled by Colliot-Th\'{e}l\`{e}ne and Sansuc in \cite{CTS78} and \cite{CTS87}. The latter reference uses flasque resolutions of tori reviewed in \Cref{sec-surjection-quasi-trivial} below to reduce to ``simpler'' tori. The toral case is used often in arguing other cases of \Cref{conj:Grothendieck-Serre}.

\item \label{pp:BGS-2}
The case when $\dim( R)  \le 1$, that is, when $R$ is either a field (trivial case) or a discrete valuation ring was settled by Nisnevich in \cite{Nis82}, \cite{Nis84}, with clarifications and complements given by Guo in \cite{Guo20}. Roughly, the idea is to replace $R$ by its completion via approximation arguments that go back to Harder and to then exploit the Bruhat--Tits theory to conclude. Guo's result also gives the semilocal case: the statement of \Cref{conj:Grothendieck-Serre} holds when $R$ is a semilocal Dedekind domain. One question that seems to still be open in the discrete valuation ring case is whether \eqref{eqn:GS} remains injective when $G$ is merely a parahoric $R$-group scheme. 

By induction on $\dim( R)$, the $1$-dimensional case implies that any generically trivial torsor over a regular local ring $R$ trivializes over the residue field of any prime $\fp \subset R$. In particular, it trivializes over the residue field of the maximal ideal of $R$, and Hensel's lemma \cite{EGAIV4}*{Th\'{e}or\`{e}me~18.5.17} then implies that \Cref{conj:Grothendieck-Serre} holds in the case when the regular local ring $R$ is Henselian, for instance, complete (compare with \cite{CTS79}*{Assertion~6.6.1}).


Thus, one possible point of view is that for a general $R$ the main difficulty lies in passing to the completion. Such passage remains out of reach when $\dim( R) \ge 2$: in effect, in this higher-dimensional case, geometric approaches have so far been more fruitful.

\m \label{pp:BGS-3}
The case when $R$ is of equicharacteristic, that is, when $R$ contains a field, was settled by Fedorov--Panin \cite{FP15} when the field is infinite and by Panin \cite{Pan20a} when the field is finite. In spite of numerous group-theoretic subtleties that accompany an arbitrary $G$, crudely speaking, the overall structure of the Fedorov--Panin strategy is somewhat similar to the approach to the Bass--Quillen conjecture discussed in \S\ref{sec:Bass-Quillen} and \S\ref{sec:unramified-BQ} below: more precisely, it uses the Popescu Theorem \ref{thm:Popescu} to pass to local rings of smooth algebras over a field, it then combines Artin's results on good neighborhoods from \cite{SGA4III}*{Expos\'{e} XI} with Voevodsky's ideas that appear through Panin's notion of ``nice triples'' (which are smooth relative curves over $R$ equipped with a section and an $R$-finite closed subscheme) to pass via excision to studying torsors over the relative affine line $\bA^1_R$, and it concludes via Horrocks-style results aided by insights from the geometry of the affine Grassmannian. The split of the argument into the cases of an infinite versus a finite base field was primarily due to technical difficulties caused by the Bertini theorem over finite fields (these difficulties have since been resolved).

In their strategy, one first reduced to semisimple, simply-connected $G$, although the need for this has since been eliminated by Fedorov by refining the part that concerns the affine Grassmannian, see \cite{Fed21a} and \S\ref{sec:A1-analysis} below. In the general case of \Cref{conj:Grothendieck-Serre}, a reduction to semisimple, simply connected $G$ remains unavailable; however, one may at least reduce to those $G$ whose derived group $G^\der$ is simply connected, see \Cref{prop:GS-app} below.

\m
In mixed characteristic, the case when $R$ is unramified and the group $G$ has a Borel $R$-subgroup was settled in \cite{split-unramified}.  The argument builds on the Panin--Fedorov strategy, in fact, it simultaneously reproves the equal characteristic case. The main novelties in comparison to their strategy are in the ``middle part'' of the argument: the role of Artin's good neighborhoods got replaced by a presentation lemma in the style of Gabber (see \Cref{thm:geometric-presentation}), and ``nice triples'' were replaced by a more direct analysis of relative curves (see \Cref{prop:embed}); this simplified the argument to the point that, modulo circumventing some technical difficulties caused by mixed characteristic, it could work over discrete valuation rings in place of fields. We refer to \cite{split-unramified}*{especially, Section 1.6} for more details.

For the moment, the unramifiedness assumption seems difficult to bypass in any ``geometric'' approach that eventually reduces to the relative affine line $\bA^1_R$ (see the end of \S\ref{sec:trichotomy}). In turn, the Borel $R$-subgroup helps by ensuring that a generically trivial $G$-torsor $E$ over $R$ reduces to a $B$-torsor away from some closed subset $Z \subset \Spec(R)$ of codimension $\ge 2$ (apply the valuative criterion of properness to $E/B$), and this codimension aspect is used in a crucial way for extending the presentation lemma to mixed characteristic. Other ways in which a Borel helps are that it allows one to reduce to the semisimple, simply connected case, even without $R$ being unramified, and that it allows one to bypass the compactification question discussed in \S\ref{sec:conj-compactify} below (but both of these seem less essential at the cost of further work).

\m
Beyond the cases above, some sporadic cases were settled in \cite{Oja82}, \cite{Nis89}, \cite{Fir22}, \cite{BFFP20}.

\m
The analogue of the Grothendieck--Serre conjecture for valuation rings was established by Guo \cite{Guo21}: for a valuation ring $V$ and a reductive $V$-group $G$, no nontrivial $G$-torsor trivializes over $\Frac(V)$. A desirable further step in this direction would be to show the same with $V$ replaced by a local ring of a smooth scheme over a valuation ring (say, of equicharacteristic). 
\eenum

To sum up, the Grothendieck--Serre conjecture is known in equal characteristic but remains open in mixed characteristic, especially, over ramified regular local rings, for which one may need a substantially different approach. As for the unramified mixed characteristic case beyond quasi-split $G$, we feel that it may, in principle, be approachable, perhaps by finding some way to improve or to bypass the presentation lemma, but we do not know of a precise way to attack it fruitfully.
\epp


The Grothendieck--Serre conjecture has the following consequence for uniqueness of reductive group schemes with a fixed generic fiber over a regular local ring.

\bprop \label{prop:unique-reductive-model}
For a regular local ring $R$, its fraction field $K$, and a reductive $R$-group scheme $G$ such that the Grothendieck--Serre conjecture \uref{conj:Grothendieck-Serre} holds for every form of $G^\ad$, up to isomorphism $G$ is the unique reductive $R$-group scheme with the generic fiber isomorphic to $G_K$. 
\eprop

In particular, for a regular local ring $R$ that is either of equicharacteristic or of dimension $\le 1$,  nonisomorphic reductive $R$-group schemes do not become isomorphic over $K$ (see \S\ref{pp:basic-GS}~\ref{pp:BGS-2}--\ref{pp:BGS-3}).

\bpf
One uses the extension structure \eqref{eqn:Aut-extn} to argue that the map
\[
H^1(R, \underline{\Aut}_\gp(G)) \ra H^1(K, \underline{\Aut}_\gp(G))
\]
has a trivial kernel, see \cite{Guo20}*{Proposition 14} for a detailed argument.
\epf

\brem
For unramified $R$, the known cases of \Cref{conj:Grothendieck-Serre} suffice for showing that a reductive $R$-group scheme $G$ is split if and only if so is its generic fiber $G_K$, see \cite{split-unramified}*{Theorem~9.3}.
\erem

\brem \label{rem:abelian-schemes}
In the context of the Grothendieck--Serre conjecture, one may consider the analogy with an abelian scheme $A$ over regular base scheme $S$. A key simplifying difference is that in this case, for any dense open $U \subset S$, one has
$A(S) \isomto A(U)$ (the injectivity follows from the separatedness of $A$ (see \S\ref{pp:extend}) and the surjectivity follows by considering $A$ as its own double dual and by extending line bundles via the regularity assumption, see \cite{BLR90}*{Section 8.4, Corollary 6}). By applying this \'{e}tale locally on $S$, the same  holds for a torsor under an abelian scheme, so that, in particular,
\[
H^1(S, A) \hra H^1(U, A).
\]
Consequently, the Grothendieck--Serre conjecture \ref{conj:Grothendieck-Serre} holds if $G$ is replaced by an abelian $R$-scheme. Similarly, \Cref{prop:unique-reductive-model} holds if $G$ is replaced by an abelian scheme, see \cite{Fal83}*{Lemma 1} (or recall that the moduli scheme of suitably polarized abelian schemes with level structure is separated).
\erem

\csub[The Colliot-Th\'{e}l\`{e}ne--Sansuc purity conjecture] \label{sec:CTS}

The Grothendieck--Serre conjecture \ref{conj:Grothendieck-Serre} predicts that, in its notation, $H^1(R, G)$ is a subset of $H^1(K, G)$. The following purity conjecture posed as \cite{CTS79}*{Question~6.4} characterizes this subset.

\bconj[Colliot-Th\'{e}l\`{e}ne--Sansuc]\label{conj:CTS}
For a regular local ring $R$, its fraction field $K$, and a reductive $R$-group scheme $G$, a $G$-torsor over $K$ that extends to a $G$-torsor over $R_\fp$ for every height $1$ prime $\fp \subset R$ extends uniquely to a $G$-torsor. In other words, we have
\[
H^1(R, G) = \bigcap_{\fp\x{ \upshape{of height} }1} H^1(R_\fp, G) \qxq{inside} H^1(K, G).
\]
\econj


It is instructive to contrast this conjecture with the purity for the Brauer group, according to which
\[
\Br(R) = \bigcap_{\fp\x{ of height }1} \Br(R_\fp) \qxq{inside} \Br(K),
\]
see \cite{brauer-purity}*{Theorem 6.2}. This Brauer group variant had been conjectured by Auslander--Goldman in \cite{AG60}, established by Gabber in most cases, and completed in the remaining cases of  mixed characteristic in \cite{brauer-purity} using a perfectoid method, see \emph{op.~cit.}~for an overview of prior literature.  The principal reason why the Brauer group version is more approachable is the relation to derived functor cohomology, namely, to the cohomological Brauer group via the isomorphism 
\[
\Br(R) \cong H^2_\et(R, \bG_m)_\tors
\]
due to Gabber \cite{Gab81}*{Chapter II, Theorem 1}. This dramatically broadens the range of available techniques, basically, because abelian cohomology classes are simpler to manipulate than torsors.

In \Cref{conj:CTS}, by spreading out and glueing in the Zariski topology, any $G$-torsor over $K$ that extends to a $G$-torsor over $R_\fp$ for every height $1$ prime $\fp \subset R$ also extends to a $G$-torsor $E$ over $U$ for some nonempty open $U \subset \Spec(R)$ whose complement is of codimension $\ge 2$. However, there may be many ways to glue, and so many possible $E$ with the same generic fiber $E_{K}$---the key point, and the main difficulty, is to be able to glue in such a way that $E$ extends to a $G$-torsor over $\Spec(R)$ (in the Brauer group case, the analogue of $E$ is automatically unique for any fixed $U$). In contrast, thanks to \S\ref{pp:extend}, this further extension to a $G$-torsor will have to be unique.

\bpp[Basic reductions and known cases of \Cref{conj:CTS}] \label{pp:basic-CTS} \hfill
\benuma
\item \label{BCTS-1}
The case when either $\dim( R) \le 2$ or $G$ is a torus follow from general principles reviewed in \S\ref{pp:extend}: these assumptions ensure that any $E$ as above extends to a $G$-torsor over $\Spec(R)$. 



\m
Some cases with split $G$ were settled by Chernousov--Panin, Panin, Panin--Pimenov, and Antieau--Williams in \cite{CP07}, \cite{Pan10}, \cite{PP10}, \cite{CP13}, and \cite{AW15}, see \cite{EKW21}*{Remark 4.3}.

\m
Antieau and Williams showed in \cite{AW15} that, once the regular domain $R$ is no longer assumed to be local, the statement of \Cref{conj:CTS} does not hold even in the case $G = \PGL_n$.
\eenum
\epp

To sum up, the Colliot-Th\'{e}l\`{e}ne--Sansuc purity \cref{conj:CTS} remains widely open beyond somewhat restrictive special cases. In the rest of this section, we turn to  
its following consequence for extending reductive group schemes. This theme has also been investigated by Vasiu from the point of view of extending their associated Lie algebras, see \cite{Vas16} for details.

\bconj \label{conj:purity-model}
For a regular local ring $R$ and its fraction field $K$, a reductive $K$-group extends to a reductive $R$-group scheme if and only if it extends to a reductive $R_\fp$-group scheme for every prime $\fp \subset R$ of height $1$, in which case this extension is unique up to isomorphism.
\econj

\bpp[\Cref{conj:CTS} implies \Cref{conj:purity-model}] \label{pp:first-implies}
We assume that the Colliot-Th\'{e}l\`{e}ne--Sansuc purity conjecture \ref{conj:CTS} holds for adjoint semisimple $R$-group schemes, and we will argue that then every reductive $K$-group $G$ that extends to a reductive $R_\fp$-group scheme for every prime $\fp \subset R$ of height $1$ also extends to a reductive $R$-group scheme. By spreading out and glueing, we may take advantage of the assumption on the $R_\fp$ to arrange that $G$ begins life as a reductive $U$-group scheme for some open $U \subset \Spec(R)$ whose complement is of codimension $\ge 2$. Letting $\bbG$ be the split reductive $R$-group scheme of the same type as $G$, we then use the dictionary of \S\ref{pp:automorphisms} to reduce to showing that for every $\underline{\Aut}_{\gp}(\bbG)$-torsor $E$ over $U$ there is an $\underline{\Aut}_{\gp}(\bbG)$-torsor $\cE$ over $R$ such that $\cE_K \simeq E_K$ as $\underline{\Aut}_{\gp}(\bbG)$-torsors over $K$. For this, we first show the same for $\underline{\Out}_\gp(\bbG)$-torsors. 

Let $F$ be an $\underline{\Out}_\gp(\bbG)$-torsor over $U$. By \cite{SGA3II}*{Expos\'e X, Corollaire 5.14}, the connected components of $F$ are open and finite \'{e}tale over $U$. Thus, Auslander--Nagata purity \cite{SGA2new}*{Expos\'e~X, Th\'{e}or\`{e}me 3.4} ensures that they extend uniquely to finite \'{e}tale $R$-schemes. In this way, $F$ extends to an \'{e}tale locally constant $R$-scheme $\cF$, and the maps describing the torsor structure likewise extend and make $\cF$ an $\underline{\Out}_\gp(\bbG)$-torsor over $R$. 

We now let $F$ to be the $\underline{\Out}_\gp(\bbG)$-torsor induced by the  $\underline{\Aut}_{\gp}(\bbG)$-torsor $E$ via \eqref{eqn:Aut-extn}. We use the splitting \eqref{eqn:semidirect} to view the resulting $\underline{\Out}_\gp(\bbG)$-torsor $\cF$ as an $\underline{\Aut}_{\gp}(\bbG)$-torsor. We twist by this $\underline{\Aut}_{\gp}(\bbG)$-torsor $\cF$ and combine the twisting bijection \eqref{eqn:change-origin} with the cohomology exact sequence \eqref{eqn:coho-seq} to note that the image of $E$ under this bijection comes from a torsor $E'$ over $U$ under an $R$-form $\cG$ of $\bbG^\ad$. By the assumed \Cref{conj:CTS} for $\cG$, there is a $\cG$-torsor $\cE'$ over $R$ with $\cE'_K \simeq E'_K$. This $\cE'$ gives rise to a torsor under the twist of $\underline{\Aut}_{\gp}(\bbG)$ in question, and its preimage under the twisting bijection is a desired $\underline{\Aut}_{\gp}(\bbG)$-torsor $\cE$ with $\cE_K \simeq E_K$. \QED
\epp



\csub[The Grothendieck--Serre conjecture for Levi reductions and parabolic subgroups]

We wish to draw attention to two variants of the Grothendieck--Serre conjecture, one for Levi reductions of torsors and another one for parabolic subgroups of reductive groups. As we show, the parabolic variant implies the Levi variant and both follow from the conjectures discussed in \S\S\ref{sec:Grothendieck-Serre}--\ref{sec:CTS}. One could hope that these variants may be more amenable to direct attack. We begin with the less general variant that concerns Levi reductions of torsors under reductive groups.


\bconj \label{conj:Levi}
For a regular local ring $R$, its fraction field $K$, a reductive $R$-group scheme $G$, a parabolic $R$-subgroup $P \subset G$, and an $R$-Levi $M \subset P$, a $G$-torsor $E$ reduces \up{necessarily uniquely} to an $M$-torsor \up{equivalently, a $P$-torsor} iff $E_{K}$ reduces to an $M_{K}$-torsor \up{equivalently, a $P_{K}$-torsor}. 
\econj

The parenthetical aspects follow from the basic review of \uS\uref{pp:parabolic-torsors}, especially, from \eqref{eqn:M-P}--\eqref{eqn:Levi-H1}.

\bpp[\Cref{conj:Grothendieck-Serre,conj:CTS} imply \Cref{conj:Levi}] \label{pp:second-implies}
We assume that the Grothendieck--Serre map \eqref{eqn:GS} is injective for $G$ and that the Colliot-Th\'{e}l\`{e}ne--Sansuc conjecture \ref{conj:CTS} holds for $M$. Since the `only if' is obvious, 
we seek the converse, so we let $F$ be the unique $P_{K}$-torsor that induces $E_{K}$. By the Grothendieck--Serre conjecture for $G$, namely, by \eqref{eqn:GS}, we need to show that $F$ extends to an $P$-torsor over $R$. 
For this, \eqref{eqn:M-P} and the Colliot-Th\'{e}l\`{e}ne--Sansuc conjecture~\ref{conj:CTS} for $M$ reduce us to showing that $F$ extends to a $P$-torsor over $R_\fp$ for every height $1$ prime $\fp \subset R$, so we may assume that $R$ is a discrete valuation ring. But then the valuative criterion of properness applied to $E/P$ suffices: $F$ amounts to a $K$-point of $E/P$ (see see \S\ref{pp:basic-subgroups}), which automatically extends to an $R$-point of $E/P$, which amounts to the desired extension of $F$ to a $P$-torsor. \QED
\epp

The following variant for parabolic subgroups is more general, goes back to ideas of Colliot-Th\'{e}l\`{e}ne and Panin, and was stated as a conjecture in \cite{split-unramified}*{Conjecture 9.4}. 

\bconj[Colliot-Th\'{e}l\`{e}ne--Panin] \label{conj:CTP}
For a regular local ring $R$, its fraction field $K$, and a reductive $R$-group scheme $G$, if $G_K$ has a proper parabolic subgroup, then so does $G$\uscolon more precisely, if $G_K$ has a parabolic $K$-subgroup of a fixed type, then $G$ has a parabolic $R$-subgroup of the same~type.
\econj

The $R$-scheme that parametrizes the types of parabolic subgroups of (base changes of) $G$ is finite \'{e}tale (see \cite{SGA3IIInew}*{Expos\'e XXVI, Section 3.1, D\'{e}finition 3.4}), so its $K$-points extend uniquely to $R$-points. Thus, the aspect of the conjecture about parabolics of the same type is well posed.

\bpp[\Cref{conj:CTP} implies \Cref{conj:Levi}]
Let $R$, $G$, and $P$ be as in \Cref{conj:Levi}, and let $E$ be a $G$-torsor. By \S\ref{pp:basic-subgroups} and \S\ref{pp:parabolic}, the quotient $E/P$ parametrizes both reductions of $E$ to a $P$-torsor and also parabolic subgroups of the same type as $P$ of the inner form $\Aut_G(E)$ of $G$. Thus, if $E_K$ reduces to a $P_K$-torsor, then $\Aut_G(E)_K$ has a parabolic subgroup of the same type as $P$. By \Cref{conj:CTP}, then $\Aut_G(E)$ itself has a parabolic subgroup of the same type as $P$, so that $E$ reduces to a $P$-torsor over $R$, as predicted by \Cref{conj:Levi}. \QED
\epp

\bpp[Basic reductions and known cases of \Cref{conj:CTP}] \label{pp:known-CTP} \hfill
\benuma
\m \label{pp:BCTP-1}
The case when $\dim( R) \le 1$ follows from the valuative criterion of properness: in fact, since the scheme that parametrizes parabolic subgroups of $G$ of a fixed type is proper (see \S\ref{pp:parabolic}), every parabolic $K$-subgroup of $G_K$ extends to a parabolic $R$-subgroup of $G$. Similarly, for a general regular local $R$, any parabolic $K$-subgroup of $G_K$ extends to a parabolic $U$-subgroup of $G_U$ for some open $U \subset \Spec(R)$ whose complement is of codimension $\ge 2$. The difficulty lies in arguing the existence of a $K$-parabolic of $G_K$ for which even $U = \Spec(R)$.

\m
As observed by Sean Cotner, the case when $\dim( R) \le 1$ implies the case when $R$ is Henselian. Indeed, analogously to \S\ref{pp:basic-GS}~\ref{pp:BGS-2}, induction on $\dim( R)$ shows that the desired parabolic subgroup exists over the residue field of any prime $\fp \subset R$; since the scheme that parametrizes parabolic subgroups of $G$ of a fixed type is smooth (see \S\ref{pp:parabolic}), Hensel's lemma \cite{EGAIV4}*{Th\'{e}or\`{e}me~18.5.17} then lifts the parabolic from the residue field $k$ to all of $R$.

\m
Several cases in which $G$ is an orthogonal group were settled in \cite{CT79}, \cite{Pan09}, \cite{PP10}, \cite{PP15}, \cite{Scu18}. In fact, these cases related to quadratic forms suggested the general conjecture.

\m
The case when $G_k$, where $k$ is the residue field of $R$, has no proper parabolic subgroup follows from the fact that then neither does $G_K$. To see this last claim, we use induction on $\dim( R)$ to replace $R$ by some regular quotient that is a discrete valuation ring and then apply \ref{pp:BCTP-1}. As an aside, similarly, if $G_k$ is even anisotropic, then so is $G_K$. 

The anisotropicity of $G_k$ is a very stringent condition: as we now argue, it implies that 
\[
\tst \qqq G(R[\f 1r]) = G(R) \qxq{for every regular parameter $r \in R$, so also that} G(R\llp t \rrp) = G(R\llb t \rrb).
\]
For this, since $R/(r)$ is regular, the aside above implies that $G_{k_{(r)}}$ is also anisotropic. Consequently, \eqref{eqn:unique-extn} allows us to replace $R$ by its localization at the height $1$ prime $(r) \subset R$ and to thus reduce to the case of a discrete valuation ring. To then see that every $K$-point of $G$ is integral, we may even replace $R$ by its completion, at which point, since $G_K$ is still anisotropic by the aside above, $G(R) = G(K)$ by, for instance, \cite{Guo21}*{Proposition 4.4~(c)}.

\m
For minimal parabolic subgroups, that is, for Borels, the conjecture predicts that a reductive $R$-group scheme is quasi-split if and only if so is its generic fiber. By \cite{split-unramified}*{Theorem 9.5}, this consequence follows from the adjoint case of the Grothendieck--Serre conjecture \ref{conj:Grothendieck-Serre}---the latter gets used via \Cref{prop:unique-reductive-model} and the argument is similar to that of \S\ref{pp:first-implies}. 

Conversely, this consequence implies the quasi-split case of the Grothendieck--Serre conjecture over $R$: indeed, if $G$ is a reductive $R$-group scheme, $B \subset G$ is an $R$-Borel, and $E$ is a generically trivial $G$-torsor, then, since $B$-reductions of $E$ amount to Borels of the inner form $\Aut_G(E)$ of $G$ (see \S\ref{pp:parabolic}), we find that $E$ admits a generically trivial $B$-reduction; the latter must be trivial due to \eqref{eqn:M-P} and the Grothendieck--Serre conjecture for tori \S\ref{pp:basic-GS}~\ref{pp:BGS-1}.
\eenum
\epp


We now extend the result of \S\ref{pp:second-implies} by showing that \Cref{conj:CTP} also follows from the combination of the Grothendieck--Serre conjecture and the Colliot-Th\'{e}l\`{e}ne--Sansuc purity conjecture.


\bpp[\Cref{conj:Grothendieck-Serre,conj:CTS} imply \Cref{conj:CTP}]
With the notation of \Cref{conj:CTP}, assume that $G_K$ has a parabolic subgroup $P_K$ of a fixed type, let $(\bbG, \bbP)$ be a split reductive $R$-group and a parabolic subgroup of the same type as $(G_K, P_K)$ such that $(\bbG, \bbP)$ admits a pinning (see \cite{SGA3IIInew}*{Expos\'e XXVI, D\'{e}finition 1.11, Lemme 1.14}), and assume that the Grothendieck--Serre conjecture \ref{conj:Grothendieck-Serre} holds for every form of $G^\ad$. By \Cref{prop:unique-reductive-model}, then $G$ is the unique reductive $R$-group scheme with generic fiber isomorphic to  $G_K$, so, to show that $G$ has a parabolic subgroup of the desired type, all we need to do is build an $R$-form of $(\bbG, \bbP)$ with generic fiber isomorphic to $G_K$. In terms of torsors, we need to build an $\underline{\Aut}_{\gp}(\bbG, \bbP)$-torsor $\cE$ over $R$ whose $K$-fiber is isomorphic to the torsor $E$ that corresponds to $(G_K, P_K)$. 

By \S\ref{pp:known-CTP}~\ref{pp:BCTP-1}, we may assume that $E$ starts out as an $\underline{\Aut}_{\gp}(\bbG, \bbP)$-torsor over an open $U \subset \Spec(R)$ with complement of codimension $\ge 2$. Since parabolics self-normalizing (see \S\ref{pp:parabolic}), up to the center act simply transitively on pinnings adapted to them (see \cite{SGA3IIInew}*{Expos\'e~XXVI, Proposition~1.15}), and, for split groups, up to a pinning correspond to subsets of the base of positive roots (see \cite{SGA3IIInew}*{Expos\'e~XXVI, Proposition 1.4, D\'{e}finition 1.11}), the formula \eqref{eqn:semidirect} gives
\[
\underline{\Aut}_{\gp}(\bbG, \bbP) \cong \bbP^\ad \rtimes I
\]
where $\bbP^\ad \le \bbG^\ad$ is the image of $\bbP$ and $I \le \underline{\Out}_\gp(\bbG) \le \underline{\Aut}_{\gp}(\bbG)$ is the subgroup 
of  automorphisms that preserve both a fixed pinning of $\bbG$ and the subset of positive roots corresponding to $\bbP$. 

At this point, the argument becomes analogous to that of \S\ref{pp:first-implies}. Namely, $I$ is a constant $R$-group, so the $I$-torsor over $U$ induced by $E$ extends uniquely to an $I$-torsor $\cF$ over $R$. We then use the semidirect product structure to upgrade $\cF$ to an $\underline{\Aut}_{\gp}(\bbG, \bbP)$-torsor and twist by it as in \eqref{eqn:change-origin} to arrange that $E$ comes from a torsor $E'$ over $U$ under an $R$-form $\cP$ of $\bbP^\ad$. By the last paragraph of \S\ref{pp:parabolic} and the Colliot-Th\'{e}l\`{e}ne--Sansuc purity conjecture \ref{conj:CTP} applied to a Levi subgroup of this $R$-form $\cP$, there then exists a $\cP$-torsor $\cE'$ over $R$ with generic fiber $\cE'_K \simeq E'_K$. By tracing $\cE'$ back across the twisting bijection, we arrive at the desired $\underline{\Aut}_\gp(\bbG, \bbP)$-torsor $\cE$ with $\cE_K \simeq E_K$. \QED
\epp

\csub[The Nisnevich conjecture] \lab{Nisn-conj}

In search for a strategy for the Grothendieck--Serre conjecture, Nisnevich proposed to extend the Quillen conjecture \ref{conj:Quillen} to general reductive $R$-group schemes in \cite{Nis89}*{Conjecture 1.3}. However, recent examples of Fedorov \cite{Fed21c}*{Theorem 2} show that it is necessary to restrict this extension to totally isotropic groups (see \S\ref{pp:totally-isotropic}). The resulting formulation of the conjecture is as follows.

\bconj[Nisnevich] \label{conj:Nisnevich}
For a regular local ring $(R, \fm)$, a regular parameter $r \in \fm$, and a totally isotropic reductive $R$-group scheme $G$, every \emph{generically trivial} $G$-torsor over $R[\f1r]$ is trivial\ucolon
\be \label{eqn:Nisnevich-conj}
\tst \Ker(H^1(R[\f1r], G) \ra H^1(K, G)) = \{*\}, \qxq{where} K \ce \Frac(R).
\ee
\econj

We changed the original formulation by requiring generic triviality, as opposed to Zariski local triviality: this stresses the parallel with the Grothendieck--Serre conjecture \ref{conj:Grothendieck-Serre}, by which these two versions of \Cref{conj:Nisnevich} ought to be the same. One could also weaken the assumption that $r$ be a regular parameter and only require that $r \not \in \fm^2$, where $\fm \subset R$ is the maximal ideal: indeed, the new case in which $r$ is a unit is already covered by the Grothendieck--Serre conjecture \ref{conj:Grothendieck-Serre}. 




%

\bpp[Basic reductions and known cases of \Cref{conj:Nisnevich}] \label{pp:basic-Nisnevich} 
The case when $G = \GL_n$ is the Quillen conjecture (every $\GL_n$-torsor over $R[\f1r]$ is generically trivial, so one indeed recovers \Cref{conj:Quillen}), whose known cases were reviewed in \S\ref{pp:basic-Quillen}, so now we focus on~other~$G$.
\benuma
\m \label{BN-1}
The case when $G$ is a torus follows from the Grothendieck--Serre conjecture \ref{conj:Grothendieck-Serre} as follows. For a generically trivial $G$-torsor $E$ over $R[\f1r]$, we may glue it (noncanonically!) with the trivial $G$-torsor over $R_{(r)}$ and then extend the glueing to a $G$-torsor $\wt{E}$ over $R$ (see \S\ref{pp:extend}). 
By \S\ref{pp:basic-GS}~\ref{pp:BGS-1} (the Grothendieck--Serre conjecture for tori), $\wt{E}$ is trivial, so $E$ is also trivial. 

Similarly, the case when $\dim( R) = 2$ with $G$ arbitrary follows from the Grothendieck--Serre conjecture for $G$. In both of these cases, we could allow any $r \in R$, not merely $r \not\in \fm^2$, and the total isotropicity assumption is not needed.

\m \label{BN-2}
The case when $\dim( R) = 2$ and $R$ is $r$-Henselian (for instance, $r$-adically complete) may be argued as follows. As in \ref{BN-1}, the dimension assumption allows us to use \S\ref{pp:extend} to extend any generically trivial $G$-torsor over $R[\f1r]$ to a $G$-torsor over $R$ that trivializes over $R_{(r)}$. By the invariance under Henselian pairs (see \Cref{prop:Hens-pair}~\ref{HP-b}) and the Grothendieck--Serre conjecture over $R/(r)$ (see \S\ref{pp:basic-GS}~\ref{pp:BGS-2}), this extension is a trivial $G$-torsor over $R$, as desired. In this case, the total isotropicity assumption is again not needed.

\m \label{BN-3}
The case when $\dim( R) = 2$ and $G_K$ is quasi-split follows from the toral case, as we now explain (under further assumptions this case is contained in \cite{Nis89}*{Proposition 5.1}). As in \S\ref{pp:known-CTP}~\ref{pp:BCTP-1}, the dimension assumption ensures that a Borel $K$-subgroup $B \subset G_K$ extends to a Borel $U$-subgroup $\cB \subset G_U$ where $U$ is the punctured spectrum of $R$. Similarly, a generically trivial $G$-torsor over $R[\f1r]$ reduces to a generically trivial $\cB$-torsor over $R[\f1r]$, equivalently, to a generically trivial $\cT$-torsor over $R[\f1r]$ for the $U$-torus $\cT \ce \cB/\sR_u(\cB)$ (see \S\ref{pp:basic-subgroups}, \S\ref{pp:parabolic}, and \eqref{eqn:M-P}). Since $\cT$ extends to an $R$-torus (for example, by \S\ref{pp:first-implies} or simply as the ``abstract Cartan'' torus of $G$), this achieves the promised reduction to the toral case discussed in \ref{BN-1}.

\m \label{BN-2}
The case when $\dim( R) \le 3$ and $G = \PGL_n$ follows from the corresponding case of the Quillen conjecture proved by Gabber (see \S\ref{pp:basic-Quillen}~\ref{BQ-2}): indeed, granted the latter, the same argument as in \Cref{eg:Brauer-square} reduces one to the injectivity of the map $H^2(R[\f1r], \bG_m) \ra H^2(K, \bG_m)$. 

\m \label{BN-Fedorov}
The case when $R$ is of equal characteristic and either contains an infinite field or contains its own residue field was recently settled by Fedorov in \cite{Fed21c}*{Theorem 1}.



\m
Over a valuation ring $V$, one possible analogue is the statement that for any reductive $V$-group $G$, every generically trivial $G$-torsor over $V\llp t \rrp$ is trivial. This analogue was settled by Guo in \cite{Guo21}*{Corollary 7.5} (with critical input from the work of Gabber and Ramero \cite{GR18}).
\eenum

To sum up, \Cref{conj:Nisnevich} is known in almost all equicharacteristic cases but, beyond tori and some $\GL_n$ and $\PGL_n$ cases, remains widely open in mixed characteristic.
\epp


\brem
As in \S\ref{pp:Quillen-variants}~\ref{QV-1}, one may also consider a more general  variant of \Cref{conj:Nisnevich} in which $R[\f 1r]$ is replaced by $R[\f 1 {r_1\cdots r_t}]$ for a part of a regular system of parameters $r_1, \dotsc, r_t \in R$. 
\erem

Heuristically, the Nisnevich conjecture \ref{conj:Nisnevich} suggests that phenomena related to torsors under reductive groups over regular local rings may persist after inverting a regular parameter. We wish to illustrate statements of this type with the following consequence  for parabolic subgroups.

\bprop
For a regular local ring $(R, \fm)$, a regular parameter $r \in \fm$, a reductive $R$-group scheme $G$, and a parabolic $R$-subgroup $P \subset G$ such that \eqref{eqn:Nisnevich-conj} holds with $M \ce P/\sR_u(P)$ in place of $G$ \up{as is the case when $P$ is a Borel, see \uS\uref{pp:basic-Nisnevich}~\uref{BN-1}}, every parabolic $R[\f1r]$-subgroup of $G$ of the same type as $P$ is conjugate to $P_{R[\f1r]}$ by an element of $G(R[\f1r])$, equivalently,
\[
\tst G(R[\f1r])\surjects (G/P)(R[\f1r]), \qxq{equivalently,} \Ker(H^1(R[\f1r], P) \ra H^1(R[\f1r], G)) = \{*\}.
\]
\eprop

It is instructive to recall from \S\ref{pp:parabolic-torsors} that analogous statements hold when $R[\f1r]$ is replaced by any semilocal ring (but, of course, $R[\f1r]$ itself is far from being semilocal when $\dim (R) \ge 2$).

\bpf
We recall from \S\ref{pp:parabolic} that $G/P$ represents the functor that parametrizes those parabolic subgroups of (base changes of) $G$ that are of the same type as $P$, so the equivalent reformulations follow by also using the cohomology sequence \eqref{eqn:coho-seq}. For the claim itself, by \S\ref{pp:parabolic} again, the subfunctor of $G$ consisting of those sections that conjugate (a base change of) $P_{R[\f1r]}$ to a fixed parabolic $R[\f1r]$-subgroup of $G$ is a $P$-torsor over $R[\f1r]$. Since parabolics of the same type are conjugate Zariski locally on the base (see \S\ref{pp:parabolic}), this torsor is generically trivial. It then remains to note that, by \eqref{eqn:Nisnevich-conj} for $M$ and \eqref{eqn:M-P}, no nontrivial $P$-torsor over $R[\f1r]$ is generically trivial. 
\epf




\csub[The Horrocks phenomenon for totally isotropic reductive group schemes] \label{sec:Guo}

For general reductive groups, studying torsors over the relative affine line $\bA^1_A$ requires some substitute for the Horrocks principle that we reviewed in \Cref{prop:Horrocks}. For general groups, the same statement as there does not hold: in \cite{Fed16a}*{Corollary 2.3; Theorem 3, Example 2.4, Lemma~2.5}, Fedorov gave examples of regular local rings $A$ and semisimple, simply-connected $A$-group schemes $G$ for which some $G$-torsor $E$ over $\bA^1_A$ does not descend to a $G$-torsor over $A$ but is trivial away from an $A$-(finite \'{e}tale) closed subscheme $Z\subset \bA^1_A$ (so that $E$ extends to a $G$-torsor over $\bP^1_A$). 

More precisely, \emph{loc.~cit.}~shows that the Horrocks phenomenon requires some isotropicity condition on $G$. The following conjecture suggested by Ning Guo seems to capture a precise desired statement.



\bconj \label{conj:Guo}
For a commutative ring $A$ and a totally isotropic reductive $A$-group scheme $G$, every $G$-torsor over $\bA^d_A$ that is trivial away from an $A$-finite closed subscheme $Z \subset \bA^d_A$ is trivial.
\econj

In terms of analogies, this conjecture is a version of Horrocks \Cref{prop:Horrocks} beyond $G = \GL_n$.

\bpp[Basic reductions and known cases of \Cref{conj:Guo}] \label{pp:basic-Guo} \hfill
\benuma
\m \label{BG-1}
Since $Z$ is also finite over $\bA^{d - 1}_A$, induction on $d$ allows us to assume that $d = 1$. Once $d = 1$, one may use Quillen patching to reduce to local $A$, see \Cref{cor:local-to-global} for details. 


\m \label{BG-2}
The case when $G$ is a torus follows from the general formula of \cite{Hitchin-torsors}*{Theorem 3.1.7}  (essentially due to Gabber), according to which
\[
\qqq H^1(A[t], G) \oplus H^1(A, X_*(G)) \isomto H^1(A\llp t\i \rrp, G) \qxq{with} X_*(G) \ce \Hom_{\gp}(\bG_m, G).
\]
Indeed, the injectivity $H^1(A[t], G) \hra H^1(A\llp t\i \rrp, G)$ means that no nontrivial $G$-torsor over $\bA^1_A$ trivializes away from an $A$-finite closed $Z \subset \bA^1_A$ (since $Z$ is also closed in $\bP^1_A$, it does not meet the infinity section, so it also does not meet the formal neighborhood of infinity in $\bP^1_A$). 

\m \label{BG-3}
The case when $G$ is semisimple, simply connected is known, see \cite{split-unramified}*{Proposition 8.4} (possibly also \cite{PSV15}*{Theorem 1.3} for an earlier special case).

\m \label{BG-4}
The case when $G$ is split follows from the discussion in the rest of this section. 

\eenum
To sum up, \Cref{conj:Guo} is known in many cases, and so has the feeling of being within reach.
\epp

To review the strategy in the case \ref{BG-3} and to simultaneously settle the split case claimed in \ref{BG-4}, we put ourselves in the key case when $A$ is local and $d = 1$, see \ref{BG-1} above. By the following lemma, the key is to extend to a torsor over $\bP^1_A$ in such a way that the extension be trivial over the special fiber.

\blem \label{lem:Fedorov}
For a semilocal ring $S$ and a reductive $S$-group $G$ such that $\rad(G)$ is isotrivial \up{for instance, such that $G$ is either semisimple or split, or such that $S$ is normal}, every $G_{\bP^1_S}$-torsor $E$ whose base change to $\bP^1_{k_\fm}$ is trivial for every maximal ideal $\fm \subset S$ is the base change of a $G$-torsor. 
\elem

\bpf
See \cite{split-unramified}*{Lemma 8.3} for a detailed argument (possibly, see also the earlier \cite{PSV15}*{Proposition~9.6}, \cite{Fed21b}*{Proposition~2.2}, \cite{Tsy19}). One uses the assumption on $\rad(G)$ to embed $G$ into $\GL_{n,\, S}$ (see \S\ref{pp:isotriviality}) and then combines the resulting exact cohomology sequence \eqref{eqn:coho-seq} with the affineness of $\GL_{n,\, S}/G$ (see \S\ref{pp:reductive}) to reduce to the key case when $G = \GL_{n,\, S}$. The latter concerns vector bundles and its argument is similar to the proof of the Horrocks \Cref{prop:Horrocks} in that it is based on  input from \cite{EGAIII1} about cohomology and base change.
\epf

With this lemma, arguing the triviality of a $G$-torsor $E$ over $\bA^1_A$ breaks up into two steps: into the case when $A$ is replaced by its residue field $k$ and into lifting the extension of $E|_{\bP^1_k}$ to the trivial $G$-torsor over $\bP^1_k$ to an extension of $E$ to a $G$-torsor over $\bP^1_A$. These correspond to the following two lemmas, the second of which requires the total isotropy and the assumptions of \ref{BG-3} or \ref{BG-4}.

\blem \label{lem:Gille}
For a reductive group $G$ over a field $k$, a generically trivial $G$-torsor over $\bA^1_k$ is~trivial.
\elem

\bpf
The Grothendieck--Serre conjecture \ref{conj:Grothendieck-Serre} holds for discrete valuation rings (see \S\ref{pp:basic-GS}~\ref{pp:BGS-2}), so generic triviality amounts to Zariski local triviality. Gille showed in \cite{Gil02}*{Corollaire~3.10~(a)}  that every Zariski locally trivial $G$-torsor over $\bA^1_k$ reduces to a torsor under a maximal $k$-split subtorus of $G$ (see also \cite{Gil05}, and possibly compare with the earlier \cite{RR84}). Since $\bA^1_k$ has no nontrivial line bundles, it follows that every generically trivial $G$-torsor over $\bA^1_k$ is trivial.
\epf

\blem \label{lem:trivial-Gr}
For a local ring $A$ with the residue field $k$ and a totally isotropic reductive $A$-group scheme $G$ that is either split or semisimple, simply connected, the following map is surjective\ucolon
\[
G(A\llp t \rrp)/G(A\llb t \rrb) \surjects G(k\llp t \rrp)/G(k\llb t \rrb).
\]
\elem

\bpf
\addtocounter{footnote}{-7}
\renewcommand{\thefootnote}{\fnsymbol{footnote}}
The semisimple, simply-connected case was discussed in \cite{split-unramified}*{proof of Proposition 8.4}. In this case, at the cost of allowing semilocal $A$ and replacing $k$ by the product of the residue fields at the maximal ideals of $A$, the canonical decomposition \eqref{eqn:ad-decomposition} allows one to assume that $G$ is, in addition, fiberwise simple.\footnote{\emph{Added after publication.} For semilocal $A$, however, one should strengthen the total isotropicity assumption as follows: each $\widetilde G_i$ that appears in the canonical decomposition \eqref{eqn:ad-decomposition} contains $\mathbb G_{m,\, S_i}$ as an $S_i$-subgroup. Indeed, by \cite{SGA3IIInew}*{Expos\'{e} XXVI, Corollaire 6.12}, this ensures that then each $\widetilde{G}_i$ contains a fiberwise proper parabolic $S_i$-subgroup (whereas total isotropicity would supply such a parabolic only after base change to each local ring of $A$), which allows us to choose the fiberwise proper parabolic $A$-subgroups $P \subset G$ and $P^- \subset G$ in the argument below (in order to apply \cite{Gil09}*{Fait 4.3,~Lemme 4.5}, the minimality assumption on $P$ is not needed). } The key point is then the unramified nature of the Whitehead group: by \cite{Gil09}*{Fait 4.3,~Lemme 4.5}, letting $G(k\llp t \rrp)^+ \subset G(k\llp t \rrp)$ be the subgroup generated by $(\sR_u(P))(k\llp t \rrp)$ and $(\sR_u(P^-))(k\llp t \rrp)$, where $P \subset G$ is a minimal parabolic subgroup and $P^- \subset G$ is an opposite parabolic subgroup (see \S\ref{pp:parabolic}), we have
\addtocounter{footnote}{6}
\[
\tst 
G(k\llp t \rrp) = 
G(k\llp t \rrp)^+G(k\llb t \rrb).
\]
 To conclude the semisimple, simply-connected case it then suffices to note the following surjectivity: since both $\sR_u(P)$ and $\sR_u(P^-)$ are isomorphic to affine spaces $\bA^d_A$ (see \S\ref{pp:parabolic}) and $A\llp t\rrp \surjects 
k\llp t \rrp$, 
\be \label{eqn:Ru-surj}
\tst (\sR_u(P))(A\llp t \rrp) \surjects 
(\sR_u(P))(k\llp t \rrp) \qxq{and} (\sR_u(P^-))(A\llp t \rrp) \surjects 
(\sR_u(P^-))(k\llp t \rrp).
\ee
The case when $G$ is split is simpler. Then there are a split maximal $A$-torus and a Borel $A$-subgroup $T \subset B \subset G$. The Iwasawa decomposition, so, in essence, the valuative criterion of properness, gives
\[
G(k\llp t \rrp) = B(k\llp t \rrp) G(k\llb t \rrb) = (\sR_u(B))(k\llp t \rrp) T(k\llp t \rrp) G(k\llb t \rrb).
\]
Thus, \eqref{eqn:Ru-surj} applied to $\sR_u(B)$ instead reduces us to the case $G = \bG_m$. It then suffices to note that, since $A$ is local, the map $A\llp t \rrp^\times \ra k\llp t \rrp^\times \cong t^{\bZ} \times k\llb t \rrb^\times$~is~surjective.
\epf

\brem \label{rem:trivial-Gr}
\addtocounter{footnote}{-7}
\renewcommand{\thefootnote}{\fnsymbol{footnote}}
As we saw, for semisimple, simply-connected, totally isotropic\footnote{\emph{Added after publication.} Here the total isotropicity assumption should be strengthened as in the previous footnote.} $G$, \Cref{lem:trivial-Gr} continues to hold when $A$ is semilocal and $k$ is the product of its residue fields at the maximal~ideals.  
\addtocounter{footnote}{6}
\erem


\csub[The Bass--Quillen conjecture for general reductive groups] \lab{Nisn-conj}

The Bass--Quillen conjecture \ref{conj:Bass-Quillen} fails for torsors under arbitrary reductive groups $G$ in place of vector bundles (that is, in place of $\GL_n$-torsors), as examples of Parimala and others show, see \cite{Par78}, \cite{Fed16a}*{Remark 2.6}, or \cite{EKW21}*{Example 1.3}. Moreover, Fedorov's \cite{Fed16a}*{Remark~2.6} suggests excluding anisotropic groups if one aims for a positive statement. In fact, due to Balwe--Sawant \cite{BS17a}*{Proposition~4.9}, the Bass--Quillen conjecture for $G$ cannot hold over all smooth algebras over a fixed infinite perfect base field over which $G$ is defined unless $G$ is totally isotropic. 
The relevance of total isotropicity was stressed already by Raghunathan in \cite{Rag89}. 

A key feature of the groups $\GL_n$ is that all of their torsors are Zariski locally trivial. In \cite{Gro58}*{Th\'{e}or\'{e}me~3}, Grothendieck classified groups that have this property and called them \emph{special}, for instance, $\SL_n$ and $\Sp_n$ are special. For general reductive groups, one could hope that some phenomena specific to special groups may be witnessed by only considering Zariski locally trivial torsors. Over regular bases, due to the Grothendieck--Serre conjecture \ref{conj:Grothendieck-Serre}, this Zariski local triviality ought to amount to generic triviality or, if one prefers, to local triviality in the Nisnevich topology.

With these observations in mind, it seems reasonable to consider the following extension of the Bass--Quillen conjecture to torsors under more general reductive group schemes.

\bconj \label{conj:BQ-gen}
For a regular ring $R$ and a totally isotropic reductive $R$-group scheme $G$, every Zariski locally trivial $G$-torsor over $\bA^d_R$ descends to a $G$-torsor over $R$. 
\econj

\brem \label{rem:BQG-1}
Induction on $d$ reduces \Cref{conj:BQ-gen} to $d = 1$ and Quillen patching of \Cref{cor:Quillen-patching} reduces further to local $R$. Once $R$ is local and $d = 1$, the conjecture may be strengthened to predict that every Zariski locally trivial $G$-torsor over $\bA^1_R$ is trivial:  by \Cref{lem:Gille} applied to $\Frac(R)$ and by the Grothendieck--Serre conjecture \ref{conj:Grothendieck-Serre}, this ought to give an equivalent statement.
\erem

Analogously to how the Quillen conjecture \ref{conj:Quillen} and the Horrocks \Cref{prop:Horrocks} imply the Bass--Quillen conjecture \ref{conj:Bass-Quillen} (see \S\ref{sec:Quillen}), the Nisnevich conjecture \ref{conj:Nisnevich} and the Horrocks phenomenon stated in \Cref{conj:Guo} imply the extension of the Bass--Quillen conjecture above as follows.

\bpp[Conjectures \ref{conj:Nisnevich} and \ref{conj:Guo} imply \Cref{conj:BQ-gen}] \label{pp:BQ-implication}
By \Cref{rem:BQG-1}, for a regular local ring $R$, a totally isotropic reductive $R$-group $G$, and a Zariski locally trivial $G$-torsor $E$ on $\bA^1_R$, we need to argue that $E$ is trivial. The Nisnevich \cref{conj:Nisnevich} implies that $E$ becomes trivial over the punctured neighborhood of infinity in $\bP^1_R$, that is, after inverting all the monic polynomials in $R[t]$ (compare with the proof of \Cref{Lam-conj}). Thus, by a limit argument, $E$ also becomes trivial  after inverting a single monic polynomial. \Cref{conj:Guo} applied with $Z$ being the vanishing locus of that polynomial then implies that $E$ is trivial to begin with. \QED
\epp


\bpp[Basic reductions and known cases of \Cref{conj:BQ-gen}] \label{pp:basic-BQ-gen} 
The case $G = \GL_n$ is the Bass--Quillen conjecture \ref{conj:Bass-Quillen}, whose known cases were reviewed in \S\ref{pp:basic-BQ}, so now we focus on~other~$G$.
\benuma


\m \label{BQG-4}
By \S\ref{pp:BQ-implication}, the case when $G$ is a torus follows from \S\ref{pp:basic-Nisnevich}~\ref{BN-1} and \S\ref{pp:basic-Guo}~\ref{BG-2}.


\m \label{BQG-3}
The case when $R$ is a smooth algebra over a field $k$ and $G$ is the base change of a totally isotropic reductive $k$-group was settled by Asok--Hoyois--Wendt in \cite{AHW18}*{Theorem 3.3.7} when $k$ is infinite and in \cite{AHW20}*{Theorem 2.4} when $k$ is finite. Their argument follows an axiomatic approach of Colliot-Th\'{e}l\`{e}ne--Ojanguren \cite{CTO92}, and they check a crucial  Nisnevich excision axiom using methods of $\bA^1$-homotopy theory of Morel--Voevodsky. One may also check this Nisnevich excision more directly by using \Cref{lem:patch}, see \cite{Li21}.

\m
By \S\ref{pp:BQ-implication}, the case when $R$ contains an infinite field and $G$ is either quasi-split, semisimple, and simply connected or split follows from \S\ref{pp:basic-Nisnevich}~\ref{BN-Fedorov} and \S\ref{pp:basic-Guo}~\ref{BG-3}--\ref{BG-4}. 

\m \label{BQG-4}
By \S\ref{pp:BQ-implication}, the case when $\dim( R) \le 2$, $d = 1$, and $G = \PGL_n$ follows from \S\ref{pp:basic-Nisnevich}~\ref{BN-2} and \S\ref{pp:basic-Guo}~\ref{BG-4}. 

\m \label{BQG-4}
By \S\ref{pp:BQ-implication}, the case when $\dim( R) = 1$, $d = 1$, and $G$ is either quasi-split, semisimple, and simply connected or split follows from \S\ref{pp:basic-Nisnevich}~\ref{BN-3} and \S\ref{pp:basic-Guo}~\ref{BG-3}--\ref{BG-4}. 

\eenum

To sum up, these known cases give evidence for \Cref{conj:BQ-gen}, especially in equal characteristic, but the general case remains open even in the setting when $R$ is a smooth algebra over a field $k$ and the totally isotropic reductive $R$-group scheme $G$ does not descend to $k$.
\epp


\section{Passage to the affine space via presentation lemmas and excision} \label{sec:ring-structure}

Most of the conjectures discussed in \S\S\ref{sec:vector-bundles}--\ref{sec:general-groups} are specific to regular bases. Crudely speaking, one broad strategy for attacking them has been to reduce to the ``geometric'' case in which the base ring is a localization of a smooth algebra over a ring $k$ that is either a field or $\bZ$,  and to then combine suitable preparation lemmas with excision to reduce further to working over a localization of some affine space $\bA^d_k$. This strategy is specific to unramified regular base rings because the reduction to the geometric case is based on the Popescu theorem \ref{thm:Popescu}, which requires the unramifiedness assumption. 

We elaborate on this strategy in the sections that follow. In \S\ref{sec:Noether}, we overview some of the presentation lemmas that have been used for building the required maps to $\bA^d_k$. In \S\ref{sec:excision}, we review the relevant excision input that allows subsequent passage to the local rings of $\bA^d_k$.



\csub[Presentation lemmas] 
\label{sec:Noether}

By \cite{EGAIV4}*{Corollaire 17.11.4}, Zariski locally on the source every smooth map of schemes factors as an \'{e}tale morphism to a relative affine space, in particular, for any ring $k$ and a smooth $k$-algebra $R$, upon localizing around a fixed point of $\Spec(R)$, there is an \'{e}tale map
\be \label{eqn:smooth-local-structure}
f\colon \Spec(R) \ra \bA^d_k.
\ee
The goal of presentation lemmas is to build such a map subject to additional requirements, for instance, it may be handy to have that $f$ induce an isomorphism on the residue fields $k_{f(s)} \isomto k_s$ at some point $s \in \Spec(R)$ of interest or that a specified closed subscheme $Z \subset \Spec(R)$ of smaller dimension would be \emph{finite} over the affine space $\bA^{d - 1}_k$ given by the first $d - 1$ coordinates. In practice, $R$ comes equipped with a torsor under some reductive $R$-group and one knows that ``something good'' happens, for instance, the torsor trivializes, over a dense open $U \subset \Spec(R)$; the control of the complement $Z \ce \Spec(R) \setminus U$ via some presentation lemma is then crucial for reducing the problem to its counterpart for the affine space. Presentation lemmas are most developed in the literature in the case when $k$ an (often infinite) field, for instance, the following result gives a broadly useful statement in this setting. It grew out of refinements due to Gabber \cite{Gab94b}*{Section 3} and Gros--Suwa \cite{GS88}*{Section 2} to a basic such lemma used by Quillen \cite{Qui73}*{Section 7, Lemma 5.12}.

\bthm[Geometric Presentation Theorem] \label{thm:geometric-presentation}
For a smooth, affine, irreducible scheme $X$ of dimension $d > 0$ over an infinite field $k$, a closed subscheme $Z \subset X$ of codimension $> 0$, and $x_1, \dotsc, x_n \in X$, there are a $k$-map $f\colon X \ra \bA^{d}_k$ making $Z$ finite over the $\bA^{d - 1}_k$ of the first $d - 1$ coordinates and an open $X' \subset X$ containing $x_1, \dotsc, x_n$ that fit into a commutative diagram 
\[
\xymatrix@R=12pt@C=8pt{
  \ar@{}[d]|-{\x{\begin{turn}{-90}\scalebox{1.25}{$\cong$}\end{turn}}} X' \cap Z \ar@{^(->}[r] &X' \ar[rrr]^{f|_{X'}}\ar[ddd] \ar@{}[rd]|-{\x{\begin{turn}{-45}\scalebox{2.25}{\,$\subset$}\end{turn}}} &&& \bA^d_k \ar[ddd]^{(t_1,\, \dotsc,\, t_d)\, \mapsto\, (t_1,\, \dotsc,\, t_{d - 1})} \\
g\i(S) \cap Z \ar[rdd] && X \ar[rrdd]_-g \ar[urr]_-f && \\ 
&&&& \\
&S \ar@{^(->}[rrr] &&& \bA^{d - 1}_k
}
\]
in which $f|_{X'}$ is \'{e}tale, $S \subset \bA^{d- 1}_k$ is an open, $f(x_i) \not \in f(Z)$ if $x_i \not \in Z$, and $f|_{X'}$ maps $X' \cap Z$ isomorphically to a closed subscheme $Z' \subset \bA^1_S$ in such a way that the following square is Cartesian\ucolon
\[
\xymatrix@C=10pt{
X' \cap Z \ar@{}[r]|-{\x{\scalebox{1.25}{$\cong$}}} & (f|_{X'})\i(Z') \ar@{^(->}[d] \ar[r]^-{\sim} & Z' \ar@{^(->}[d] \\
& X' \ar[r]^{f|_{X'}} & \bA^1_S.
}
\]
\ethm

\bpf
A detailed proof was given by Colliot-Th\'{e}l\`{e}ne--Hoobler--Kahn, see \cite{CTHK97}*{Theorem 3.1.1}. Since $X$ is affine, there is some affine space over $k$ that parametrizes maps $X \ra \bA^{d}_k$. \emph{Loc.~cit.}~shows that the desired conditions describe a nonempty open of this ``moduli space.'' Since $k$ is infinite, this nonempty open  has a $k$-rational point, which corresponds to the desired $f$.
\epf

\brem
For a version of \Cref{thm:geometric-presentation} in a setting where $k$ is replaced by a discrete valuation ring $\cO$, see \cite{split-unramified}*{Variant 3.7}. There it was convenient to even allow semilocal Dedekind $\cO$: the method was to bootstrap from \Cref{thm:geometric-presentation} applied to the closed $\cO$-fibers of $X$, but one of the difficulties was that some of the $x_i$ may lie in the generic $\cO$-fibers; to overcome it, we enlarged $\Spec(k)$ by glueing in an auxiliary discrete valuation ring along which such $x_i$ specialized well.
\erem

\Cref{thm:geometric-presentation} is particularly useful for dealing with imperfect base fields $k$, for which other techniques run into difficulties. However, it needs $k$ to be infinite; on the other hand, to treat finite $k$, one may adapt Artin's technique of good neighborhoods from \cite{SGA4III}*{Expos\'{e} XI}. Artin's method is more direct: roughly, instead of considering a moduli space of maps $X \ra \bA^d_k$ as in the proof of \Cref{thm:geometric-presentation}, it uses Bertini theorem to directly construct hypersurfaces $H_i \subset X$ whose defining functions are the images of the coordinates of the affine space under the desired map $X \ra \bA^d_k$. A detailed implementation of such an approach may be found in \cite{split-unramified}*{proof of Proposition 3.6, around equation (3.6.2)}; here we content ourselves with reviewing the relevant Bertini statement.

\bprop \label{lem:Bertini}
For a projective scheme $X$ of pure dimension over a field $k$, a nowhere dense closed subscheme $Z \subset X$, and a $t \le \dim(X) - \dim(Z)$ such that even $t \le \dim(X)$ when $Z = \emptyset$, there are hypersurfaces $H_1, \dotsc, H_t \subset X$ with respect to a fixed ample line bundle $\sO_X(1)$ whose intersection
\[
H_1 \cap \dotsc \cap H_t
\]
is of pure dimension $\dim(X) - t$, contains $Z$, and has a $k$-smooth intersection with $X^\sm\setminus Z$. Moreover, we may simultaneously achieve the following additional requirements\ucolon
\benum
\m
for closed subschemes $Y_1, \dotsc, Y_n \subset X$, the $H_i$ intersect each $Y_j \setminus Z$ transversally in the sense that $\dim((Y_j \setminus Z) \cap \bigcap_{i \in I} H_i) \le \dim(Y_j \setminus Z) - \#I$ for all $1 \le j \le n$ and $I \subset \{1, \dotsc, t\}$\uscolon

\m \label{m:Bertini-b}
if $Z = Z_1 \sqcup Z_0$ for a $0$-dimensional $Z_0 \subset X^\sm$ all of whose residue fields are separable extensions of $k$, then $H_1 \cap \dotsc \cap H_t$ is $k$-smooth even at the points in $Z_0$\uscolon

\m \label{m:Bertini-c}
iteratively on $i$, with $H_1, \dotsc, H_{i - 1}$ already fixed, $H_i$ may be chosen to have any sufficiently large degree divisible by the characteristic exponent of $k$. 
\eenum
\eprop

\bpf
The lemma is essentially a restatement of \cite{split-unramified}*{Lemma 3.2}, whose proof split into the characteristic $0$ and the positive characteristic cases. In both cases, the conclusion followed from Bertini theorem (which supplies a single $H_i$), although the slightly nonstandard requirements $Z \subset H_i$ and \ref{m:Bertini-b} required further care. In positive characteristic, the key input was Gabber's version of Bertini theorem over finite fields from \cite{Gab01}, notably, it allowed us to arrange \ref{m:Bertini-c}. The latter seems less straight-forward to obtain from Poonen's version of Bertini theorem over finite fields from \cite{Poo04} (which is sharper in other aspects that are not relevant for the present lemma).
\epf

To give a small illustration of the Geometric Presentation \Cref{thm:geometric-presentation} 
in practice, we present the following refinement of the local structure of smooth maps of schemes.

\bprop[Compare with \cite{Lin81}*{Proposition 2}] \label{cor:baby-Lindel}
For a local ring $k$ and a $k$-algebra $R$ that is smooth  of relative dimension $d > 0$ at a maximal ideal $\fm \subset R$ that lies over the maximal ideal of $k$, there are an affine open $S \subset \Spec(R)$ containing $\fm$ and an \'{e}tale $R$-map 
\[
f\colon S \ra \bA^d_k \qxq{that induces an isomorphism on residue fields} k_{f(\fm)} \isomto k_{\fm}.
\]
\eprop

\bpf
The main point is the aspect about isomorphisms on residue fields: without it and under the additional assumption that $k_{\fm}$ is a separable extension of the residue field of $k$, it would suffice to choose global sections $t_1, \dotsc, t_d \in R$ that form a regular system of parameters at $\fm$ in the closed $k$-fiber of $R$, to define $f$ by sending the standard coordinates of $\bA^d_k$ to the $t_i$, and to combine  \cite{EGAIV4}*{Proposition 17.15.8} with the fibral criterion of flatness \cite{EGAIV3}*{Th\'{e}or\`{e}me 11.3.10} to conclude that this $f$ is \'{e}tale at $\fm$. Even though this method does not suffice for us, the fibral criterion of flatness of \emph{loc.~cit.}~and the openness of the \'{e}tale locus do allow us to replace $k$ by its residue field, and hence to assume for the rest of the proof that $k$ is a field. 

Once $k$ is a field, if it is also infinite, then the conclusion is a special case of \Cref{thm:geometric-presentation} with $Z \subset X$ there being our $\{\fm\} \subset \Spec(R)$. In the remaining case when $k$ is a finite field, the finite extension $k_{\fm}/k$ is automatically separable, so the method of the previous paragraph at least gives us an $R$-morphism $\Spec(R) \ra \bA^{d - 1}_k$ that is smooth of relative dimension $1$ at $\fm$ (this time send the standard coordinates to $t_1, \dotsc, t_{d - 1}$). The technique of the previous paragraph now allows us to replace $\Spec(R)$ by its fiber over the origin of $\bA^{d - 1}_k$ to reduce further to the case when $d = 1$ (and $k$ is still a finite field). In this case, $\Spec(R)$ is an affine curve over $k$ and the first infinitesimal neighborhood $\eps_\fm$ of its point $\fm$ is isomorphic to $k_\fm[t]/(t^2)$. In particular, we may embed $\eps_\fm$ into $\bA^1_k$ and then choose an $r \in R$ whose image in the coordinate ring of $\eps_\fm$ agrees with the image of the standard coordinate of $\bA^1_k$. By sending this standard coordinate to $r$, we obtain a desired $k$-morphism $f \colon \Spec(R) \ra \bA^1_k$ that is \'{e}tale at $\fm$ and induces an isomorphism $k_{f(\fm)} \isomto k_\fm$. 
\epf


 The morphism $f$ as in \Cref{cor:baby-Lindel} is not only an isomorphism on residue fields but is automatically even an isomorphism along a germ of a local hypersurface as follows.

\blem[Lindel] \label{lem:Lindel}
For an \'{e}tale, local homomorphism $R_0 \ra R$ of local rings that induces an isomorphism on residue fields, there is a nonunit $r \in R_0$ such that 
\[
R_0/r^nR_0 \isomto R/r^nR \qxq{for all} n > 0.
\]
\elem

\bpf
We follow Lindel's argument from \cite{Lin81}*{Lemma on p.~321}. Namely, by \cite{SP}*{Proposition~\href{https://stacks.math.columbia.edu/tag/00UE}{00UE}}, the morphism is standard \'{e}tale, more precisely, $R$ is a localization of the quotient $R_0[T]/(h(T))$ for some $h(T) \in R_0[T]$ whose derivative $h'(T) \in R_0[T]$ is a unit in $R$. Since $R_0 \ra R$ is an isomorphism on residue fields, we may change variables and arrange that $T$ lie in the maximal ideal of $R$. It then follows that 
\[
r \ce h(0)
\]
lies in the maximal ideal of $R_0$. The condition on the derivative $h'(T)$ ensures that $R/rR \cong R_0/rR_0$ and, by \'{e}taleness, then $R_0/r^nR_0 \isomto R/r^nR$ for all $n > 0$ (see \cite{SP}*{Theorem~\href{https://stacks.math.columbia.edu/tag/039R}{039R}}).
\epf

The idea of embedding $\eps_\fm$ into the affine line at the end of the proof of \Cref{cor:baby-Lindel} is also central in the proof of the following preparation result for relative curves. This result is used in proving cases of the Grothendieck--Serre conjecture \ref{conj:Grothendieck-Serre}: its role is to produce an excision square, which one then combines with patching discussed in \S\ref{sec:excision} to reduce to studying torsors over $\bA^1_R$.

\bprop \label{prop:embed}
Let $R$ be a semilocal ring, let $C$ be an affine $R$-scheme that is smooth of pure relative dimension $1$ \up{an $R$-curve}, and let $Z \subset C$ be an $R$-finite closed subscheme such that, for every maximal ideal $\fm \subset R$ whose residue field is finite, $Z_{k_\fm}$ is connected. There are an affine open $C' \subset C$ containing $Z$ and an \'{e}tale morphism $f \colon C' \ra \bA^1_R$ that maps $Z$ isomorphically to a closed subscheme of $\bA^1_R$ whose scheme-theoretic preimage in $C'$ is $Z$, so that we have Cartesian squares
\be\label{eqn:Cpr-square}\ba 
\xymatrix{
C' \setminus Z \ar[d] \ar@{^(->}[r] & C' \ar[d]^-f & \ar@{_(->}[l] f\i(Z) \ar[d]^{\sim} \\
\bA^1_R \setminus Z \ar@{^(->}[r] & \bA^1_R & \ar@{_(->}[l] Z.
}
\ea\ee
\eprop

\bpf
The claim is a special case of \cite{split-unramified}*{Lemma 6.3}, which is itself a generalization of earlier versions given by Panin and his collaborators in \cite{OP99}*{Section 5}, \cite{PSV15}*{Theorem 3.4}, \cite{Pan19a}*{Theorem~3.8}. To the best of our knowledge, it was Panin who introduced this type of statement.    

The argument given in \cite{split-unramified}*{Lemma 6.3} is not long (and is self-contained):  the idea is to embed each $Z_{k_\fm}$ into $\bA^1_{k_\fm}$, to then lift these to an $R$-embedding $Z \hra \bA^1_R$ using the Nakayama lemma (the $R$-finiteness of $Z$ is critical for this), and, finally, to build $C'$ using prime avoidance. The construction of the embedding $Z_{k_\fm} \hra \bA^1_{k_\fm}$ rests on the following corollary that one shows first and that deserves to be known more widely: for a closed point of a \emph{smooth} curve over a field $k$, 
its residue field $k'$ is also the residue field of some point of $\bA^1_k$ (equivalently, $k'$ is generated by a single element as a field extension of $k$; equivalently, $k'/k$ has only finitely many field subextensions).
\epf

\brem \label{rem:number-pts}
The version in \cite{split-unramified}*{Lemma 6.3} is more general, for instance, the assumption on $Z_{k_\fm}$ may be weakened to that, for every $d \ge 1$, the number of points of $Z_{k_\fm}$ with residue field of degree $d$ over $k_\fm$ is at most its counterpart for $\bA^1_{k_\fm}$ (as is automatically the case if $k_\fm$ is infinite).
\erem


\csub[Excision and patching techniques] 
\label{sec:excision}

An important and often used technique for studying torsors is excision, sometimes also called patching or formal glueing. For instance, it is used in conjunction with squares such as \eqref{eqn:Cpr-square} to descend a torsor over $C'$ to a torsor over $\bA^1_R$. The following proposition gives a general basic excision result.


\bprop \label{lem:patch}
Let $U \subset X$ be an open immersion of schemes, let $Z \subset X$ be a complementary closed subscheme that is locally cut out by a finitely generated ideal, and consider Cartesian squares
\be\label{eqn:formal-square}\ba 
\xymatrix{
f\i(U) \ar[d] \ar@{^(->}[r] & X' \ar[d]^-f & \ar@{_(->}[l] f\i(Z) \ar[d]^{\sim} \\
U \ar@{^(->}[r] & X & \ar@{_(->}[l] Z
}
\ea\ee
in which the morphism $f$ is affine, flat, and induces an isomorphism over $Z$ as indicated. For a quasi-affine, flat, finitely presented $X$-group scheme $G$, base change induces an equivalence of~categories
\[
\{\x{$G$-torsors over $X$}\} \isomto \{\x{$G$-torsors over $X'$}\} \times_{\{\x{$G$-torsors over $f\i(U)$}\}} \{\x{$G$-torsors over $U$}\},
\]
in other words, giving a $G$-torsor over $X$ amounts to giving $G$-torsors over $X'$ and $U$ together with an isomorphism between their base changes to $f\i(U)$ \up{and likewise for $G$-torsor \up{iso}morphisms}.
\eprop

\bpf
The claim is a special case of general patching results of Moret-Bailly \cite{MB96}*{Corollaire~6.5.1~(a)} applied to the classifying stack $\bbB G$, although it also follows from earlier Ferrand--Raynaud patching of modules \cite{FR70}*{Proposition 4.2}. See also \cite{split-unramified}*{Lemma 7.1} for comments on why $\bbB G$ with $G$ as in the statement satisfies the general assumptions of Moret-Bailly.
\epf

The preceding result is general but its requirement that $f$ be flat is too restrictive in some situations, especially, in non-Noetherian settings. In some such cases, one may instead use the following result, which refines the widely-known Beauville--Laszlo patching from \cite{BL95}.

\bprop\label{lem:patch-nonflat}
For a ring $A$, an $a \in A$, a ring map $f\colon A \ra A'$ that is an isomorphism on derived $a$-adic completions \up{concretely, this means that $f$ induces an isomorphism both modulo $a^n$ and on $a^n$-torsion for every $n > 0$}, and a quasi-affine, flat $A$-group scheme $G$, we have an equivalence
\[
\tst \{\x{$G$-torsors over $A$}\} \isomto \{\x{$G$-torsors over $A'$}\} \times_{\{\x{$G$-torsors over $A'[\f1a]$}\}} \{\x{$G$-torsors over $A[\f1a]$}\},
\]
in other words, giving a $G$-torsor over $A$ amounts to giving $G$-torsors over $A'$ and $A[\f1a]$ together with an isomorphism between their base changes to $A'[\f1a]$ \up{and likewise for $G$-torsor \up{iso}morphisms}. 
\eprop

\bpf
The claim is \cite{Hitchin-torsors}*{Lemma 2.2.11~(b)} that is due to de Jong. Its main input are the results from \cite{SP}*{Section~\href{https://stacks.math.columbia.edu/tag/0F9M}{0F9M}} that were inserted into the Stacks Project to facilitate this proof. 
\epf


\section{The analysis of torsors over  the relative affine line} 
\label{sec:affine-space}

The final stages of the known approaches to problems about torsors over regular rings usually involve the analysis of torsors over the relative affine line $\bA^1_R$ (so, by changing $R$, also over the relative affine space $\bA^d_R$). The techniques for studying torsors over $\bA^1_R$  tend to work for any base ring $R$, although they sometimes require $R$ to be local. In practice, this last requirement is not stringent: one reduces to local $R$ via Quillen patching, which we review in its general form in \S\ref{sec:Quillen-patch}, see \Cref{cor:Quillen-patching}. We illustrate the utility of Quillen patching by reviewing the proof of the unramified case of the Bass--Quillen conjecture in \S\ref{sec:unramified-BQ}. Once $R$ is local, a key technique that permits the analysis of torsors over $\bA^1_R$ beyond semisimple, simply-connected groups is due to Fedorov and is based on the geometry of the affine Grassmannian; we review it in \S\ref{sec:A1-analysis}.





\csub[Quillen patching for general groups] \label{sec:Quillen-patch}

A central technique for studying torsors over $\bA^1_R$ is a local-to-global principle known as Quillen patching. In \Cref{cor:Quillen-patching}, we show that Quillen patching holds for any locally finitely presented group scheme $G$. The argument is not long and its key insight is due to Gabber. The result is much more general than what has appeared in the literature: the Quillen case \cite{Qui76}*{Theorem 1} is $G = \GL_n$ and his proof was extended to arbitrary finitely presented closed subgroups of $\GL_n$ by Moser \cite{Mos08}*{Satz 3.5.1} and Asok--Hoyois--Wendt \cite{AHW18}*{Theorem 3.2.5}. Quillen's insightful and in essence elementary technique was axiomatized by Bass--Connell--Wright in \cite{BCW76}, who isolated a crucial ``axiom Q'' that ensures that patching still holds for $G$-torsors: for any $r \in R$ and any $g \in G(R[\f 1r][T])$ that reduces to the identity section modulo $T$, there ought to be an $n > 0$ such that the $R[\f1r]$-algebra automorphism of $R[\f1r][T]$ determined by $T \mapsto r^nT$ brings $g$ to the image of some $\wt{g} \in G(R[T])$ whose reduction modulo $T$ is the identity section. It seems to have been overlooked in the literature for a long time that this axiom is straight-forward to verify when $G$ is affine and of finite presentation by considering generators and relations of the coordinate ring. Gabber's insight is deeper: he noticed that ``axiom Q'' holds even when $G$ is a locally finitely presented $R$-algebraic space because one may check it by using a result of Temkin--Tyomkin \cite{TT16}*{Theorem 4.3}, according to which the functor $G(-)$ commutes with fiber products of rings $A_1 \times_{A_0} A_2$ provided that one of the maps $A_1 \ra A_0$ and $A_2 \ra A_0$ is surjective (one applies this to the maps $R[\f1r][T] \surjects R[\f1r]$ and $R \ra R[\f1r]$). We explain his observation in more detail in the following lemma.

%

\blem \label{lem:axiom-Q}
For a ring $R$, an $R$-algebraic space $X$ locally of finite presentation, an $r \in R$, an 
\[
\tst x(T) \in X(R[\f1r][T]) \qxq{whose pullback} x_0 \in X(R[\f 1r]) \qxq{along} T \mapsto 0 \qxq{lifts to an} \wt{x}_0 \in X(R),
\]
and every large enough $n \ge 0$ \up{that depends on $x(T)$}, the section $x(r^nT) \in X(R[\f1r][T])$ lifts to an 
\[
\wt{x}(T) \in X(R[T]) \qxq{whose pullback along} T \mapsto 0 \qxq{is} \wt{x}_0 \in X(R).
\]
\elem

Here and below we let $x(r^nT)$ denote the pullback of $x(T)$ along the map $R[\f1r][T] \xra{T\,\mapsto\, r^nT} R[\f1r][T]$. 

\bpf
By \cite{TT16}*{Lemma 4.1 and Theorem 4.3}, $\Spec(-)$ transforms fiber products of rings into pushouts in the category of algebraic spaces granted that one of the two maps of which the fiber product is formed is surjective. We apply this to the fiber product $R[\f 1r][T] \times_{R[\f1r]} R$, in which the first map is the surjection given by $T \mapsto 0$, to conclude that $x$ and $\wt{x}_0$ assemble to a unique section
\[
\tst x'(T) \in X(R[\f 1r][T] \times_{R[\f1r]} R). 
\]
Concretely, $R[\f 1r][T] \times_{R[\f1r]} R$ may be thought of as the ring of polynomials in $R[\f 1r][T]$ whose constant coefficient is equipped with a lift to an element of $R$, in other words,
\[
\tst R[\f 1r][T] \times_{R[\f1r]} R \cong \varinjlim_{T \mapsto rT} R[T]. 
\]
Since $X$ is locally of finite presentation, applying $X(-)$ commutes with this filtered direct limit, to the effect that $x'(T)$ lifts to the resulting $n$-th copy of $X(R[T])$ for every large $n \ge 0$ in such a way that this lift restricts to $\wt{x}_0$ along $T \mapsto 0$. Then $x'(r^nT)$ lifts to a desired $\wt{x}(T)$ in the $0$-th (that is, in the initial) copy of $X(R[T])$ in the filtered direct limit $\varinjlim_{T \mapsto rT} X(R[T])$. 
\epf

\bpp[The notation $G_0$]
In the rest of this section, for a ring $R$ and a group-valued functor $G$ on the category of $R$-algebras $R'$, we will let $G_0$ be the group-valued functor defined by
\[
G_0(R') \ce \Ker(G(R'[T]) \xra{T\,\mapsto\, 0} G(R')). 
\]
This shorthand notation is nonstandard, but it will be convenient. The following generalization will also be useful: for a $\bZ_{\ge 0}$-graded $R$-algebra $A = A_0 \oplus A_1 \oplus \dotsc$, we let $G_0^A$ be the functor defined by
\[
G_0^A(R') \ce  \Ker(G(A \tensor_R R') \ra G(A_0 \tensor_R R')). 
\]
As a concrete example, $A$ could be some $R$-subalgebra of $R[T]$ generated by monomials.
\epp

\bprop \label{prop:Zariski-patch}
For a ring $R$, a group-valued functor $G$ on the category of $R$-algebras $R'$ that has a locally finitely presented $R$-subgroup algebraic space $G' \subset G$ as an open subfunctor that is $R$-fiberwise clopen in $G$ \up{main example\ucolon $G' = G$}, and elements $r, r' \in R$ that generate the unit ideal,
\[ 
\tst G_0(R[\f 1{rr'}]) = G_0(R[\f1r]) G_0(R[\f1{r'}])  \qxq{and, more generally,} G_0^A(R[\f 1{rr'}]) = G_0^A(R[\f1r]) G_0^A(R[\f1{r'}])
\]
for every $\bZ_{\ge 0}$-graded $R$-algebra $A = A_0 \oplus A_1 \oplus \dotsc$ that satisfies $R \isomto A_0$ \up{our notation is slightly abusive, since localization maps such as $G_0(R[\f1r]) \ra G_0(R[\f 1{rr'}])$ need not be injective}.
\eprop


\bpf
We focus on the claim about $A$ because it includes the case when $A = R[T]$ with the grading given by the degree. Any idempotent in a $\bZ_{\ge 0}$-graded ring is homogeneous of degree $0$, so $\Spec(A)$ has connected $R$-fibers. Thus, since the identity section of $G$ lies in $G'$, we have $G_0^A = G^{\prime A}_0$. In conclusion, we may replace $G$ by $G'$ and assume that $G' = G$. 

The claim now follows from \cite{BCW76}*{Corollary 2.7}. Indeed, by \Cref{lem:axiom-Q}, the functor $G$ satisfies ``axiom Q'' formulated in \cite{BCW76}*{Axiom~1.1}: explicitly, for an $R$-algebra $R'$, an $r \in R'$, and a
\[
\tst g(T) \in \Ker\p{G(R'[\f1{r}][T]) \xra{T\,\mapsto\,0} G(R'[\f1{r}])}\!,
\]
there is an $n \ge 0$ such that $g(r^nT)$ lifts to an element of $\Ker\p{G(R'[T]) \xra{T\,\mapsto\,0} G(R')}$.
\epf

\bcor\lab{cor:H1Zar-G0}
For a ring $R$ and a group functor $G$ as in Proposition \uref{prop:Zariski-patch}, we have
\[
H^1_\Zar(R, G_0) = \{*\} \qxq{and, with $A$ as there, also} H^1_\Zar(R, G_0^A) = \{*\}.
\]
\ecor

\bpf
With \Cref{prop:Zariski-patch} in hand, we merely need to follow the argument of \cite{Qui76}*{proof of Theorem 1}. Namely, we need to show that every Zariski locally trivial $G_0$-torsor $X$ is trivial, so we let $S \subset R$ be the subset of those $r \in R$ such that $X$ trivializes over $R[\f1r]$. It suffices to show that $S$ is an ideal, since then the Zariski local triviality will imply that $S = R$, and for this it is enough to argue that $r + r' \in S$ whenever $r, r' \in S$. Moreover, by replacing $R$ by $R[\f1{r + r'}]$, we may assume that $r, r'$ generate the unit ideal. However, then \Cref{prop:Zariski-patch} applies and, in terms of Zariski descent, implies that no nontrivial $G_0^A$-torsor trivializes over both $R[\f1r]$ and $R[\f1{r'}]$, as desired.
\epf

\bcor[Gabber] \label{cor:Quillen-patching}
Let $R$ be a ring and let $G$ be a group functor as in Proposition \uref{prop:Zariski-patch} \up{for instance, $G$ could be any locally finitely presented $R$-group scheme}.
\benum
\m \label{QP-a}
For a $G$-torsor $X$ over $R[t_1, \dotsc, t_d]$, 
the set $S \subset R$ of those $r \in R$ such that $X|_{(R[t_1,\dotsc, t_d])[\f1r]}$ descends to a $G$-torsor over $R[\f1r]$ is an ideal. 
 
\m \label{QP-b}
A $G$-torsor over $R[t_1, \dotsc, t_d]$ descends to a $G$-torsor over $R$ iff it does so Zariski locally on~$R$. 
 
 \eenum
 More generally, the analogues of \ref{QP-a} and \ref{QP-b} hold with $R[t_1, \dotsc, t_d]$ replaced by any $\bZ^{\oplus d}_{\ge 0}$-graded $R$-algebra $A \cong \bigoplus_{i_1, \dotsc, i_d \ge 0} A_{i_1, \dotsc, i_d}$ such that $R \isomto A_{0, \dotsc, 0}$. 
\ecor

Of course, the main case of interest is when $G$ is locally finitely presented $R$-group scheme. However, better descent properties of algebraic spaces (see \S\ref{pp:basic-representability}) make the added generality quite useful: for instance, if one wishes to pass to a form of $G$, one needs not worry about the form being a scheme. 

\bpf
Evidently, \ref{QP-a} implies \ref{QP-b}, so we only focus on \ref{QP-a} and follow the proof of \cite{AHW18}*{Proposition~3.2.4}. Namely, as in the proof of \Cref{cor:H1Zar-G0}, it suffices to show that $r + r' \in S$ whenever $r, r' \in S$, and we may assume that $r, r'$ generate the unit ideal, so that we seek to show that $X$ descends to a $G$-torsor over $R$. Induction on $d$ then allows us to assume that $d = 1$, that is, that we are dealing with a $G$-torsor $X$ over a $\bZ_{\ge 0}$-graded $R$-algebra $A \cong A_0 \oplus A_1 \oplus \dotsc$ with $R \isomto A_0$. The sought descent of $X$ to a $G$-torsor over $R$ will have to be $X|_{A_0}$, and we have $G$-torsor isomorphisms 
\[
\gA \colon (X|_{A_0})_{A[\f1r]} \isomto X|_{A[\f1r]} \qxq{and} \gA' \colon (X|_{A_0})_{A[\f1{r'}]} \isomto X|_{A[\f1{r'}]}.
\]
By adjusting these isomorphisms by elements of $G(A_0[\f1r])$ and $G(A_0[\f1{r'}])$, we may assume that both $\gA|_{A_0}$ and $\gA'|_{A_0}$ are the identity isomorphisms. The isomorphisms $\gA$ and $\gA'$ glue to a desired $G$-torsor isomorphism $(X|_{A_0})_A \isomto X$ if and only if their restrictions to $A[\f 1{rr'}]$ agree. The difference of these restrictions is given by an element $g \in G_0^A(R[\f 1{rr'}])$, and our flexibility of adjusting the choices of $\gA$ and $\gA'$ amounts to the fact that $g$ only matters through its class in the double coset 
\[
\tst G_0^A(R[\f1r])\backslash G_0^A(R[\f 1{rr'}]) /G_0^A(R[\f 1{r'}]).
\]
By \Cref{prop:Zariski-patch}, this double coset is trivial, so we may adjust $\gA$ and $\gA'$ to ensure that they agree over $A[\f 1{rr'}]$, as desired.
\epf

\brems
\remi
\Cref{prop:Zariski-patch} and \Cref{cor:H1Zar-G0,cor:Quillen-patching} hold for any group sheaf $G$ that commutes both with filtered direct limits of rings and with fiber products $A_1 \times_{A_0} A_2$ in which one of the ring homomorphisms $A_i \ra A_0$ is surjective: indeed, \Cref{lem:axiom-Q} holds (with the same proof) for sheaves $X$ satisfying these properties and  the arguments then continue to work.

\remi
Quillen patching fails beyond affine bases. For instance, the universal extension of $\sO$ by $\sO(-2)$ on $\bP^1_\bC$ is a vector bundle of rank $2$ on $\bA^1_{\bP^1_\bC}$ that does not descend~to~$\bP^1_\bC$ in spite of the fact that it does descend Zariski locally on $\bP^1_\bC$ (because every vector bundle on $\bA^2_\bC$ is trivial). 

\erems


We obtain the following consequence for descending reductive group schemes defined over $\bA^d_R$. 

\bcor
For a ring $R$, a reductive group scheme $H$ over $R[t_1, \dotsc, t_d]$ descends to a reductive group scheme over $R$ if and only if it does so Zariski locally on $R$\uscolon moreover, the same holds with $R[t_1, \dotsc, t_d]$ replaced by any $\bZ^{\oplus d}_{\ge 0}$-graded $R$-algebra $A \cong \bigoplus_{i_1, \dotsc, i_d \ge 0} A_{i_1, \dotsc, i_d}$ such that $R \isomto A_{0, \dotsc, 0}$.
\ecor

\bpf
We focus on the `if,' since the converse is obvious. The type of the geometric fibers of a reductive group scheme is locally constant on the base (see \cite{SGA3IIInew}*{Expos\'e XXII, Proposition~2.8}) and $\Spec(A)$ has connected $R$-fibers (compare with the proof of \Cref{prop:Zariski-patch}), so we may replace $R$ by a direct factor to assume that this type is constant for $H$. We let $\bbH$ be the split reductive group over $R$ of the same type as $H$ and use the same references as in the proof of \Cref{cor:lift-reductive}~\ref{LR-c} to argue that $H$ corresponds to a torsor under the $R$-group scheme $\underline{\Aut}_\gp(\bbH)$. The claim then follows from \Cref{cor:Quillen-patching} applied with $G = \underline{\Aut}_\gp(\bbH)$.
\epf

The following consequence of Quillen patching reduces \Cref{conj:Guo} to the local case.


\bcor \label{cor:local-to-global}
Let $R$ be a ring and let $G$ be a group functor as in Proposition \uref{prop:Zariski-patch} that is locally of finite presentation \cite{SP}*{Definition \href{https://stacks.math.columbia.edu/tag/049J}{049J}} \up{for instance, $G$ could be an $R$-group algebraic space locally of finite presentation}. Every $G$-torsor over $\bA^1_R$ that is trivial away from an $R$-finite closed subscheme is trivial as soon as the same holds with $R$ replaced by $R_\fm$  for each maximal ideal $\fm \subset R$.
\ecor

\bpf
Let $X$ be a $G$-torsor over $\bA^1_R$ that is trivial away from an $R$-finite closed subscheme $Z \subset \bA^1_R$. In the category of sets, filtered direct limits commute with finite inverse limits, so, by descent, $X$ is also locally of finite presentation. Thus, the assumption about the $R_\fm$ implies that $X$ is trivial Zariski locally on $R$. \Cref{cor:Quillen-patching}~\ref{QP-b} then ensures that $X$ descends to a $G$-torsor $X_1$ over $R$, which is simply the pullback of $X$ along the section $t \mapsto 1$, and we need to show that $X_1$ is trivial. 

Since $X$ is trivial away from $Z$, we may glue it with the trivial $G$-torsor over $\bP^1_R \setminus Z$, and so extend $X$ to a $G$-torsor $\ov{X}$ over $\bP^1_R$. By the assumption on the $R_\fm$ again, the restrictions $X|_{\bA^1_{R_\fm}}$ are all trivial, so the restriction $\ov{X}|_{\bP^1_{R} \setminus \{ 0 \}}$ becomes trivial away from the section at infinity after base change to each $R_\fm$. The assumption on the $R_\fm$, this time applied to the $G$-torsor $\ov{X}|_{\bP^1_{R} \setminus \{ 0 \}}$ (note that $\bP^1_{R} \setminus \{ 0 \}$ is isomorphic to $\bA^1_R$), now implies that each $\ov{X}|_{\bP^1_{R_\fm} \setminus \{ 0 \}}$ is trivial. Thus, by local finite presentation as before, $\ov{X}|_{\bP^1_{R} \setminus \{ 0 \}}$ is trivial Zariski locally on $R$. \Cref{cor:Quillen-patching}~\ref{QP-b} now ensures that $\ov{X}|_{\bP^1_{R} \setminus \{ 0 \}}$ descends to a $G$-torsor over $R$. By pulling back along the section $t \mapsto 1$, we find that this descended $G$-torsor is $X_1$ while, on the other hand, by pulling back along the infinity section, we find that it is trivial. In conclusion, $X_1$ is a trivial $G$-torsor, as desired. 
\epf

The following ``inverse'' to Quillen patching is more elementary but is also useful. 
Its case when $G = \GL_n$ and $A = R[t_1, \dotsc, t_d]$ is due to Roitman \cite{Roi79}*{Proposition 2}.

\bprop \label{prop:inverse-patching}
Let $R$ be a ring, let $G$ be a quasi-affine, flat, finitely presented $R$-group scheme, let $A \cong \bigoplus_{i_1,\dotsc, i_d \ge 0} A_{i_1,\dotsc, i_d}$ be a $\bZ_{\ge 0}^{\oplus d}$-graded $R$-algebra such that $R \isomto A_{0, \dotsc, 0}$ \up{for instance, $A$ could be $R[t_1, \dotsc, t_d]$}, and suppose that every $G$-torsor over $A$ \up{resp.,~whose pullback to $A_{0, \dotsc, 0} \cong R$ is trivial} descends to a $G$-torsor over $R$. Then, for any multiplicative subset $S \subset R$, every $G$-torsor over $A_S$ whose restriction to each local ring of $(A_{0, \dotsc, 0})_S \cong R_S$ extends to a $G$-torsor over $R$ \up{resp.,~whose restriction to $(A_{0, \dotsc, 0})_S \cong R_S$ is Zariski locally trivial} descends to a $G$-torsor over $R_S$. 
\eprop

\bpf
We let $X$ be a $G$-torsor over $A_S$ as in the statement that we wish to descend to $R_S$, and we use \Cref{cor:Quillen-patching} (with a limit argument) to enlarge $S$ and reduce to the case when $R_S$ is local. Then, by our assumption, the restriction of $X$ to $(A_{0, \dotsc, 0})_S \cong R_S$ extends to a (resp.,~trivial) $G$-torsor $X_0$ over $R$. Granted this, we use a limit argument to reduce to the case when $S$ is a singleton $\{r \}$  at the cost of $R_S$ no longer being local. Consider the projection map
\[
\tst R \oplus \p{\bigoplus_{ (i_1, \dotsc, i_d) \neq (0, \dotsc, 0)} A_{i_1, \dotsc, i_d}[\f1r]} \cong A[\f1r] \times_{R[\f1r]} R \surjects R,
\]
which, by the snake lemma, induces an isomorphism both modulo $r^n$ and also on $r^n$-torsion for every $n > 0$. By patching of \Cref{lem:patch-nonflat}, we may use this map to glue up a $G$-torsor $\wt{X}$ over $A[\f1r] \times_{R[\f1r]} R$ from the $G$-torsor $X$ over $A[\f1r]$ and the $G$-torsor $X_0$ over $R$. By construction, the base change of $\wt{X}$ to $A[\f1r]$ is $X$, so it suffices to descend $\wt{X}$ to a $G$-torsor over $R$. However, 
\[
\tst A[\f1r] \times_{R[\f1r]} R \cong \varinjlim A
\]
where the direct limit is indexed by $\bN$ and its transition maps $A \ra A$ are given by multiplication by $r^{i_1 + \dotsc + i_d}$ on the degree $(i_1, \dotsc, i_d)$ piece $A_{i_1, \dotsc, i_d}$. A limit argument then shows that $\wt{X}$ descends to a $G$-torsor over some copy of $A$ in this direct limit, and hence, by the assumption on $A$, even descends further to a $G$-torsor over $R$, as desired.
\epf


\csub[The unramified case of the Bass--Quillen conjecture] \label{sec:unramified-BQ}

We wish to illustrate the utility of Quillen patching by reviewing the proof of the unramified case of the Bass--Quillen conjecture, in which it plays a central role. The starting point of the proof is the following $\dim( R) \le 1$ case, which is susceptible to an inductive argument thanks to Quillen patching. 

\bthm[Quillen, Suslin; \cite{Qui76}*{Theorem~4$'$}] \label{thm:Quillen-Suslin}
For a regular ring $R$ of dimension $\le 1$, every vector bundle $\sV$ on $\bA^d_R$ descends to $R$, in particular, $\sV$ is free if $R$ is a principal ideal domain.
\ethm

\bpf
We follow \emph{loc.~cit.} The last assertion follows from the rest and from the structure theorem for finitely generated modules over a principal ideal domain \cite{SP}*{Lemma \href{https://stacks.math.columbia.edu/tag/0ASV}{0ASV}}. Moreover, by Quillen patching, we may  assume that our $R$ is a principal ideal domain, and we will induct on $d$. The key insight for attacking $d > 0$ is the observation that the localization $R(T)$ of $R[T]$ with respect to the multiplicative set of \emph{monic} polynomials is again a principal ideal domain: indeed, $R(T)$ is a regular ring in which every prime is of height $\le 1$ (any prime of higher height would have to lie over a maximal ideal of $R$, but the closed $R$-fibers of $R(T)$ are fields), and it is a unique factorization domain by \cite{SP}*{Lemmas \href{https://stacks.math.columbia.edu/tag/0BC1}{0BC1} or \href{https://stacks.math.columbia.edu/tag/0AFT}{0AFT}}. 
This and the inductive hypothesis show that the finite projective $R[T_1, \dotsc, T_d]$-module that corresponds to $\sV$ becomes free over $R(T_1)[T_2, \dotsc, T_d]$. Thus,  formal glueing, that is, \Cref{lem:patch}, ensures that $\sV$ extends to a vector bundle $\wt{\sV}$ on $\bP^1_R \times_R \bA^{d - 1}_R$. Horrocks \Cref{prop:Horrocks} 
then implies that $\sV$ descends to vector bundle on the second factor $\bA^{d - 1}_R$. Consequently, the inductive hypothesis applies and shows that $\sV$ descends to $R$, as desired.
\epf

\bthm[Quillen, Suslin, Lindel, Popescu] \label{thm:BQ}
For a regular ring $R$ whose localizations at maximal ideals are \emph{unramified} regular local rings, every vector bundle on $\bA^d_R$ descends to $R$. 
\ethm

\bpf
The assumption implies that the localization of $R$ at any prime ideal is a regular local ring. Thus, we may apply \S\ref{pp:basic-BQ}~\ref{BBQ-1} to reduce to $d = 1$ and then apply Quillen patching, namely, \Cref{cor:Quillen-patching}, to also assume that $R$ is an unramified regular local ring. The Popescu \Cref{thm:Popescu} and a limit argument then reduce to $R$ being a local ring of a scheme $X$ that is smooth 
 over a ring $k$ that is either a field or some $\bZ_{(p)}$ (this was one of the first successes of Popescu's theorem!). 

By shrinking $X$, we may assume that the vector bundle in question is defined over all of $\bA^1_X$ and, by specializing if needed, we may assume that $R = \sO_{X,\, x}$ for a \emph{closed} point $x \in X$. By Lindel's \Cref{cor:baby-Lindel} and \Cref{lem:Lindel}, then there are
\bitem
\m
a local ring $R_0$ of an affine space over $k$;

\m
a local ring homomorphism $R_0 \ra R$ and a nonunit $r \in R_0$ such that $R_0/r R_0 \isomto R/r R$.
\eitem
\addtocounter{footnote}{-7}
\renewcommand{\thefootnote}{\fnsymbol{footnote}}
By induction on $\dim( R)$, the base case $\dim (R) \le 1$ being \Cref{thm:Quillen-Suslin}, we may assume that our vector bundle $\sV$ on $\bA^1_R$ trivializes over $\bA_{R[\f1r]}^1$: indeed, by the inductive assumption (with \Cref{cor:Quillen-patching} to pass to the local rings of $R[\f 1r]$), the restriction $\sV|_{\bA^1_{R[\f1r]}}$ descends to a vector bundle on $R[\f1r]$, and this descent is trivial because it extends to a vector bundle on $R$ given by the restriction of $\sV$ to the origin of $\bA^1_R$.\footnote{\emph{Added after publication.} This sentence contains a small gap (that seems to be inherited from the literature): as written, the inductive hypothesis does not apply to the local rings of $R[\f1r]$ because they are not local rings of $X$ at \emph{closed} points. To fix this, assume instead that $R = \sO_{X,\, x}$ for a point $x \in X$ that is \emph{nongeneric} in its $k$-fiber of $X$ (so we no longer assume that $x$ is a closed point) and then note that, for the same argument to work, it suffices to instead use the following generalization of  \Cref{cor:baby-Lindel}: for a smooth scheme $\sX$ over a ring $\kappa$ and a point $y \in \sX$ that is nongeneric in its $\kappa$-fiber of $\sX$, there are an affine open $U \subset \sX$ containing $y$ and an \'{e}tale map $f \colon U \ra \bA^{d}_\kappa$ that induces an isomorphism $k_{f(y)} \isomto k_y$ on the residue fields as indicated. To argue this generalization, since mapping $U$ to $\bA^d_\kappa$ amounts to giving $d$ global sections of $U$, the fibral criterion of flatness \cite{EGAIV3}*{Th\'{e}or\`{e}me~11.3.10} and spreading out reduce us to the case when $\kappa$ is a field (compare with the proof of \Cref{cor:baby-Lindel}). If this field $\kappa$ is infinite, then the Geometric Presentation Theorem~\ref{thm:geometric-presentation} applied with $Z \ce \ov{\{y\}} \subset \sX$ supplies the desired $U$ and $f$. For a general field $\kappa$, the variant of the Geometric Presentation Theorem given by \cite{split-unramified}*{Proposition 4.1 (with Remark~4.3 to ensure that codimension $\ge 1$ suffices there when $\cO$ is a field)} supplies an affine open $U' \subset \sX$ containing $y$ and a smooth map $g \colon U' \ra \bA^{d - 1}_k$ of relative dimension $1$ such that $y$ is a nongeneric point in its fiber of $g$. The same technique based on the fibral criterion of flatness now allows us to pass to this fiber of $g$ to reduce to the case when $\sX$ is a curve and $y$ is its closed point. In this case, however, \Cref{cor:baby-Lindel} supplies the desired $U$ and $f$. } Thus, by formal glueing of \Cref{lem:patch} applied to the square
\addtocounter{footnote}{6}
\[
\xymatrix{
\bA^1_{R[\f1r]} \ar[d] \ar@{^(->}[r] & \bA^1_R \ar[d] \\
\bA^1_{R_0[\f1r]} \ar@{^(->}[r] & \bA^1_{R_0},
}
\]
the vector bundle $\sV$ descends to a vector bundle on $\bA^1_{R_0}$. In effect, we may replace $R$ by $R_0$ to reduce to $R$ being a local ring of an affine space over $k$. The inverse patching, namely, \Cref{prop:inverse-patching}, then reduces us further to when $R$ is a polynomial algebra over $k$. In conclusion, we are left with showing that, for any $d \ge 0$, every vector bundle on $\bA^d_k$ is free, which follows from \Cref{thm:Quillen-Suslin}.
\epf

\brem
More generally, \Cref{thm:BQ} holds with the same argument when each local ring of $R$ at a maximal ideal is merely flat, with geometrically regular fibers over some Dedekind ring~$k$.
\erem

\brem
In the ramified case of the Bass--Quillen conjecture, it is difficult to envision any reduction to vector bundles on $\bA^d_\bZ$ via formal glueing: we recall from \cite{EGAIV2}*{Proposition 6.1.5} that a quasi-finite morphism between regular schemes of the same dimension is flat, so a formal glueing square between regular schemes of the same dimension always involves a morphism that is \'{e}tale at the points of interest. This automatic \'{e}taleness is a major obstacle hindering any passage from a non-smooth regular $\bZ$-scheme to a smooth one via formal glueing.
\erem




\csub[The analysis of torsors over $\bA^1_R$] \label{sec:A1-analysis}

In \cite{Fed21a} and \cite{Fed21b}, Fedorov developed a technique for analyzing torsors over $\bA^1_R$ via the geometry of the affine Grassmannian. This simplified prior approaches to the Grothendieck--Serre conjecture by eliminating the need for an initial reduction to semisimple, simply-connected groups (for which it is simpler to analyze torsors over $\bA^1_R$, as we already saw in \S\ref{sec:Guo}, see, especially, \Cref{lem:trivial-Gr}). We review his ideas in this section, in particular, we show that they continue to work beyond the equicharacteristic setting. The main statement is the useful in practice \Cref{prop:A1-improve} below. 


The geometric input about affine Grassmannians that is relevant for the study of torsors over $\bA^1_R$ is the surjectivity of the map $\Gr_{(G^\der)^{\mathrm{sc}}} \ra \Gr_G^0$ on field-valued points. After reviewing basic definitions and setup in \S\S\ref{pp:aff-Gr}--\ref{pp:Schubert}, we follow an argument suggested by Timo Richarz to establish this surjectivity in \Cref{cor:mult-stable}. For context, it is helpful to recall that from \cite{Zhu17b}*{Theorem~1.3.11~(3)} that if $G$ is semisimple and the degree of the isogeny $G^{\mathrm{sc}} \ra G$ is invertible on the base, then even $\Gr_{(G^\der)^{\mathrm{sc}}} \isomto \Gr_G^0$. Thus, we are grappling with a ``bad characteristics'' phenomenon, knowing from \cite{HLR20} that the geometry 
of the affine Grassmannian $\Gr_G$ in such characteristics is delicate.


\bpp[The affine Grassmannian] \label{pp:aff-Gr}
For a reductive group $G$ over a field $k$, the \emph{affine Grassmannian} $\Gr_G$ is the functor that to a $k$-algebra $R$ associates the set of isomorphism classes of pairs $(E, \tau)$ consisting of a $G$-torsor $E$ over $R\llb t \rrb$ and its trivialization $\tau \colon E_{R\llp t \rrp} \isomto G_{R\llp t \rrp}$ over $R\llp t \rrp$. By, for instance, \cite{Zhu17b}*{Theorem 1.2.2}, the functor $\Gr_G$ is representable by an ind-projective ind-scheme. 

Concretely, consider the loop and the positive loop groups of $G$ defined as the respective functors
\[
LG \colon R \mapsto G(R \llp t \rrp) \qxq{and} L^+G \colon R \mapsto G(R \llb t \rrb),
\]
which are representable by a group ind-affine ind-scheme (resp.,~by an affine group scheme) over $k$. The subfunctor of $\Gr_G$ that parametrizes those pairs in which $E$ is trivial is the presheaf quotient 
\be \label{eqn:quot-inc}
LG/L^+G \subset \Gr_G.
\ee
A general $E$ trivializes over $R'\llb t \rrb$ for a faithfully flat, \'{e}tale $R$-algebra $R'$ (see \Cref{prop:Hens-pair}~\ref{HP-a} below), so this inclusion exhibits $\Gr_G$ as the \'{e}tale sheafification of $LG/L^+G$. Whenever no nontrivial $G$-torsor over $R\llb t \rrb$ trivializes over $R\llp t \rrp$, the inclusion \eqref{eqn:quot-inc} induces an equality on $R$-points: 
\[
\Gr_G(R) \cong G(R \llp t \rrp)/G(R\llb t \rrb);
\]
this happens, for instance, for a field $R$ (see \S\ref{pp:basic-GS}~\ref{pp:BGS-2}), or for any $R$ when $G$ is either a torus or a pure inner form of $\GL_n$ (combine \Cref{prop:Hens-pair}~\ref{HP-a} with the formula in \S\ref{pp:basic-Guo}~\ref{BG-2} or with \Cref{prop:GLn-noniso} (with \eqref{eqn:change-origin})). In general, $L^+G$ acts on $\Gr_G$ by left multiplication, and $\Gr_G$ is the increasing union of $L^+G$-invariant projective subschemes (for this one fixes an embedding $G \hra \GL_n$ and uses the resulting closed immersion $\Gr_G \hra \Gr_{\GL_n}$, see \cite{Zhu17b}*{proof of Theorem~1.2.2}).

The scheme $L^+G$ is connected, see \cite{CLNS18}*{Chapter 3, Proposition 4.1.1}. 
By \cite{PR08}*{Theorem~5.1}, the map $LG \ra \Gr_G$ induces a bijection on sets of geometric connected components, these components are all clopen, and, if $G$ is semisimple and simply connected, then both $LG$ and $\Gr_G$ are geometrically connected. In general, the \emph{neutral components}, that is, the connected components $LG^0 \subset LG$ and $\Gr_G^0 \subset \Gr_G$ containing the class of the identity, are geometrically connected (as is any connected $k$-scheme $X$ with $X(k) \neq \emptyset$, see \cite{EGAIV2}*{Proposition 4.5.13}). Since $L^+G$ is geometrically connected, its left multiplication action on $LG$ and $\Gr_G$ respects connected components. The map
\be \label{eqn:Gr-cc}
\Gr_{(G^\der)^{\mathrm{sc}}} \ra \Gr_G^0
\ee
is surjective on topological spaces, in fact, it is even surjective on $K$-points for every algebraically closed field extension $K$ of $k$.\footnote{We justify the assertion about $K$-points as follows. Since $LG \ra \Gr_G$ is surjective on $K$-points and a bijection on sets of connected components, by \cite{PR08}*{Theorem 5.1 and the end of the proof of Lemma 17 on page~198 (with $G(L)_1$ defined after Remark 2 on page 189)} (their $G(L)_1$ is our $(LG)^0(K)$), we may replace $G$ by a $z$-extension (see \Cref{theo-resolution-of-reductive-groups}) to reduce to $G^\der$ being simply connected. For such $G$, however, the surjectivity of 
\[
\Gr_{G^\der}(K) \ra \Gr_G^0(K)
\]
follows from \cite{PR08}*{last line on page 197 and proof of Lemma~5 on page 191} (by the latter, $T(L)_1$ there is our~$T(K\llb t \rrb)$).} By \cite{Zhu17b}*{Theorem 1.3.11~(3)}, if $G$ is semisimple with $G^{\mathrm{sc}} \ra G$ of degree prime to $\Char k$, then the map \eqref{eqn:Gr-cc} is even an isomorphism.
\epp

\bpp[Schubert cells] \label{pp:Schubert}
With $G$ over $k$ as in \S\ref{pp:aff-Gr}, let $T \subset G$ be a maximal $k$-torus with 
\[
X_*(T) \ce \underline{\Hom}_\gp(\bG_m, T).
\]
 By \cite{SGA3Inew}*{Expos\'e VI$_{\x{\upshape{A}}}$, Th\'{e}or\`{e}me 3.3.2}, the $L^+G$-orbit of any $x \in \Gr_G(k)$ is a smooth $k$-subscheme of $\Gr_G$. When $x$ is the image of $t$ under  the base change to $k\llp t \rrp$ of the $k$-morphism given by a $\lambda \in X_*(T)(k)$, the resulting subscheme is the \emph{Schubert cell} 
\[
\Gr_G^\lambda \subset \Gr_G.
\]
Its closure (schematic image) in $\Gr_G$ is the \emph{Schubert variety} 
\[
\Gr_G^{\le \lambda} \subset \Gr_G,
\]
which is a reduced, projective $k$-scheme containing $\Gr_G^\lambda$ as a dense open. In the case when $T$ is split, the $\Gr_G^\lambda$ topologically exhaust $\Gr_G$: then, by \cite{PR08}*{Appendix, Proposition 8}, every field-valued (equivalently, (algebraically closed field)-valued) point of $\Gr_G$ factors through some $\Gr_G^\lambda$. In general, the same holds for the $k$-subschemes 
\[
\tst \Gr_G^{[\lambda]} \ce \bigcup_{\lambda' \in \Gal(k^\sep/k) \cdot \lambda} \Gr_G^{\lambda'} \subset \Gr_G \qxq{with} \lambda \in X_*(T)(k^\sep).
\]
Thus, letting $T^{\mathrm{sc}} \subset (G^\der)^{\mathrm{sc}}$ be the maximal torus induced by $T \subset G$, we see from \eqref{eqn:Gr-cc} that the $\Gr_G^{[\lambda]}$ with $\lambda \in X_*(T^{\mathrm{sc}})(k^\sep) \subset X_*(T)(k^\sep)$ topologically exhaust the neutral component~$\Gr_G^0$. 
\epp

We now argue that these $k$-subschemes $\Gr_G^{[\lambda]}$ are insensitive to replacing $G$ by $(G^\der)^{\mathrm{sc}}$. 

\bprop \label{prop:GrG-sc}
For a reductive group $G$ over a field $k$, a maximal $k$-torus $T \subset G$, the resulting maximal $k$-torus $T^{\mathrm{sc}} \subset (G^\der)^{\mathrm{sc}}$, and a $\lambda \in X_*(T^{\mathrm{sc}})(k^\sep)$, 
the $k$-morphism 
\[
\Gr_{(G^\der)^{\mathrm{sc}}}^{[\lambda]} \isomto \Gr_G^{[\lambda]} \qxq{induced by} \Gr_{(G^\der)^{\mathrm{sc}}} \ra \Gr_G \qx{is an isomorphism.}
\]
\eprop

\bpf
The argument is similar to that of \cite{Fed21a}*{Proposition 2.8} and was suggested to us by Timo Richarz. The claim is insensitive to enlarging $k$, so we reduce to $k$ being algebraically closed and then, by passing to individual Schubert cells, to showing that 
\[
\Gr^{\lambda}_{(G^{\der})^{\mathrm{sc}}} \isomto \Gr_G^\lambda.
\]
This last isomorphism, however, is a special case of \cite{HR21}*{Lemma 3.8}.
\epf

We turn to the promised conclusion about the behavior of $\Gr_{(G^\der)^{\mathrm{sc}}} \ra \Gr_G$ on field-valued points. 

\bcor \label{cor:mult-stable}
For a reductive group $G$ over a field $k$, 
the following map is surjective on $k$-points\ucolon
\[
\Gr_{(G^\der)^{\mathrm{sc}}} \xra{\eqref{eqn:Gr-cc}} \Gr_G^0,
\]
in particular, the image of the following map is stable under left multiplication by $G(k\llb t \rrb)$\ucolon
\[
\Gr_{(G^\der)^{\mathrm{sc}}}(k) \overset{\S\ref{pp:aff-Gr}}{\cong} (G^\der)^{\mathrm{sc}}(k\llp t \rrp)/(G^\der)^{\mathrm{sc}}(k\llb t \rrb) \ra G(k\llp t \rrp)/G(k\llb t \rrb) \overset{\S\ref{pp:aff-Gr}}{\cong} \Gr_G(k).
\]
\ecor

\bpf

By \S\ref{pp:aff-Gr}, the ind-scheme $\Gr_{(G^\der)^{\mathrm{sc}}}$ is connected, so the map $\Gr_{(G^\der)^{\mathrm{sc}}} \ra \Gr_G$ factors through the clopen $\Gr_G^0 \subset \Gr_G$. Moreover, by \S\ref{pp:Schubert}, a $k$-point of $\Gr_G^0$ factors through some $\Gr_G^{[\lambda]}$ for a $\lambda \in X_*(T^{\mathrm{sc}})(k^\sep)$, where $T \subset G$ is a maximal torus and $T^{\mathrm{sc}} \subset (G^\der)^{\mathrm{sc}}$ is the corresponding maximal torus of $(G^\der)^{\mathrm{sc}}$. Thus, by \Cref{prop:GrG-sc}, every such point lifts to $\Gr_{(G^\der)^{\mathrm{sc}}}$, as claimed.

By \S\ref{pp:aff-Gr}, the source of the left multiplication map $L^+G \times_k \Gr^0_G \ra \Gr_G$ is connected, so this map factors through $\Gr_G^0 \subset \Gr_G$. Consequently, its image  on $k$-points is $\Gr_G^0(k)$, which, by the above, agrees with the image of $\Gr_{(G^\der)^{\mathrm{sc}}}(k) \ra \Gr_G(k)$. In particular, the latter is $G(k\llb t \rrb)$-stable.
\epf

Before turning to the consequence for torsors over $\bA^1_R$ in \Cref{prop:A1-improve}, we record the following lemma, which clarifies one of the hypotheses appearing there and is a minor generalization of \cite{Fed21a}*{Proposition 2.3}. By this lemma, the hypothesis in question may be arranged by pulling back the torsor under study along the map $\bP^1_R \xra{t\, \mapsto \, t^d} \bP^1_R$ for any sufficiently divisible $d$.

\blem \label{lem:d-power}
For a field $k$, a semisimple $k$-group $G$, opens $U, U' \subset \bP^1_k$, and a generically trivial $G$-torsor $E$ over $U$, the pullback of $E$ along any finite $k$-morphism $U' \ra U$ of degree divisible by the degree of the isogeny $G^{\mathrm{sc}} \ra G$ \up{or merely by the exponent of the quotient $X_*(T)/X_*(T^{\mathrm{sc}})$ for a maximal split $k$-torus $T^{\mathrm{sc}} \subset G^{\mathrm{sc}}$ with image $T \subset G$} lifts to a Zariski locally trivial $G^{\mathrm{sc}}$-torsor~over~$U'$. 
\elem

\bpf
The kernel of the isogeny $T^{\mathrm{sc}} \ra T$ is a subgroup of the kernel of $G^{\mathrm{sc}} \ra G$, so the degree $d_T$ of the former divides that of the latter. Since $d_T$ is simply the order of $X_*(T)/X_*(T^{\mathrm{sc}})$, it is divisible by the exponent $e_T$ of this quotient. Thus, the parenthetical assertion is indeed more general, and we need to show the claim under the assumption that $e_T$ divides the degree $d$ of the finite $k$-morphism $U' \ra U$. For this, we first note that, by \S\ref{pp:basic-GS}~\ref{pp:BGS-2}, the $G$-torsor $E$ over $U$ is Zariski locally trivial.

The key input to the proof is \cite{Gil02}*{Corollaire 3.10 (a)}, according to which $E$ is the extension of $\sO(1)|_U$ (viewed as a $\bG_m$-torsor) along some cocharacter $\mu \colon \bG_m \ra T$. The pullback of $\sO(1)|_U$ to $U'$ is $\sO(d)|_{U'}$, so the pullback of $E$ to $U'$ is the extension of $\sO(1)_{U'}$ along the cocharacter $d\mu \colon \bG_m \ra T$. However, the assumption $e_T\mid d$ ensures that $d$ kills $X_*(T)/X_*(T^{\mathrm{sc}})$, so $d\mu$ factors through a cocharacter $\bG_m \ra T^{\mathrm{sc}}$. Consequently, the pullback of $E$ to $U'$ lifts to a $G^{\mathrm{sc}}$-torsor over $U'$ that comes from a $\bG_m$-torsor over $U'$, and hence is Zariski locally trivial, as desired. 
\epf

We are ready to present the following sharpening of the core result of \cite{Fed21a}; the latter refined \cite{FP15}*{Theorem 3}, which was the centerpiece technical novelty of \emph{op.~cit}. It may be viewed as a Horrocks-type statement, namely, it is in the spirit of extending \Cref{conj:Guo} beyond totally isotropic groups. However, the price of allowing anisotropic groups is a weaker conclusion: instead of the $G$-torsor being trivial over all of $\bA^1_R$ as in \Cref{conj:Guo}, one only concludes that it is trivial away from a fixed $R$-(finite \'{e}tale) closed subscheme along which $G$ is sufficiently isotropic.

\bprop \label{prop:A1-improve}
Let $R$ be a semilocal ring, 
let $G$ be a reductive $R$-group, write the canonical decomposition \eqref{eqn:ad-decomposition} of $G^\ad$ as
\[
\tst G^\ad \cong \prod_i G_i \qxq{with} G_i \ce \Res_{\wt{R}_i/R}(\wt{G}_i),
\] 
where $\wt{R}_i$ \up{resp.,~$\wt{G}_i$} is a finite \'{e}tale $R$-algebra \up{resp.,~an adjoint $\wt{R}_i$-group with simple geometric fibers}, and let $Y_i \subset Y \subset \bA^1_R$ be nonempty  closed subschemes such that
\benumr
\m
$Y$ and each $Y_i$ are all finite \'{e}tale over $R$\uscolon

\m\label{A1I-ii}
\addtocounter{footnote}{-8}
\renewcommand{\thefootnote}{\fnsymbol{footnote}}
$(G_i)_{Y_i}$ is totally isotropic\footnote{\emph{Added after publication.} The assumption \ref{A1I-ii} should be strengthened to: $\mathbb G_{m,\, Y_i} \subset (G_i)_{Y_i}$ for every $i$. This is needed in order to be able to apply \Cref{rem:trivial-Gr} in the proof below.} for every $i$\uscolon
\addtocounter{footnote}{7}

\m \label{A1I-i}
$\sO(1)$ is trivial on $\bP^1_R \setminus Y$\uscolon and

\m\label{A1I-iii}
$\sO(1)$ is trivial on $\bP^1_{k_\fm} \setminus (Y_i)_{k_\fm}$ for every $i$ and every maximal ideal $\fm \subset R$ with $(G_i)_{k_\fm}$~isotropic.
\eenum
For  a $G$-torsor $E$ over $\bP^1_R$ that is trivial away from an $R$-finite closed subscheme $Z \subset \bA^1_R \setminus Y$, if for every maximal ideal $\fm \subset R$ the $G^\ad$-torsor over $\bP^1_{k_\fm}$ induced by $E_{\bP^1_{k_\fm}}$ lifts to a \emph{generically trivial} $(G^\ad)^{\mathrm{sc}}$-torsor over $\bP^1_{k_\fm}$ \up{see Lemma \uref{lem:d-power}}, then $E_{\bP^1_R \setminus Y}$ is a trivial $G$-torsor over $\bP^1_R \setminus Y$.
\eprop

The assumptions become simpler when $G$ is totally isotropic, for instance, quasi-split: then \ref{A1I-ii} is automatic and one may choose $Y_i = Y$ to make \ref{A1I-iii} follow from \ref{A1I-i}. Since \Cref{conj:Guo} remains open in general, \Cref{prop:A1-improve} is useful even for totally isotropic $G$. Regardless of what $G$ is, it is typically straight-forward to arrange \ref{A1I-i} by making sure that either $Y$ contains an $R$-point of $\bP^1_R$ or that $Y$ contains both an $R$-(finite \'{e}tale) closed subscheme of degree $n$ and one of degree $n + 1$ for some $n > 0$ (this ensures that both $\sO(n)$ and $\sO(n + 1)$ are trivial on $\bP^1_R \setminus Y$, so that so is $\sO(1)$).

\bpf
The statement is mild generalization of \cite{Fed21a}*{Theorem 6} and the proof is similar, even if we present it slightly differently. 
It combines the techniques reviewed in \S\ref{sec:Guo} with the analysis of the geometry of the affine Grassmannian that we carried out in the beginning of this section.

By the Cayley--Hamilton theorem, the $R$-(finite \'{e}tale) closed subscheme $Y \subset \bA^1_R$ is cut out by a monic polynomial with coefficients in $R$ (see \cite{split-unramified}*{Remark 6.4}). Thus, the coordinate ring of the formal completion of $\bA^1_R$ along $Y$ is $R'\llb t \rrb$ for a finite \'{e}tale $R$-algebra $R'$ that is the coordinate ring of $Y$ (where $t$ is a monic polynomial in the coordinate of $\bA^1_R$). Likewise, for each $i$, the formal completion of $\bA^1_R$ along $Y_i$ is, compatibly, $R_i\llb t \rrb$ for a finite \'{e}tale $R$-algebra $R_i$ that is the coordinate ring of $Y_i$; this $R_i$ is a direct factor of $R'$, so that $R' \cong R_i \times R_i'$. 

We fix a trivialization $\tau \in E(\bP^1_R\setminus Z)$ of $E|_{\bP^1_R \setminus Z}$. Since $Y \subset \bP^1_R \setminus Z$, this $\tau$ trivializes the restriction of $E$ to $R'\llb t \rrb$ and  we use $\tau$ to regard $E$ as the glueing corresponding to $1 \in G(R'\llp t \rrp)/G(R'\llb t \rrb)$ of $E|_{\bP^1_R \setminus Y}$ and the trivial~$G$-torsor over $R'\llb t \rrb$ (see \Cref{lem:patch-nonflat}).

We let $\fm$ range over the maximal ideals of $R$, set $k \ce \prod_\fm k_\fm$, let $E_i$ be the $G_i$-torsor over $\bP^1_{ k}$ induced by $E$, and let $E^{\mathrm{sc}}_i$ be a generically trivial $G^{\mathrm{sc}}_i$-torsor over $\bP^1_{ k}$ that lifts $E_i$ (such an $E^{\mathrm{sc}}_i$ was assumed to exist).  By the semilocal Dedekind case of the Grothendieck--Serre conjecture (see \S\ref{pp:basic-GS}~\ref{pp:BGS-2}), the generic triviality implies that $E_i^{\mathrm{sc}}$ is trivial on a formal neighborhood of $(Y_i)_{k}$ in $\bP^1_{k}$. 
We fix a trivialization $\tau_i$ over such a neighborhood and use it to regard $E_i^{\mathrm{sc}}$ as the glueing corresponding to 
\[
\tst 1 \in G_i^{\mathrm{sc}}((R_i \tensor_R k)\llp t \rrp)/G^{\mathrm{sc}}_i((R_i \tensor_R  k)\llb t \rrb)
\]
of $E_i^{\mathrm{sc}}|_{\bP^1_{k} \setminus (Y_i)_{ k}}$ and the trivial $(G^{\mathrm{sc}}_i)_{(R_i \tensor_R k)\llb t \rrb}$-torsor (see \Cref{lem:patch-nonflat}).

Of course, the trivializations $\tau$ and $\tau_i$ need not be compatible, in other words, using $\tau$ as the reference, the image of $\tau_i$ in  $G_i((R_i \tensor_R  k)\llp t \rrp)$ need not be the identity. Nevertheless, this image of $\tau_i$ and that of $\tau$ both describe the same $G_i$-torsor over $\bP^1_{ k}$ (the one induced by $E$) as the glueing of the same $G^\ad$-torsor $E_i|_{(\bP^1 \setminus Y_i)_{ k}}$ over $(\bP^1 \setminus Y_i)_{ k}$ and the trivial $G_i$-torsor over $(R_i \tensor_R k)\llb t \rrb$. Concretely, this identification of the glueings means that the image of $\tau_i$ lies in 
\[
\tst G_i((R_i \tensor_R k)\llb t \rrb) \subset G_i((R_i \tensor_R  k)\llp t \rrp),
\]
in other words, that the images of $\tau$ and $\tau_i$ are $G_i((R_i \tensor_R  k)\llb t \rrb)$-translates of each other. Thus, \Cref{cor:mult-stable} implies---and this is a crucial point---that, at the cost of $E_i^{\mathrm{sc}}$ only lifting $E_i$ over $\bP^1_k\setminus (Y_i)_k$, we may change the glueings $E_i^{\mathrm{sc}}$ and the trivializations $\tau_i$ to arrange that they be compatible with $\tau$: namely, still with $\tau$ as the reference, that the image of $\tau_i$ be the class of the identity $1$ in 
\[
G_i((R_i \tensor_R  k)\llp t \rrp)/G_i((R_i \tensor_R  k)\llb t \rrb).
\]
By \cite{Gil02}*{Th\'{e}or\'{e}me 3.8~(b)} (with \S\ref{pp:basic-GS}~\ref{pp:BGS-2}), the generic triviality of $E_i^{\mathrm{sc}}$ means that this torsor comes from a torsor under a split subtorus of $(G_i^{\mathrm{sc}})_k$, and hence, thanks to \ref{A1I-iii}, that $E_i^{\mathrm{sc}}|_{\bP^1_{ k} \setminus (Y_i)_{k}}$ is a trivial $G_i^{\mathrm{sc}}$-torsor. In particular, the trivial $G_i^{\mathrm{sc}}$-torsor over $\bP^1_{ k}$ is a glueing of $E_i^{\mathrm{sc}}|_{(\bP^1 \setminus Y_i)_{ k}}$ and the trivial $G_i^{\mathrm{sc}}$-torsor over $(R_i \tensor_R k)\llb t \rrb$ and, continuing to use $\tau_i$ as reference, this glueing is given by an
\[
\tst \gA_i \in G_i^{\mathrm{sc}}((R_i \tensor_R  k)\llp t \rrp)/G_i^{\mathrm{sc}}((R_i \tensor_R  k)\llb t \rrb).
\]
By \ref{A1I-ii} and \Cref{rem:trivial-Gr}, this $\gA_i$ lifts to some $\wt{\gA}_i \in G_i^{\mathrm{sc}}(R_i\llp t \rrp)$. We  consider $\wt{\gA}_i$ as an element of $G_i^{\mathrm{sc}}(R'\llp t \rrp)$ by letting it be the identity on the complementary factor $G_i^{\mathrm{sc}}(R_i'\llp t \rrp)$. 

Jointly, the $\wt{\gA}_i$ assemble to an element $\wt{\gA} \in (G^\ad)^{\mathrm{sc}}(R'\llp t \rrp)$. The map $(G^\ad)^{\mathrm{sc}} \ra G^\ad$ factors through the isogeny $G^\der \ra G^\ad$, where $G^\der \subset G$ is the derived subgroup, so $\wt{\gA}$ maps to an element of $G(R'\llp t \rrp)$. With $\tau$ as the reference trivialization, this image of $\wt{\gA}$ in $G(R'\llp t \rrp)$ gives rise to a $G$-torsor $\wt{E}$ over $\bP^1_R$ that is the glueing of $E_{\bP^1_R\setminus Y}$ and the trivial $G$-torsor over $R'\llb t \rrb$. The $G^\ad$-torsor $\ov{E}$ over $\bP^1_{ k}$ induced by $\wt{E}$ is the analogous glueing over $\bP^1_{ k}$ that arises from the image of $\prod_i \wt{\gA}_i$ in 
\[
\tst \prod_i G_i((R_i \tensor_R k)\llp t \rrp).
\] 
Thus, by construction and by the prearranged compatibility of $\tau$ and $\tau_i$, this $\ov{E}$ is a trivial torsor. 

\Cref{lem:Fedorov} now implies that $\wt{E}$ induces a $G^\ad$-torsor over $\bP^1_R$ that is the pullback of a $G^\ad$-torsor over $R$. Thus, since $\wt{E}|_{\bP^1_R \setminus (Y \cup Z)}$ is trivial and since the infinity section factors through $\bP^1_R \setminus (Y \cup Z)$, we conclude that $\wt{E}$ induces a trivial $G^\ad$-torsor over $\bP^1_R$, to the effect that $\wt{E}$ comes from a $Z_G$-torsor $F$ over $\bP^1_R$. It now suffices to argue that $F|_{\bP^1_R \setminus Y}$ is the pullback of a $Z_G$-torsor over $R$: then 
\[
\wt{E}|_{\bP^1_R \setminus Y} \cong E|_{\bP^1_R \setminus Y}
\]
will be the pullback of a $G$-torsor over $R$, so, by considering pullbacks at $\infty$, it will be trivial. 

For showing that $F|_{\bP^1_R \setminus Y}$ descends to a $Z_G$-torsor over $R$, we twist to assume that the pullback of $F$ along the infinity section is trivial, and we then fix a trivialization of this pullback. With this rigidification in place, \cite{MFK94}*{Proposition 6.1} (applied to the morphism $\bP^1_R \ra Z_G$, where $(Z_G)_{\bP^1_R}$ is viewed as the automorphism functor of $F$) ensures that $F$ has no nontrivial automorphisms. We now consider the line bundle $\sO(1)$ on $\bP^1_R$, rigidify it by trivializing its pullback along the infinity section, and use \ref{A1I-i} to reduce to showing that there is a unique cocharacter 
\[
\mu \colon \bG_{m,\, R} \ra Z_G
\]
such that $F$ is isomorphic to the extension along $\mu$ of $\sO(1)$ regarded as a $\bG_{m}$-torsor. By what we already observed, such an isomorphism is unique granted that we require it to be compatible with rigidifications at infinity, so the claim is \'{e}tale local on $R$. Thus, we may assume that the multiplicative $R$-group $Z_G$ is split and reduce to when $Z_G$ is either $\bG_{m}$ or $\mu_{n}$. In the first case, the uniqueness of $\mu$ follows from the classification of line bundles on $\bP^1_R$ that results from \Cref{lem:Fedorov} and \cite{BLR90}*{Section 9.1, Proposition 2}. In the second case, since $\Pic(\bP^1_R)$ is torsion free and $R^\times \isomto \Gamma(R, \sO_{\bP^1_R}^\times)$, our $F$ descends to a $\mu_{n}$-torsor over $R$ that, by checking at infinity, is necessarily trivial, and the unique choice $\mu = 0$ works. 
\epf



\section{Techniques for equating reductive groups and their torsors} \label{sec:fix-G}

For a ring $A$ and an ideal $I \subset A$, we wish to discuss when two reductive $A$-group schemes or two torsors that are isomorphic over $A/I$ are also isomorphic over $A$ (or can be made so by replacing $A$ by a cover that has a section over $A/I$). A simple such setup is a local ring and its maximal ideal, but there are others: for instance, in arguing cases of the Grothendieck--Serre conjecture one arrives at such a setup with $A$ being a smooth curve over a semilocal regular ring $R$ such that $R\isomto A/I$ (so that $I$ cuts out an $R$-point of $A$), and one needs to replace $A$ by a \emph{finite} \'{e}tale cover to equate two reductive $A$-group schemes while preserving the $R$-point $A/I$ over which they are already equal. 

As we review in \S\ref{sec:invariance}, a basic case in which positive answers to such questions are available is when the pair $(A, I)$ is Henselian (see \cite{SP}*{Section \href{https://stacks.math.columbia.edu/tag/09XD}{09XD}} for the definition and basic properties of Henselian pairs, which we use freely; an instructive elementary example is a complete local ring $A$ equipped with its maximal ideal $I$). For a general pair $(A, I)$, this means that the answer becomes positive upon replacing $A$ by an \'{e}tale $A$-algebra $A'$ with $A/I \isomto A'/IA'$. This does not suffice for more delicate applications, for instance, for arguing cases of the Grothendieck--Serre conjecture: there one needs $A'$ to be \emph{finite} \'{e}tale over $A$ at the expense of the induced map $A/I \ra A'/IA'$ merely admitting a section instead of being an isomorphism. Arranging such finiteness tends to be delicate and to involve working with compactifications and using the Bertini theorem, we discuss it in \S\S\ref{sec:conj-compactify}--\ref{sec:compactify}.



\csub[Invariance under Henselian pairs for isomorphism classes of reductive groups] \label{sec:invariance}

 We wish to show in \Cref{cor:lift-reductive} that reductive group schemes lift uniquely across Henselian pairs. This  generalizes \cite{SGA3IIInew}*{Expos\'e XXIV, Proposition 1.21}, which treated the Henselian local case. The argument is based on the following broadly useful invariance properties of torsors.

\bprop \label{prop:Hens-pair}
Let $(A, I)$ be a Henselian pair and let $G$ be an $A$-group algebraic space.
\benum
\m \label{HP-c}
For a smooth $A$-algebraic space $X$ that is either quasi-separated or a scheme,
\[
X(A) \surjects X(A/I), \qxq{and, if $X$ is constant, then even} X(A) \isomto X(A/I).
\]

\m \label{HP-b}
If $G$ is smooth and quasi-separated, 
then
\[
H^1(A, G) \hra H^1(A/I, G).
\]

\m \label{HP-a}
If $G$ is quasi-affine, flat, and of finite presentation, then
\[
H^1(A, G) \surjects H^1(A/I, G).
\]

\m \label{HP-d}
If $G \simeq H \rtimes \ov{G}$ is a semidirect product of an $A$-group $\ov{G}$ that becomes constant over a finite \'{e}tale cover of $A$ and a smooth, quasi-affine, normal $A$-subgroup $H$, then every $G$-torsor over $A/I$ whose induced $\ov{G}$-torsor is isotrivial lifts to a $G$-torsor~over~$A$.

\eenum
\eprop

\bpf
In the constant case of \ref{HP-c}, every $A$-point (resp.,~every $A/I$-point) of $X$ factor through a quasi-compact open, so we lose no generality by assuming that $X$ is quasi-compact, so that it is a finite union of copies of $\Spec(A)$. For such $X$, by \cite{SP}*{Lemma \href{https://stacks.math.columbia.edu/tag/09XI}{09XI}}, the clopen subsets of $X$ are identified with those of $X_{A/I}$ via base change and, by \cite{SP}*{Lemma \href{https://stacks.math.columbia.edu/tag/09ZL}{09ZL}}, this identification respects the property of mapping isomorphically to $\Spec(A)$ (resp.,~to $\Spec(A/I)$). Since sections of $X$ (resp.,~of $X_{A/I}$) are precisely the clopens with this property, the constant case of \ref{HP-c} follows.

Part \ref{HP-a} and the case of \ref{HP-c} when $X$ is quasi-separated are special cases of \cite{Hitchin-torsors}*{Example~2.1.5, Theorem 2.1.6} (whose key input is Tannaka duality for algebraic stacks supplied by \cite{HR19}*{Corollary~1.5~(ii)} or \cite{BHL17}*{Corollary 1.5}). This case of \ref{HP-c} applied to 
\[
X \ce \Isom_G(E, E')
\]
for $G$-torsors $E$ and $E'$ over $A$  implies \ref{HP-b} (see \S\ref{pp:basic-pointed} and \S\ref{pp:basic-representability}).

For a scheme $X$ in \ref{HP-c}, we only seek the surjectivity and, by passing to an open, we may again assume that $X$ is quasi-compact. We will then reduce further to when $X$ is also quasi-separated, a case covered by the previous paragraph. For this reduction, we use a technique of Gabber that appeared in \cite{Bha16}*{Remark 4.6}. Namely, by \cite{SP}*{Lemma \href{https://stacks.math.columbia.edu/tag/03K0}{03K0}}, there is a filtered direct system of \'{e}tale $X$-schemes $X_i $ that are quasi-compact and quasi-separated, are such that Zariski locally on $X_i$ the structure map $X_i \ra X$ is an open immersion, and are such that
\[
\tst \varinjlim_iX_i(R) \isomto X(R) \qx{for every $A$-algebra $R$.}
\]
In particular, a fixed $A/I$-point of $X$ lifts to an $A/I$-point of some $X_i$. Since $X_i$ inherits  $A$-smoothness from $X$, we may replace $X$ by $X_i$ and achieve the desired reduction to quasi-separated $X$.


We turn to the remaining part \ref{HP-d}, in which the morphism of short exact sequences of pointed sets
\[
\xymatrix{
H^1(A, H) \ar[r] \ar[d]^-{\sim} & H^1(A, G) \ar@{->>}[r] \ar[d] & H^1(A, \ov{G}) \ar[d] \\
H^1(A/I, H) \ar[r] & H^1(A/I, G) \ar@{->>}[r] & H^1(A/I, \ov{G})
}
\]
(see \S\ref{pp:basic-subgroups}) will allow us to replace $G$ by $\ov{G}$ as follows. The semidirect product decomposition ensures the displayed surjectivity of the right horizontal arrows and, by \ref{HP-b} and \ref{HP-a}, the analogue of the left vertical map stays bijective for every smooth, quasi-affine $A$-group, for instance, for every form of $H$ for the fppf topology. Moreover, any inner form of $G$ is an extension of an inner form of $\ov{G}$ by a form of $H$. A diagram chase and the twisting bijections \eqref{eqn:change-origin} then show that a $G$-torsor over $A/I$ lifts to a $G$-torsor if the same holds for its induced $\ov{G}$-torsor. Thus, we have reduced to the case $G = \ov{G}$. 

In the remaining case in which $G$ becomes constant over a finite \'{e}tale cover $A'$ of $A$, we fix an isotrivial $G$-torsor $E$ over $A/I$ that is to be lifted to a $G$-torsor. The isotriviality means that $E$ trivializes over some finite \'{e}tale cover $B$ of $A/I$, and we may take $B$ to even be a finite \'{e}tale cover of $A'/IA'$. Consequently, $E$ is described by a section $g \in G(B \tensor_{A/I} B)$ that satisfies the cocycle condition. We use \cite{SP}*{Lemma \href{https://stacks.math.columbia.edu/tag/09ZL}{09ZL}} to lift $B$ to a finite \'{e}tale cover $A' \ra \wt{B}$, and we apply \ref{HP-c} over $\wt{B} \tensor_A \wt{B}$ and over $\wt{B} \tensor_A \wt{B} \tensor_A \wt{B}$ to lift $g$ to a section $\wt{g} \in G(\wt{B} \tensor_A \wt{B})$ that satisfies the cocycle condition with respect to $A \ra \wt{B}$. This $\wt{g}$ gives rise to the desired $G$-torsor $\wt{E}$ that lifts $E$. 
\epf

\brem
In \ref{HP-d}, every $\ov{G}$-torsor over $A/I$ is isotrivial if $A/I$ is Noetherian and its localizations at prime ideals are geometrically unibranch, see \S\ref{pp:geom-uni-split} and \cite[Expos\'e X, Corollaire~5.14]{SGA3II}. In general, however, nonisotrivial $\ov{G}$-torsors exist even when $\ov{G} = \bZ$, see \cite{Hitchin-torsors}*{Remark 2.1.8}. 
\erem

We are ready for the promised invariance under Henselian pairs for reductive group schemes. 


\bprop \label{cor:lift-reductive}
Let $(A, I)$ be a Henselian pair and let $G$ and $G'$ be reductive $A$-group schemes.
\benum
\m\label{LR-b}
Every $A/I$-group isomorphism $\iota\colon G_{A/I} \isomto G'_{A/I}$ lifts to an $A$-group isomorphism $\wt{\iota}\colon G \isomto G'$. 

\m\label{LR-c}
A reductive $A/I$-group $H$ with $\rad(H)$ isotrivial lifts \up{uniquely, by \ref{LR-b}} to a reductive $A$-group.




\eenum
\eprop


\bpf \hfill
\benum
\m
By \cite{SGA3IIInew}*{Expos\'e XXIV, Corollaire 1.9}, the functor $\underline{\Isom}_\gp(G, G')$ that parametrizes group isomorphisms is a torsor under the automorphism functor $\underline{\Aut}_\gp(G)$. Thus, \cite{SGA3IIInew}*{Expos\'e XXIV, Th\'{e}or\`{e}me 1.3} and \Cref{foot:ind-quasi-affine} ensure that $\underline{\Isom}_\gp(G, G')$ is representable by an ind-quasi-affine, smooth $A$-scheme. In particular, by \Cref{prop:Hens-pair}~\ref{HP-c}, every $A/I$-point $\iota$ of $\underline{\Isom}_\gp(G, G')$ 
lifts to a desired $A$-point $\wt{\iota}$. 

\m
By decomposing into clopens and lifting idempotents via \cite{SP}*{Lemma \href{https://stacks.math.columbia.edu/tag/09XI}{09XI}}, we may assume that the type of the geometric fibers of $H$ is constant (see \S\ref{pp:reductive}). We let $\bbH$ be a split reductive group over $A$ of the same type as $H$, so that $H$ is a form of $\bbH_{A/I}$, and hence corresponds to an $\underline{\Aut}_\gp(\bbH)$-torsor $E$ over $A/I$ (see \S\ref{pp:automorphisms}). Since $\rad(H)$ is isotrivial, so is $E$, see \S\ref{pp:isotriviality}. By \Cref{prop:Hens-pair}~\ref{HP-d} and the structure of $\underline{\Aut}_\gp(\bbH)$ described in \eqref{eqn:semidirect}, this $E$ lifts to an $\underline{\Aut}_\gp(\bbH)$-torsor over $A$ that corresponds to the desired lift of $H$.
\qedhere
\eenum

\epf

\brem 
In \Cref{cor:lift-reductive}~\ref{LR-c}, some condition on $H$ is necessary: it is not true that for every Henselian pair $(A, I)$, every reductive $A/I$-group lifts to a reductive $A$-group. Indeed, this fails already for tori: if it held, then, by considering those pairs in which $A$ is normal (or even in which $A$ is a Henselization of some affine space), we could conclude from \S\ref{pp:geom-uni-split} that every torus over an affine base is isotrivial, contradicting \Cref{rem:nonisotrivial-tori} or \cite{SGA3II}*{Expos\'e~X, Section 1.6}.
\erem



\csub[A conjecture about compactifying reductive groups and consequences for torsors] \label{sec:conj-compactify}

For a ring $A$ and its ideal $I \subset A$, arranging the finer lifting property mentioned in the introduction of this chapter amounts to finding situations in which a functor $\sF$ has the following property:
\be\label{eqn:property-star} \tag{$\star$}\ba 
&\x{for every $x \in \sF(A/I)$, there are a faithfully flat, finite, \'{e}tale $A$-algebra $\wt{A}$,}\\
&\x{an $A/I$-point $a \colon \wt{A} \surjects A/I$, and a $\wt{x} \in \sF(\wt{A})$ whose $a$-pullback is $x$.}
\ea
\ee
If $I$ lies in every maximal ideal of $A$, then the faithful flatness requirement follows from the rest. Granted this further condition on $I$,  the results presented in \S\ref{sec:invariance} arrange the same without $\wt{A}$ being finite over $A$ but with $\wt{A}/I\wt{A} \cong A/I$ instead. In contrast, getting \eqref{eqn:property-star} instead typically requires finer techniques that we discuss in this section. We begin with the following simple example.

\beg \label{eg:locally-constant}
Any faithfully flat, finite, \'{e}tale $A$-scheme has property \eqref{eqn:property-star} because we may choose $\wt{A}$ to be its coordinate ring. Somewhat more interestingly, if $A$ is Noetherian and its local rings are geometrically unibranch (see \S\ref{pp:geom-uni-split}), then any faithfully flat $A$-scheme $X$ that becomes constant \'{e}tale locally on $A$ has property \eqref{eqn:property-star}: by \cite{SGA3II}*{Expos\'e X, Corollaire~5.14} (with \cite{EGAI}*{Corollaire~6.1.9}), these assumptions ensure the connected components of $X$ are clopen subschemes that are finite \'{e}tale over $A$, so, by considering a sufficiently large union of them, we reduce to the finite \'{e}tale case. 
\eeg

In practice, the key source of property \eqref{eqn:property-star} is following lemma that we learned from the argument of \cite{OP01}*{Lemma 7.2} and that was also pointed out by Uriya First.

\blem \label{lem:quasi-section}
For a semilocal ring $A$ and a projective, finitely presented $A$-scheme $X$, any $A$-smooth open $U \subset X$ that is dense in the closed $A$-fibers of $X$, is of pure relative dimension $d \ge 0$, and is faithfully flat over $A$ has property \eqref{eqn:property-star} with respect to any ideal $I \subset A$.
\elem

\bpf
Fix a $x$ as in \eqref{eqn:property-star}. By replacing $U$ by a finite union of some of its open affines (that cover the image of $x$), we may assume that $U$ is quasi-compact. Then a limit argument allows us to assume (mostly for comfort) that $A$ is Noetherian. By decomposing into connected components, we may also assume that $\Spec(A)$ is connected. 
Finally, we fix a projective embedding $X \hra  \bP^n_A$ and postcompose it with a linear change of projective coordinates if necessary to arrange that $x$ is the origin $A/I$-point 
\[
[0 : \dotsc : 0 : 1] \in \bP^n_A(A/I).
\]
Let $S \subset \Spec(A)$ be the union of the closed points. Since $(X\setminus U)_S$ is of dimension less than $d$, we may apply \Cref{lem:Bertini}, with $Z = Z_0$ there being the image of $x_S$, to find hypersurfaces $H_1, \dotsc, H_d \subset X_S$ of large enough and constant on $S$ degrees such that $H_1 \cap \dotsc \cap H_d$ lies in $U_S$, is finite \'{e}tale over $S$, is fiberwise nonempty, and contains $x_S$. Granted that these degrees are sufficiently large, \cite{EGAIII1}*{Corollaire~2.2.4} allows us to lift the $H_i$ to hypersurfaces $H_1', \dotsc, H_d' \subset X_{S'}$ where the closed subscheme $S' \subset \Spec(A)$ is the union of $\Spec(A/I)$ and $S$. We may choose these lifts in such a way that they contain $x$: indeed, $x = [0 : \dotsc : 0 : 1]$, so ensuring that $x \in H'_i(A/I)$ amounts to lifting a defining equation of $H_i$ in such a way that the coefficient of the monomial that is a power of the last variable stays zero. Once such $H_i'$ of large  degrees are fixed, we apply \cite{EGAIII1}*{Corollaire~2.2.4} again to lift them to hypersurfaces $\wt{H}_1, \dotsc, \wt{H}_d \subset X$, which, by construction, contain $x$. 

By construction, the scheme-theoretic intersection $\wt{Z} \ce \wt{H}_1 \cap \dotsc \cap \wt{H}_d$ lies in $U$ and contains $x$. By the openness of the quasi-finite locus \cite{SP}*{Lemma \href{https://stacks.math.columbia.edu/tag/01TI}{01TI}} and the finiteness of proper, quasi-finite morphisms \cite{SP}*{Lemma \href{https://stacks.math.columbia.edu/tag/02OG}{02OG}}, the $A$-scheme $\wt{Z}$ is finite. By \cite{EGAIV3}*{Th\'eor\`eme 11.3.8}, it is $A$-flat at its closed points, so the openness of the flat locus \cite{EGAIV3}*{Th\'eor\`eme 11.3.1} ensures that it is $A$-flat. Thus, we check over $S$ that $\wt{Z}$ is faithfully flat and \'{e}tale over $A$. In conclusion, $\wt{Z} = \Spec(\wt{A})$ for a faithfully flat, finite, \'{e}tale $A$-algebra $\wt{A}$ that is equipped with an $A/I$-point $a \colon \wt{A} \surjects A/I$ that corresponds to $x$. The inclusion $\wt{Z} \subset U$ is the desired $\wt{A}$-point $\wt{x} \in U(\wt{A})$ whose $a$-pullback is $x$. 
\epf

The lemma above reveals that checking property \eqref{eqn:property-star} in practice hinges on compactifying a smooth scheme $U$ in question in such a way that it be \emph{fiberwise} dense in its projective compactification $X$. Such density is straight-forward to arrange over a field, basically, because any quasi-projective variety is dense in its closure in a projective space, but the question becomes significantly more delicate over a general base, for instance, in mixed characteristic. In this regard, it would be useful to resolve the following conjecture about compactifying reductive group schemes. In \cite{split-unramified}, we bypassed it by taking advantage of the quasi-splitness assumption.


\bconj \label{conj:compactify}
For an isotrivial reductive group $G$ over a Noetherian scheme $S$, there are a projective, finitely presented $S$-scheme $\ov{G}$ equipped with a left $G$-action and a $G$-equivariant $S$-fiberwise dense open immersion 
\[
G \hra \ov{G}.
\]
\econj

We restricted to Noetherian $S$ for concreteness, although it is plausible that a sufficiently natural argument may work more generally. It may be worthwhile to also require that $\ov{G}$ be  equipped with a commuting \emph{right} $G$-action and then that the open immersion $G \hra \ov{G}$ be equivariant with respect to both actions. In \S\ref{sec:compactify}, we check that even this finer variant holds in the case of a torus. We now show that the conjecture implies that any isotrivial $G$-torsor over a semilocal base has property \eqref{eqn:property-star}.

\bprop \label{prop:torsors-star}
For an isotrivial reductive group $G$ over a Noetherian scheme $S$, if Conjecture~\uref{conj:compactify} holds for $G$, then any isotrivial $G$-torsor $E$ admits an $S$-fiberwise dense open immersion 
\[
E \hra \ov{E}
\]
into a projective, finitely presented $S$-scheme $\ov{E}$\uscolon in particular, if, in addition, $S = \Spec(A)$ for a semilocal ring $A$, then any isotrivial $G$-torsor has property \eqref{eqn:property-star} with respect to any ideal $I \subset A$. 
\eprop

\bpf
The final assertion follows from the rest and from \Cref{lem:quasi-section}. For the rest, we let $\ov{G}$ be the compactification of $G$ supplied by \Cref{conj:compactify} and consider the contracted product
\[
\ov{E} \ce  E \times_G \ov{G}.
\]
By general results on quotients, $\ov{E}$ is an algebraic space (see \S\ref{pp:basic-representability}). Moreover, by construction, it comes equipped with an open immersion $E \hra \ov{E}$ that is \'{e}tale locally on $S$ isomorphic to $G \hra \ov{G}$. Thus, all that remains for us to check is that this $\ov{E}$ is a projective $S$-scheme. 

For this, we will only use a finite \'{e}tale cover $S' \surjects S$ such that $\ov{E}_{S'}$ is a projective $S'$-scheme, for instance, this may be a finite \'{e}tale cover trivializing $E$ (see \S\ref{pp:isotriviality}). Consider the restriction of scalars 
\[
\ov{E}' \ce \Res_{S'/S}(\ov{E}_{S'}).
\]
Its base change to a larger finite \'{e}tale cover of $S$ decomposes as a product of copies of $\ov{E}$, so, by \cite{CGP15}*{Proposition A.5.8 and its proof} and \cite{EGAII}*{Remarques 5.5.4~(i)}, this $\ov{E}'$ is a projective $S$-scheme. By checking \'{e}tale locally on $S$, the adjunction morphism $\ov{E} \hra \ov{E}'$ is a closed immersion (compare with \cite{CGP15}*{Proposition A.5.7}), so $\ov{E}$ is a projective $S$-scheme, as promised. 
\epf

We conclude the section with a consequence of \Cref{conj:compactify} for equating reductive group schemes. This consequence is used in proving cases of the Grothendieck--Serre conjecture \ref{conj:Grothendieck-Serre}, and it is close in spirit to \cite{OP01}*{Proposition~7.1}, \cite{PSV15}*{Proposition~5.1}, or \cite{Pan20b}*{Theorem~3.4}. In \cite{split-unramified}, we used the quasi-splitness assumption to avoid tackling \Cref{conj:compactify}, see \cite{split-unramified}*{Lemma~5.1}.


\bprop \label{prop:finer-invariance}
For  a Noetherian semilocal ring $A$ whose local rings are geometrically unibranch, an ideal $I \subset A$, reductive $A$-groups $G$ and $G'$ that on geometric $A$-fibers have the same type and such that Conjecture~\uref{conj:compactify} holds for $G^\ad$, and an $A/I$-group isomorphism 
\[
\iota \colon G_{A/I} \isomto G'_{A/I},
\]
there are a faithfully flat, finite, \'{e}tale $A$-algebra $\wt{A}$ equipped with an $A/I$-point $a \colon \wt{A} \surjects A/I$ and an $\wt{A}$-group isomorphism $\wt{\iota} \colon G_{\wt{A}} \isomto G'_{\wt{A}}$ whose $a$-pullback is $\iota$. 
\eprop

It is key that $A \ra \wt{A}$ be finite: without this, the assertion would be a special case of \Cref{cor:lift-reductive}.

\bpf
Similarly to the proof of \Cref{cor:lift-reductive}, we consider the smooth, ind-quasi-affine scheme 
\[
X \ce \underline{\Isom}_\gp(G, G'),
\]
and need to show that $X$ has property \eqref{eqn:property-star} with respect to the ideal $I \subset A$. The condition on the geometric $A$-fibers ensures that $G$ and $G'$ are isomorphic \'{e}tale locally on $A$ (see \S\ref{pp:reductive}). Thus,  by \S\ref{pp:automorphisms} (see also \S\ref{pp:basic-representability}), the adjoint group $G^\ad$ acts freely on $X$ by conjugation and the quotient 
\[
\ov{X} \ce X/G^\ad
\]
is a faithfully flat $A$-scheme that becomes constant \'{e}tale locally on $A$. By \Cref{eg:locally-constant}, this $\ov{X}$ has property \eqref{eqn:property-star}, so there are a faithfully flat, finite, \'{e}tale $A$-algebra $A'$ equipped with an $A/I$-point $a' \colon A' \surjects A/I$ and an $A'$-point $\iota' \in \ov{X}(A')$ whose $a'$-pullback is the $A/I$-point $\ov{\iota} \in \ov{X}(A/I)$ induced by $\iota \in X(A/I)$. In effect, by base changing the $G^\ad$-torsor $X \ra \ov{X}$ along $\iota'$, we reduce to showing that every $G^\ad$-torsor over $A'$ has property \eqref{eqn:property-star} with respect to the ideal $\Ker(a') \subset A'$. However, by \S\ref{pp:isotriviality}, the adjoint group $G^\ad$ is isotrivial on every connected component of $\Spec(A)$ and every $G^\ad$-torsor over $A'$ is also isotrivial. Consequently, our assumption about \Cref{conj:compactify} allows us to apply \Cref{prop:torsors-star} to conclude that the $\iota'$-pullback of $X \ra \ov{X}$ has property \eqref{eqn:property-star}, as desired.
\epf


\csub[Compactifying torsors under tori] \label{sec:compactify}

We show that \Cref{conj:compactify} holds in the case when $G$ is a torus. For this, we adapt the work of Colliot-Th\'{e}l\`{e}ne--Harari--Skorobogatov \cite{CTHS05}*{Corollaire 1} (so also previous work of Brylinski and K\"{u}nnemann), who used toric techniques to build the required compactification over a field.





\bthm \label{thm:torsor-compactify}
For an isotrivial torus $T$ over a Noetherian scheme $S$,  there are a projective, smooth $S$-scheme $\ov{T}$, commuting left and right $T$-actions on $\ov{T}$, and an $S$-fiberwise dense open~immersion 
\[
\iota\colon T \hra \ov{T}
\]
that is equivariant with respect to both the left and the right translation actions
 of $T$.
\ethm


\bpf[Proof of Theorem \uref{thm:torsor-compactify}]
We show how to construct the desired $\iota\colon T \hra \ov{T}$ by using the results of \cite{CTHS05}, where $\iota$ was constructed when the base is a field by using the theory of toric varieties. We decompose $S$ into connected components to assume that it is connected and let $S'$ be a finite \'{e}tale cover of $S$ splitting $T$.  
We may assume that $S'$ is connected and then enlarge it to ensure that it is Galois over $S$ with group $\Gamma$. We claim that it suffices to construct an analogous equivariant compactification $\iota' \colon T_{S'} \hra \ov{T}'$ over $S'$ granted that $\ov{T}'$ is equipped with a $\Gamma$-action (compatibly with the $\Gamma$-action on $S'$, so that the action will be free on $\ov{T}'$ because it so already on $S'$) that commutes with the left and right actions of $T_{S'}$  and $\iota'$ is $\Gamma$-equivariant. Indeed, by \cite{SP}*{Lemma \href{https://stacks.math.columbia.edu/tag/07S7}{07S7}}, the projectivity of $\ov{T}'$ will then ensure that the quotient $\ov{T} \ce \ov{T}'/\Gamma$ is an $S$-scheme. Moreover, by \cite{SP}*{Lemma \href{https://stacks.math.columbia.edu/tag/0BD0}{0BD0} (with Lemmas \href{https://stacks.math.columbia.edu/tag/0BD2}{0BD2}, \href{https://stacks.math.columbia.edu/tag/0AH6}{0AH6}, and \href{https://stacks.math.columbia.edu/tag/05B5}{05B5})}, this $\ov{T}$ will automatically be projective and smooth over $S$. Thus, we will be able to choose $\iota$ to be
\[
T \cong T_{S'}/\Gamma \hra \ov{T}'/\Gamma = \ov{T}.
\] 
For the remaining construction of $\ov{T}'$, we will use \cite{CTHS05}*{Th\'{e}or\`{e}me~1} and the theory of toric varieties, and we begin by noting that, by functoriality, $\Gamma$ acts on the cocharacter lattice $L \ce X_*(T_{S'})$, as well as on $L_\bR \ce L \tensor_\bZ \bR$. Let $\cF$ be a (rational, polyhedral) fan in $L_\bR$ whose associated toric variety is $\bP^{\, \rk(L)}$ (see, for instance, \cite{Dan78}*{Example 5.3}). This fan need not be $\Gamma$-invariant but, by \cite{CTHS05}*{Th\'{e}or\`{e}me~1}, there is a (rational, polyhedral) fan $\cF'$ in $L_\bR$ that is $\Gamma$-invariant, is a subdivision of $\cF$, and is projective and smooth in the sense that its associated toric variety is projective and smooth (these properties can be expressed combinatorially in terms of $\cF'$, see \cite{Dan78}*{Section 3.3} and \cite{CTHS05}*{Proposition~1}). The construction \cite{Dan78}*{Section~5.2} that builds the toric variety associated to $\cF'$ adapts to any base, so we obtain a flat, finitely presented $S'$-scheme $\ov{T}'$ equipped with commuting left and right $T_{S'}$-actions, a compatible $\Gamma$-action, and an $S'$-fiberwise dense, $T_{S'}$-biequivariant and compatibly $\Gamma$-equivariant open immersion $\iota'\colon T_{S'} \hra \ov{T}'$ over $S'$. By \cite{Dan78}*{Section 3.3} applied $S'$-fiberwise, $\ov{T}'$ is $S'$-smooth, so it remains to argue that it is projective over $S'$.

Due to its combinatorial definition, the $S'$-scheme $\ov{T}'$ descends to a scheme over $\Spec(\bZ)$, so \cite{Dan78}*{Proposition 5.5.6} and its proof, which is based on the finer than usual form \cite{EGAII}*{Corollaire~7.3.10~(ii)} of the valuative criterion of properness, imply that $\ov{T}'$ is proper over $S'$. In combinatorial terms, the fact that $\cF'$ is projective means that there exists a function $h \colon L_\bR \ra \bR$ that is \emph{strictly upper convex} in the sense that, letting $\cF'^{\, \mathrm{ top}} \subset \cF'$ denote the subset of top-dimensional cones, there are linear forms 
\[
\{ \ell_\sigma\}_{\sigma \in \cF'^{\, \mathrm{top}}} \subset \Hom_\bZ(L, \bZ) = X^*(T)
\]
satisfying $\ell_\sigma(x) \ge h(x)$ for all $x \in L_\bR$ with equality if and only if $x \in \sigma$ (see \cite{CTHS05}*{Proposition~1} and \cite{Oda88}*{Lemma~2.12}). This last requirement uniquely determines the characters $\ell_\sigma$  because each $\sigma$ is top-dimensional. Thus, as in \cite{Oda88}*{Proposition 2.1~(i) and its proof}, the function $h$, more precisely, the $\ell_\sigma$, define a line bundle $\sL_h$ on $\ov{T}'$. By \cite{EGAIV3}*{Corollaire 9.6.4}, checking that $\sL_h$ is ample over $S'$ can be done $S'$-fiberwise. Consequently, \cite{Oda88}*{Proposition 2.1~(vi), Corollary~2.14, and their proofs} imply the $S'$-ampleness of $\sL_h$, and hence the $S'$-projectivity of $\ov{T}'$. 
\epf

As an immediate consequence, \Cref{prop:torsors-star} holds in the case when $G$ is a torus. More explicitly, we obtain the following statement about torsors under tori.


\bcor\label{cor:property-star}
For a Noetherian semilocal ring $A$ whose local rings are geometrically unibranch, an ideal $I \subset A$, an $A$-torus $T$, a $T$-torsor $E$, and an $e \in E(A/I)$, there are a faithfully flat, finite, \'{e}tale $A$-algebra $\wt{A}$ equipped with an $A/I$-point $a \colon \wt{A} \surjects A/I$ and an $\wt{e} \in E(\wt{A})$ whose $a$-pullback is $e$\uscolon in particular, $E$ trivializes over some \emph{finite} \'{e}tale cover of $A$ \up{choose $I = A$}. 
\ecor

If we do not require $\wt{A}$ to be \emph{finite} over $A$, then the claim follows from \Cref{prop:Hens-pair}~\ref{HP-c}.  


\bpf
We included the last aspect of the claim for the sake of emphasis: as we already saw in \S\ref{pp:isotriviality}, the geometric unibranchedness assumption ensures that both $T$ and $E$ are isotrivial. Thanks to this isotriviality, the main assertion follows by combining \Cref{thm:torsor-compactify} with \Cref{prop:torsors-star}.
\epf



\begin{appendix}

\section{Resolutions of reductive groups} \label{appendix}

\center{Yifei Zhao\footnote{CNRS, Universit\'{e} Paris-Saclay,   Laboratoire de math\'{e}matiques d'Orsay, F-91405, Orsay, France. E-mail address: \url{yifei.zhao@universite-paris-saclay.fr}}}

\justify

This appendix is an exposition on the construction of flasque and coflasque resolutions of a reductive group $G$ over a general base scheme $S$, subject only to the condition that $\rad(G)$ be isotrivial.

The notions of flasque and coflasque tori are due to Colliot-Th\'el\`ene--Sansuc \cite{CTS87}. The existence of a coflasque resolution strengthens that of a $z$-extension of Langlands and Kottwitz \cite[Section~1]{Kot86}, which is often stated over a field of characteristic zero (\cite{DMOS82}*{Chapter V, Section 3}, \cite{BK00}, for example). When the base is a field of arbitrary characteristic, both resolutions are constructed by Colliot-Th\'el\`ene in \cite{CT04}, but, as observed by Gonz\'alez-Avil\'es \cite{GA13}, the same proof yields the existence of flasque resolutions over locally Noetherian, geometrically unibranch schemes (e.g.~a normal scheme). Our proof follows Colliot-Th\'el\`ene's strategy, but we replace the hypotheses on $S$ by the isotriviality of $\rad(G)$, which holds whenever $S$ is locally Noetherian and geometrically unibranch but could remain valid in other contexts.

The \S\S\ref{sec-tori-flasque-and-coflasque}--\ref{sec-central-isogenies-recall} are preparatory and the main construction appears as Theorem \ref{theo-resolution-of-reductive-groups}. As an application, we explain a simple reduction of the Grothendieck--Serre conjecture \ref{conj:Grothendieck-Serre} to the case when the  derived subgroup is simply connected, see \Cref{prop:GS-app}. The author thanks K.~\v{C}esnavi\v{c}ius for many helpful conversations and comments.


\csub[Group schemes of multiplicative type]


In this section, we review group schemes of multiplicative type. The most important notion for us is the isotriviality of such group schemes.

\bpp
Let $S$ be a scheme. For an fppf sheaf of abelian groups $\sF$ over $S$, one may consider the fppf sheaf $\mathbb D_S(\sF)$ whose value at an affine $S$-scheme $S'$ is $\Hom(\sF_{S'}, \mathbb G_{m, S'})$. Here, $\Hom$ is viewed in the category of fppf sheaves of abelian groups over $S'$. The fppf sheaf $\mathbb D_S(\sF)$ again takes values in abelian groups, the group structure being inherited from $\mathbb G_m$.
\epp

\bpp
\label{sec-multiplicative-type-group-schemes}
An $S$-group scheme $G$ is \emph{diagonalizable} if there exist a finitely generated abelian group $M$ and an isomorphism between $G$ and the group scheme $\mathbb D_S(M_S)$, where $M_S$ denotes the constant sheaf with values in $M$. An $S$-group scheme $G$ is \emph{of multiplicative type} if it is diagonalizable fppf locally on $S$.\footnote{This definition agrees with \cite{Con14}*{Definition B.1.1}, but it differs from \cite{SGA3II}*{Expos\'e IX, D\'{e}finition 1.1}, where $G$ is only required to be fpqc locally isomorphic to $\mathbb D_S(M)$ for an abelian group $M$ (which is not necessarily finitely generated).} In fact, every $S$-group scheme of multiplicative type is diagonalizable \'etale locally on $S$ (\cite{SGA3II}*{Expos\'e X, Corollaire 4.5} or \cite[Proposition B.3.4]{Con14}).

If an $S$-group scheme $G$ of multiplicative type becomes diagonalizable after base change along $\widetilde S\rightarrow S$, then $G$ is said to be \emph{split by $\widetilde S$}. If $S$ is connected, then any $S$-group scheme of multiplicative type is split by some fppf (equivalently, \'etale) surjection $\widetilde S\rightarrow S$ (\cite[Expos\'e IX, Remarque 1.4.1]{SGA3II}). By fppf descent, any $S$-group scheme $G$ of multiplicative type is $S$-affine.

An $S$-group $G$ of multiplicative type is a \emph{torus} if fppf (equivalently, \'etale) locally on $S$  it is of the form $\mathbb D_{S}(M_S)$ for some finitely generated free abelian group $M$.
\epp

\bpp
\label{sec-multiplicative-type-closure-properties}
Groups of multiplicative type enjoy certain closure properties:
\benuma
	\item an $S$-flat, finitely presented closed subgroup of an $S$-group scheme of multiplicative type is again of multiplicative type (\cite{SGA3II}*{Expos\'e X, Corollaire 4.7 b)} or \cite[Corollary B.3.3]{Con14});
	\item a commutative extension of group schemes of multiplicative type is again of multiplicative type (\cite{SGA3II}*{Expos\'e XVII, Proposition 7.1.1} or \cite[Corollary B.4.2]{Con14}).
\eenum
Furthermore, an $S$-group scheme $G$ of multiplicative type is \emph{reflexive} in the sense that the natural transformation $G \ra \mathbb D_S(\mathbb D_S(G))$ is an isomorphism. Indeed, this statement may be verified fppf locally on $S$, where it follows from \cite[Expos\'e VIII, Th\'eor\`eme 1.2]{SGA3II}.

\epp

\bpp \label{pp:def-isotrivial}
An $S$-group scheme $G$ of multiplicative type is called \emph{isotrivial} if $G$ is split by a finite \'etale surjection $\widetilde S\rightarrow S$. When $S$ is connected, an isotrivial $S$-group scheme of multiplicative type is split by a finite connected \'etale Galois cover $\widetilde S\rightarrow S$. We discuss some ways to obtain isotrivial $S$-group schemes of multiplicative type.
\epp

\blem
\label{lem-finite-multiplicative-type-isotrivial}
Let $S$ be a connected scheme. Then any finite $S$-group scheme $G$ of multiplicative type is isotrivial.
\elem
\begin{proof}
Let $\sF$ denote the fppf sheaf of abelian groups $\mathbb D_S(G)$. Since $S$ is connected, there is an fppf surjection $\widetilde S\rightarrow S$ such that $\sF_{\widetilde S}$ is isomorphic to the constant sheaf $M_{\widetilde S}$ for a \emph{finite} abelian group $M$. The descent data of $\sF_{\widetilde S}$ allow us to construct an $\Aut(M)$-torsor $\sP$ over $S$ such that $\sF$ is the fppf sheaf of abelian groups induced from $\sP$. Since $\Aut(M)$ is finite, $\sP$ is representable by a finite \'etale surjection $\widetilde S'\rightarrow S$. In particular, $\sP$ is trivialized by $\widetilde S'$. It follows that $G$ is split by $\widetilde S'$.
\end{proof}

The same argument proves more generally that an $S$-group scheme of multiplicative type whose maximal torus has rank $\le 1$ is isotrivial. This fact can be compared with Lemma \ref{lem-isotriviality-closure-properties}(ii) below.

\blem
\label{lem-isotriviality-closure-properties}
Let $S$ be a connected scheme. Given a short exact sequence of $S$-group schemes of multiplicative type\ucolon
$$
1 \ra G_1 \ra G \ra G_2 \ra 1,
$$
\benumr
	\item if $G$ is diagonalizable \up{resp.~isotrivial}, then both $G_1$ and $G_2$ are diagonalizable \up{resp.~isotrivial}\uscolon
	\item if $G_1$ is isotrivial and $G_2$ is finite, then $G$ is isotrivial.
\eenum
\elem
\begin{proof}
Statement (i) is established in \cite[Expos\'e IX, Proposition 2.11]{SGA3II}. To prove statement (ii), we may assume that both $G_1$ and $G_2$ are diagonalizable by replacing $S$ with a connected finite \'etale cover. Since $\mathbb D_S$ is an anti-equivalence on reflexive fppf sheaves of abelian groups \cite[Expos\'e~VIII, Proposition 1.0.1]{SGA3II}, it restricts to an exact functor on the full subcategory of $S$-group schemes of multiplicative type. In particular, we obtain a short exact sequence of fppf sheaves of abelian groups:
$$
1 \rightarrow M_{2, S} \rightarrow \mathbb D_S(G) \rightarrow M_{1, S} \rightarrow 1.
$$
Here, $M_{i, S}$ (for $i=1,2$) denotes the constant sheaf associated to a finitely generated $\mathbb Z$-module $M_i$. The finiteness of $G_2$ allows us to assume that $M_2$ is finite.

It remains to show that any class in $\Ext^1_{\fppf}(M_{1,S}, M_{2,S})$ comes from $\Ext^1_{\mathbb Z}(M_1, M_2)$ after passing to a finite \'etale cover $\widetilde S\rightarrow S$. For this statement, it suffices to treat the case where $M_1$ is a cyclic group. For $M_1 = \mathbb Z$, we have $\Ext^1_{\fppf}(\mathbb Z_S, M_{2, S}) \cong H^1_{\fppf}(S, M_2) \cong H^1_{\et}(S, M_2)$, and because $M_2$ is finite, any class in $H^1_{\et}(S, M_2)$ vanishes over a finite \'etale surjection $\widetilde S\rightarrow S$. For $M_1 = \mathbb Z/n$ for an integer $n\ge 1$, we have an exact sequence:
$$
\Hom(\mathbb Z, M_2) \xrightarrow{n} \Hom(\mathbb Z, M_2) \ra \Ext^1_{\fppf}((\mathbb Z/n)_S, M_{2,S}) \ra \Ext^1_{\fppf}(\mathbb Z_S, M_{2,S}).
$$
By the same argument as above, any class in $\Ext^1_{\fppf}((\mathbb Z/n)_S, M_{2,S})$ has zero image in $\Ext^1_{\fppf}(\mathbb Z_S, M_{2, S})$ after passing to a connected finite \'etale cover $\widetilde S\rightarrow S$. Equivalently, this means that over $\widetilde S$, the class comes from $\Ext^1_{\mathbb Z}(\mathbb Z/n, M_2)$.
\end{proof}

\brem \label{rem:nonisotrivial-tori}
In Lemma \ref{lem-isotriviality-closure-properties}(ii), the finiteness hypothesis on $G_2$ cannot be dropped. Indeed, whenever $H^1_{\et}(S, \mathbb Z)\neq 0$, there exist self-extensions of $\mathbb G_m$ which are not isotrivial. To see this, we use the isomorphism $\Ext^1_{\fppf}(\mathbb G_m, \mathbb G_m) \cong H^1_{\et}(S, \mathbb Z)$ and the fact that any class of $H^1_{\et}(S, \mathbb Z)$ which vanishes on a finite \'etale cover of $S$ must already be zero (because $H^1_\et(S, \bZ) \hra H^1_\et(S, \bQ)$).
\erem

\bpp
\label{pp:geom-uni-split}
There is a convenient condition on the base scheme $S$ which guarantees that all multiplicative type $S$-group schemes are isotrivial.

A local ring $R$ is \emph{geometrically unibranch} if its strict Henselization $R^\sh$ has a unique minimal prime, see \cite{SP}*{Definition \href{https://stacks.math.columbia.edu/tag/0BPZ}{0BPZ} and Lemma \href{https://stacks.math.columbia.edu/tag/06DM}{06DM}}. A scheme $S$ is \emph{geometrically unibranch} if so are its local rings. For example, a normal scheme (in the usual sense that its local rings are normal domains) is geometrically unibranch. Every connected component of a locally Noetherian, geometrically unibranch scheme is irreducible (see \cite{GW20}*{Exercise~3.16~(a)}).

Let $S$ be a locally Noetherian, geometrically unibranch scheme. By \cite[Expos\'e X, Th\'eor\`eme~5.16]{SGA3II}, every $S$-group scheme $G$ of multiplicative type splits over a finite \'etale surjection $\widetilde S \ra S$. When $S$ is connected, we may further assume that $\widetilde S \ra S$ is a connected Galois cover.
\epp

\csub[Flasque and coflasque tori]
\label{sec-tori-flasque-and-coflasque}
In this section, we focus on isotrivial tori. The study of these objects is equivalent to that of Galois modules with integral coefficients. We discuss several conditions on isotrivial tori (quasi-trivial, flasque, and coflasque) which are ``of Galois cohomology nature''.

\bpp
\label{sec-galois-cover-equivalence}
Suppose that $S$ is a connected scheme and let $\widetilde S \rightarrow S$ be a connected finite \'etale Galois cover. By \cite{SGA3II}*{Expos\'e X, Proposition 1.1}, the construction $\mathbb D_{\widetilde S}$ induces an equivalence of categories between
\benuma
	\item group schemes $G\rightarrow S$ of multiplicative type split by $\widetilde S$;
	\item finitely generated $\mathbb Z$-modules $M$ equipped with a $\Gal(\widetilde S/S)$-action.
\eenum
Under this equivalence, tori $T\rightarrow S$ split by $\widetilde S$ correspond to finitely generated free $\mathbb Z$-modules $M$ equipped with a $\Gal(\widetilde S/S)$-action---such an $M$ is called the \emph{character lattice} of $T$, and we denote it by $\check{\Lambda}_{T,\, \widetilde S}$. Its $\mathbb Z$-linear dual is called the \emph{cocharacter lattice} of $T$, which we denote by $\Lambda_{T,\, \widetilde S}$.
\epp

\bpp
Let $\Gamma$ be a finite group. A $\Gamma$-lattice $\Lambda$, i.e.,~a finitely generated free $\mathbb Z$-module equipped with a $\Gamma$-action, is called \emph{quasi-trivial} if it has a $\Gamma$-stable $\mathbb Z$-basis. Clearly, $\Lambda$ is quasi-trivial if and only if its $\mathbb Z$-linear dual $\check{\Lambda} \ce \Hom_\bZ(\Lambda, \bZ)$, equipped with the contragredient $\Gamma$-action, is quasi-trivial.
\epp

The following lemma describes the conditions that end up defining flasque tori.

\blem
\label{lemm-flasque-tori}
Let $\Gamma$ be a finite group and let $\Lambda$ be a lattice equipped with a $\Gamma$-action. The following conditions are equivalent\ucolon
\benumr
	\item $H^1(\Gamma', \Lambda) = 0$ for any subgroup $\Gamma' \subset \Gamma$;
	\item $\Ext^1_{\mathbb Z[\Gamma]}(P, \Lambda) = 0$ for any quasi-trivial $\Gamma$-lattice $P$.
\eenum
\elem

\begin{proof}
The key observation is as follows. Suppose that $P$ is a quasi-trivial lattice with a basis $X$ that consists of a single $\Gamma$-orbit. Fix an $x\in X$ and let $\Gamma'\subset\Gamma$ be the stabilizer of $x$. Then we have an isomorphism $P \cong \mathbb Z[\Gamma/\Gamma']$ of $\mathbb Z[\Gamma]$-modules, and so an isomorphism $\Ext^1_{\mathbb Z[\Gamma]}(P, \Lambda) \cong H^1(\Gamma', \Lambda)$.
\end{proof}

\bpp
Let $\Gamma$ be a finite group. A $\Gamma$-lattice $\Lambda$ is called
\benuma
	\item \emph{coflasque} if it satisfies the equivalent conditions of Lemma \ref{lemm-flasque-tori}; and
	\item \emph{flasque} if its $\mathbb Z$-linear dual $\check{\Lambda}$ satisfies the equivalent conditions of Lemma \ref{lemm-flasque-tori}.\footnote{We refer to the original paper of Colliot-Th\'el\`ene--Sansuc \cite[Section 0.5]{CTS87} for other equivalent characterizations of flasque lattices, including the one involving Tate cohomology, which often appears in the literature.}
\eenum
By Shapiro lemma, any quasi-trivial $\Gamma$-lattice is both flasque and coflasque.
\epp

\bpp
In the setting of \S\ref{sec-galois-cover-equivalence}, a torus $T\rightarrow S$ split by $\widetilde S$ is called \emph{quasi-trivial} (resp.~\emph{flasque}; resp.,~\emph{coflasque}) if its character lattice $\check{\Lambda}_{T,\, \widetilde S}$ is quasi-trivial (resp.~flasque; resp.,~coflasque). By \cite[Lemma 1.1]{CTS87}, these notions are independent of the choice of the Galois cover $\widetilde S$.

For any scheme $S$, a torus $T\rightarrow S$ is called \emph{quasi-trivial} (resp.~\emph{flasque}; resp.,~\emph{coflasque}) if every connected component of $S$ admits a connected finite \'etale Galois cover $\widetilde S$ such that $T$ is split by $\widetilde S$ and quasi-trivial (resp.~flasque; resp.,~coflasque) with respect to $\widetilde S$ (again, these notions do not depend on $\wt{S}$).  If a torus $T\rightarrow S$ is quasi-trivial (resp.~flasque, coflasque), then so is its base change along any morphism $S'\rightarrow S$ with $S'$ still connected.

Quasi-trivial tori are both flasque and coflasque, and they can be made more explicit as follows.
\epp

\blem
\label{lemm-quasitrivial-torus-explicit}
Let $S$ be a connected scheme. A torus $T\rightarrow S$ is quasi-trivial if and only if it is a finite product of Weil restrictions of $\mathbb G_m$ along finite \'etale surjections $S' \rightarrow S$.
\elem
\begin{proof}
Suppose that $T\rightarrow S$ is quasi-trivial. Let $\widetilde S\rightarrow S$ be a connected finite \'etale Galois cover such that $T$ is split by $\widetilde S$. Without loss of generality, we may assume that $\Lambda_{T,\, \widetilde S}$ has a basis $X$ that consists of a single $\Gal(\widetilde S/S)$-orbit. Then the $\Gal(\widetilde S/S)$-set $X$ gives rise to a finite \'etale cover $S' \rightarrow S$, and, by \S\ref{sec-galois-cover-equivalence}, we have an isomorphism $T \simeq \Res_{S'/S}(\mathbb G_m)$. The converse is analogous.
\end{proof}

Flasque and coflasque tori enjoy the following pleasant splitting property.

\blem
\label{lemm-flasque-coflasque-split}
In the setting of \uS\uref{sec-galois-cover-equivalence}, a short exact sequence of $S$-tori split by $\widetilde S$
$$
1 \rightarrow T_1 \rightarrow T_2 \rightarrow T_3 \rightarrow 1
$$
is split if either of the following conditions holds\ucolon
\benumr
	\item $T_1$ is quasi-trivial and $T_3$ is coflasque\uscolon
	\item $T_1$ is flasque and $T_3$ is quasi-trivial.
\eenum
\elem
\begin{proof}
By considering character lattices, we translate the problem to splitting the exact sequence
\begin{equation}
\label{eq-character-lattice-short-exact-sequence}
0 \rightarrow \check{\Lambda}_{T_3,\, \widetilde S} \rightarrow \check{\Lambda}_{T_2,\, \widetilde S} \rightarrow \check{\Lambda}_{T_1,\, \widetilde S} \rightarrow 0
\end{equation}
of $\Gal(\widetilde S/S)$-lattices. 
Suppose that $T_1$ is quasi-trivial and $T_3$ is coflasque. Then, by definition,
$$
\Ext^1_{\mathbb Z[\Gal(\widetilde S/S)]}(\check{\Lambda}_{T_1,\, \widetilde S}, \check{\Lambda}_{T_3,\, \widetilde S}) = 0,
$$
so \eqref{eq-character-lattice-short-exact-sequence} splits. Suppose that $T_1$ is flasque and $T_3$ is quasi-trivial. Then the dual of \eqref{eq-character-lattice-short-exact-sequence} splits for the same reason, so, by dualizing again, \eqref{eq-character-lattice-short-exact-sequence} splits as well.
\end{proof}

Next, we shall construct ``resolutions'' of $S$-group schemes of multiplicative type split by $\widetilde S$ in terms of flasque and coflasque tori. The following Lemma of Colliot-Th\'el\`ene--Sansuc \cite{CTS87} will be the basis of our construction of resolutions of reductive $S$-group schemes.

\blem
\label{lemm-resolution-multiplicative-type}
In the setting of \uS\uref{sec-galois-cover-equivalence}, let $G\rightarrow S$ be a group scheme of multiplicative type split by $\widetilde S$. There exist $S$-tori $T_1$ and $T_2$ split by $\widetilde S$ that fit into a short exact sequence of $S$-group schemes
\begin{equation}
\label{eq-resolution-multiplicative-type}
1 \rightarrow G \rightarrow T_1 \rightarrow T_2 \rightarrow 1.
\end{equation}
 Furthermore, we may arrange \eqref{eq-resolution-multiplicative-type} so that either of the following conditions is satisfied\ucolon
\benumr
	\item $T_1$ is flasque and $T_2$ is quasi-trivial\uscolon
	\item $T_1$ is quasi-trivial and $T_2$ is coflasque.
\eenum
\elem
\begin{proof}
By \S\ref{sec-galois-cover-equivalence}, the problem translates into one concerning finitely generated $\mathbb Z$-modules equipped with a $\Gal(\widetilde S/S)$-action, which is addressed in \cite[Lemma 0.6]{CTS87}.
\end{proof}

\brem
\label{sec-surjection-quasi-trivial}
In the setting of \S\ref{sec-galois-cover-equivalence}, let $T \rightarrow S$ be a torus split by $\widetilde S$. In the same vein as Lemma \ref{lemm-resolution-multiplicative-type}, \cite[Lemma 0.6]{CTS87} implies the existence of resolutions by tori split by $\widetilde S$:
$$
1 \rightarrow T_1 \rightarrow T_2 \rightarrow T \rightarrow 1,
$$
such that either 
\benumr
\m
$T_1$ is flasque and $T_2$ is quasi-trivial; or 
\m
$T_1$ is quasi-trivial and $T_2$ is coflasque.
\eenum
These resolutions are the flasque, respectively coflasque resolutions of the torus $T$. The main result we shall prove (Theorem \ref{theo-resolution-of-reductive-groups}) can be viewed as its generalization where $T$ is replaced by a reductive $S$-group with isotrivial radical. In its proof, however, we will only need a special case of the result for $T$: the existence of a surjection $P \twoheadrightarrow T$ from a quasi-trivial torus $P$ split by $\widetilde S$.
\erem










\csub[Central isogenies and the simply connected cover]
\label{sec-central-isogenies-recall}

Before proceeding to construct the promised resolutions of reductive groups in \S\ref{sec:resolutions}, we review the notion of central isogenies that plays an important role there. Recall the notion of the center $Z_G$ of a reductive group scheme $G\rightarrow S$ as defined in \S\ref{pp:center}.

\bpp \label{pp:central-isogeny}
For a scheme $S$, a morphism $f \colon G' \rightarrow G$ of reductive $S$-group schemes is called a \emph{central isogeny} if
\benuma
	\item $f$ is finite, flat, and surjective;
	\item $\ker(f)$ lies in the center of $G'$.
\eenum
We only define the notion of central isogenies for reductive $S$-group schemes, as is done in \cite[Expos\'e XXII, D\'efinition 4.2.9]{SGA3IIInew}. One may generalize this notion to other $S$-group schemes, but it may become pathological: for example, the composition of two central isogenies may fail to be central, see \cite[Exercise 3.4.4(ii)]{Con14}. We now show that such phenomena do not occur for reductive group schemes and then we use central isogenies to define the simply connected cover of a semisimple group scheme in \Cref{lemm-simply-connected-cover}. 
\epp

\blem \label{lem:central-isogeny}
Let $f \colon G' \rightarrow G$ be a central isogeny of reductive $S$-group schemes. 
\benumr
\m\label{lemm-central-isogeny-inverse-of-center}
The induced map $Z_{G'} \rightarrow f^{-1}(Z_G)$ is an isomorphism.

\m\label{lemm-central-isogeny-composition}

For any other central isogeny $g \colon G'' \rightarrow G'$ of reductive group schemes over $S$, the composition $f \circ g \colon G'' \rightarrow G$ is also a central isogeny.
\eenum
\elem
\begin{proof}
In \ref{lemm-central-isogeny-inverse-of-center}, the problem is \'etale local on $S$, so we may assume that $G$ contains a split maximal torus $T$, whose inverse image $T' \ce f^{-1}(T)$ is then a split maximal torus of $G'$. Since $f$ is a central isogeny, the induced map on character lattices $\check{\Lambda}_T \rightarrow \check{\Lambda}_{T'}$ restricts to a bijection between the roots of $(G, T)$ and $(G', T')$, see \cite[Example 6.1.9]{Con14}. The result then follows from the characterization of $Z_G$ as the kernel of the adjoint action $T \rightarrow \GL(\Lie(G))$ (see \S\ref{pp:center}), that is, as the intersection of the $\ker(\alpha)$ over all the roots $\alpha \colon T \rightarrow \mathbb G_m$ of $(G, T)$.

In \ref{lemm-central-isogeny-composition}, $f \circ g$ is finite, flat, and surjective, so we need to verify that $\ker(f \circ g) \subset Z_{G''}$. Indeed, we have
$$
\ker(f \circ g) \cong g^{-1}(\ker(f)) \subset g^{-1}( Z_{G'}) \cong Z_{G''},
$$
where the last isomorphism comes from \ref{lemm-central-isogeny-inverse-of-center}.
\end{proof}

\brem
\label{rem-central-isogeny-kernel-multiplicative}
Suppose that $f : G' \rightarrow G$ is a central isogeny of reductive $S$-group schemes. Then $\ker(f)$ is an $S$-group scheme of multiplicative type. Indeed, $Z_{G'}$ is of multiplicative type (see \S\ref{pp:center}) so this assertion follows from the closure property in \S\ref{sec-multiplicative-type-closure-properties}.
\erem




\bprop
\label{lemm-simply-connected-cover}
Let $S$ be a scheme and let $G$ be a semisimple $S$-group scheme. Consider the category of pairs $(G', f)$ consisting of a semisimple $S$-group scheme $G'$ and a central isogeny $f \colon G' \rightarrow G$, with morphisms $(G'_1, f_1) \rightarrow (G'_2, f_2)$ being given by central isogenies $\alpha \colon G_1' \rightarrow G_2'$ such that $ f_1 = f_2 \circ \gA$. This category has an initial object $(G^{\mathrm{sc}}, f)$, the \emph{simply connected cover} of $G$.
\eprop

\begin{proof}
The proof relies on the classification of pinned reductive groups by root data (\cite[Expos\'e XXV, Th\'eor\`eme 1.1]{SGA3IIInew} or \cite[Theorem 6.1.16]{Con14}). The universal property allows us to work \'etale locally on $S$, so we may assume that $G$ is split with respect to a split maximal torus $T\subset G$.

The split maximal torus $T$ allows us to extract the root data $(\Lambda_T, \Phi, \check{\Lambda}_T, \check{\Phi})$. Let $\Lambda_T^r\subset\Lambda_T$ denote the sublattice generated by the coroots $\Phi$. There is a morphism of root data:
\begin{equation}
\label{eq-root-data-universal-covering}
(\Lambda_T^r, \Phi, \check{\Lambda}_T^r, \check{\Phi}) \rightarrow (\Lambda_T, \Phi, \check{\Lambda}_T, \check{\Phi})
\end{equation}
which induces the identity maps on $\Phi$ and $\check{\Phi}$. The root data $(\Lambda_T^r, \Phi, \check{\Lambda}_T^r, \check{\Phi})$ define a pinned reductive $S$-group $G^{\mathrm{sc}}$ with split maximal torus $T^{\mathrm{sc}}$ and \eqref{eq-root-data-universal-covering} comes from a central isogeny $f \colon G^{\mathrm{sc}} \rightarrow G$ compatible with the splitting (i.e., mapping $T^{\mathrm{sc}}$ to $T$), but $f$ is only unique up to conjugation by $(T/Z_G)(S)$ (\cite[Theorem 6.1.16(1)]{Con14}). The pair $(G^{\mathrm{sc}}, f)$, however, is canonically defined thanks to the isomorphism $T^{\mathrm{sc}}/Z_{G^{\mathrm{sc}}} \cong T/Z_G$ induced by $f$. Next, we argue that $(G^{\mathrm{sc}}, f)$ is canonically independent of the choice of the split maximal torus $T\subset G$. Indeed, conjugation defines an isomorphism between $G/N_G(T)$ and the scheme parametrizing maximal tori of $G$ (\cite[Theorem~3.1.6]{Con14}) so the claim follows from the isomorphism $G^{\mathrm{sc}}/N_{G^{\mathrm{sc}}}(T^{\mathrm{sc}}) \cong G/N_G(T)$ induced by $f$.

To show that the pair $(G^{\mathrm{sc}}, f)$ satisfies the universal property of an initial object, we suppose being given another central isogeny $f' \colon G'\rightarrow G$. For a split maximal torus $T\subset G$ as above, we write $T' \subset G'$ for the induced maximal torus. Arguing with root data as above, we find a central isogeny $\alpha_1 \colon G^{\mathrm{sc}} \rightarrow G'$ such that $f$ and $f'\circ \alpha_1$ differ by conjugation by an element of $(T/Z_G)(S)$. The isomorphism $T'/Z_{G'} \cong T/Z_G$ then allows us to construct the unique central isogeny $\alpha \colon G^{\mathrm{sc}} \rightarrow G'$ which satisfies $f = f'\circ\alpha$.
\end{proof}

\brem
\label{sec-simply-connected-cover-pointwise-definition}
Another definition of the simply connected cover is given in \cite[Exercise~6.5.2(i)]{Con14}, which characterizes the central isogeny $f \colon G^{\mathrm{sc}} \rightarrow G$ by the fact that the geometric $S$-fibers of $G^{\mathrm{sc}}$ are simply connected, i.e., they admit no nontrivial central isogenies from semisimple groups. It is easy to see that the two definitions agree. In particular, the formation of the simply connected cover $G^{\mathrm{sc}}$ of $G$ commutes with arbitrary base change $S' \rightarrow S$.
\erem




\csub[Existence of resolutions] \label{sec:resolutions}

In this section, we construct flasque and coflasque resolutions of reductive group schemes with isotrivial radical tori. Recall that to a reductive $S$-group scheme $G$, we have associated several other reductive $S$-group schemes in the main text: the derived subgroup $G^{\der}$, which is semisimple, and the tori $\rad(G)$ and $\corad(G) := G/G^{\der}$ (see \S\ref{pp:center}).


\bthm
\label{theo-resolution-of-reductive-groups}
Let $S$ be a connected scheme and let $G$ be a reductive $S$-group scheme such that $\rad(G)$ is isotrivial. Fix a central isogeny $f \colon G'^{\der} \rightarrow G^{\der}$. There exists a central extension
\begin{equation}
\label{eq-resolution-of-reductive-groups}
1 \rightarrow T_1 \rightarrow G' \rightarrow G \rightarrow 1
\end{equation}
of reductive $S$-group schemes such that $\rad(G')$ is isotrivial and $G' \rightarrow G$ induces $f$ on derived subgroups. Furthermore, setting $T_2 := \corad(G')$, we may arrange \eqref{eq-resolution-of-reductive-groups} so that one of the following conditions is satisfied\ucolon
\benum
	\item $T_1$ is a flasque torus and $T_2$ is a quasi-trivial torus\uscolon
	\item $T_1$ is a quasi-trivial torus and $T_2$ is a coflasque torus.
\eenum
\ethm

\brem
The most typical application of Theorem \ref{theo-resolution-of-reductive-groups} is with $f \colon G'^{\der} \rightarrow G^{\der}$ being the simply connected cover reviewed in Lemma \ref{lemm-simply-connected-cover}. In this case, we obtain a resolution \eqref{eq-resolution-of-reductive-groups} where $G'$ has a simply connected derived subgroup and the tori $T_1$, $T_2$ satisfy the conditions above. These are called \emph{flasque}, respectively \emph{coflasque resolutions} of $G$.

Note that if an $S$-torus $T$ admits a flasque (or coflasque) resolution, then it must be isotrivial (Lemma \ref{lem-isotriviality-closure-properties}(i)). Thus the hypothesis that $\rad(G)$ be isotrivial cannot be dropped.

Finally, we remark that if $S$ is locally Noetherian and geometrically unibranch (such as a normal scheme), then the isotriviality condition on $\rad(G)$ is automatically satisfied (see \S\ref{pp:geom-uni-split}).
\erem


\begin{proof}[Proof of Theorem \uref{theo-resolution-of-reductive-groups}]
By composing the canonical central isogeny $G^{\der} \times \rad(G) \rightarrow G$ of \eqref{eq-reductive-group-derived-radical} with $f \times \id_{\rad(G)}$, we obtain a central isogeny of reductive $S$-group schemes:
\begin{equation}
\label{eq-central-isogeny-radical}
1 \rightarrow H_2 \rightarrow G'^{\der} \times \rad(G) \rightarrow G \rightarrow 1.
\end{equation}
In particular, $H_2$ is a finite $S$-group scheme of multiplicative type (Remark \ref{rem-central-isogeny-kernel-multiplicative}).

Let us denote by $\widetilde S \rightarrow S$ a connected finite \'etale Galois cover which splits $\rad(G)$. By Remark \ref{sec-surjection-quasi-trivial}, we may choose a short exact sequence of tori split by $\widetilde S$:
\begin{equation}
\label{eq-quasi-trivial-resolution-of-radical}
1 \rightarrow H_1 \rightarrow P \rightarrow \rad(G) \rightarrow 1.
\end{equation}
where $P$ is quasi-trivial. Compose the central isogeny \eqref{eq-central-isogeny-radical} with the surjection $P\rightarrow \rad(G)$, we obtain a central extension of reductive $S$-groups:
\begin{equation}
\label{eq-resolution-of-reductive-groups-extension}
1 \rightarrow M \rightarrow G'^{\der} \times P \rightarrow G \rightarrow 1.
\end{equation}
Let us study the commutative $S$-group scheme $M$. By construction, it is an extension of $H_2$ by $H_1$. Since $H_1$ is an isotrivial torus and $H_2$ is a finite $S$-group scheme of multiplicative type, $M$ is of multiplicative type (\S\ref{sec-multiplicative-type-closure-properties}) and even isotrivial (Lemma \ref{lem-isotriviality-closure-properties}(ii)). Thus, we may take another connected finite \'etale Galois cover $\widetilde S'\rightarrow\widetilde S$ and assume that $M$ is split by $\widetilde S'$. Using Lemma \ref{lemm-resolution-multiplicative-type}, we find a resolution of $M$ by $S$-tori which are also split by $\widetilde S'$:
\begin{equation}
\label{eq-resolution-of-reductive-groups-multiplicative-type}
1 \rightarrow M \rightarrow T_1 \rightarrow Q \rightarrow 1,
\end{equation}
where either
\benuma
	\item $T_1$ is flasque and $Q$ is quasi-trivial; or
	\item $T_1$ is quasi-trivial and $Q$ is coflasque.
\eenum

Let us form the push-out of the extension \eqref{eq-resolution-of-reductive-groups-extension} along the map $M \rightarrow T_1$. This gives rise to a central extension of $G$ by $T_1$ that fits into a commutative diagram
$$
\begin{tikzcd}
	&  & 1 \ar[d] & & \\
	1 \ar[r] & M \ar[r]\ar[d] & G'^{\der} \times P \ar[r]\ar[d, "\alpha"] & G \ar[r]\ar[d, "\cong"] & 1 \\
	1 \ar[r] & T_1 \ar[r] & G' \ar[r]\ar[d] & G \ar[r] & 1 \\
	& & Q \ar[d] & & \\
	& & 1 & &
\end{tikzcd}
$$
By construction, the map $\alpha$ induces an isomorphism on derived subgroups. Hence, the morphism $G' \rightarrow G$ induces the given central isogeny $f \colon G'^{\der} \rightarrow G^{\der}$ on derived subgroups. Recall that the formation of radicals is preserved under quotient maps. (This statement may be verified over geometric points, where it is \cite[Expos\'e XIX, Section 1.7]{SGA3IIInew}.) Hence $\rad(G')$ is a quotient of the torus $T_1\times P$. Since the latter is split by $\widetilde S'$, so is $\rad(G')$  (Lemma \ref{lem-isotriviality-closure-properties}(i)).

Finally, we show that the two types of resolutions \eqref{eq-resolution-of-reductive-groups-multiplicative-type} give rise to the two conditions in the statement of Theorem \ref{theo-resolution-of-reductive-groups}. Indeed, write $T_2 := \corad(G')$. Since $T_2$ is a quotient of $\rad(G')$, it is also split by $\widetilde S'$. We have a short exact sequence of $S$-tori split by $\widetilde S'$:
\begin{equation}
\label{eq-resolution-of-reductive-groups-torus-quotient}
1 \rightarrow P \rightarrow T_2 \rightarrow Q \rightarrow 1.
\end{equation}
Since $P$ is quasi-trivial and $Q$ is at least coflasque, Lemma \ref{lemm-flasque-coflasque-split} shows that \eqref{eq-resolution-of-reductive-groups-torus-quotient} splits. In particular, $T_2$ is quasi-trivial (resp.~coflasque) whenever $Q$ is.
\end{proof}


\csub[An application to the Grothendieck--Serre conjecture]






We use \Cref{theo-resolution-of-reductive-groups} to reduce the Grothendieck--Serre \cref{conj:Grothendieck-Serre} to the case when the group $G$ has a simply connected derived subgroup. This argument is suggested to me by K.~\v{C}esnavi\v{c}ius.

\bprop \lab{prop:GS-app}
Let $R$ be a regular local ring, let $K \ce \Frac(R)$, let $G$ be a reductive $R$-group scheme, and consider the pullback map
\begin{equation}
\label{eq-grothendieck-serre-map}
 H^1_{\et}(R, G) \rightarrow  H^1_{\et}(K, G).
\end{equation}
If this map has trivial kernel whenever $G$ is replaced by some central extension $G'$ of $G$ whose derived subgroup $G'^{\der}$ is simply connected, then it has trivial kernel for $G$ itself.
\eprop

\begin{proof}
By Theorem \ref{theo-resolution-of-reductive-groups}, we may find a central extension of reductive $R$-group schemes
$$
1 \rightarrow T_1 \rightarrow G' \rightarrow G \rightarrow 1,
$$
such that the map $G' \rightarrow G$ induces the simply connected cover $G'^{\der} \cong (G^{\der})^{\mathrm{sc}} \rightarrow G^{\der}$ on the derived subgroups and $T_1$ is a quasi-trivial torus.\footnote{Information about $T_2 := \corad(G')$ is not needed for this proof.} By \eqref{eqn:coho-seq}, this extension gives rise to the following map of exact sequences of pointed sets:
$$
\begin{tikzcd}
	 H^1(R, T_1) \ar[d, "\alpha_1"] \ar[r] &  H^1(R, G') \ar[d, "\beta'"]\ar[r] &  H^1(R, G) \ar[r]\ar[d, "\beta"] &  H^2(R, T_1) \ar[d, "\alpha_2"] \\
	 H^1(\Frac(R), T_1) \ar[r] & H^1(\Frac(R), G') \ar[r] &  H^1(\Frac(R), G) \ar[r] &  H^2(\Frac(R), T_1).
\end{tikzcd}
$$
Since $T_1$ is a quasi-trivial torus, by Lemma \ref{lemm-quasitrivial-torus-explicit}, it is isomorphic to a finite product of tori of the form $\Res_{R'/R}(\mathbb G_m)$ for some finite \'etale maps $R \ra R'$. In particular, Hilbert 90 implies that $\alpha_1$ is an isomorphism between singletons. Grothendieck's theorem on the Brauer group \cite[Corollaire~1.8]{Gro68b} implies that $\alpha_2$ is injective. Therefore, if $\beta'$ has trivial kernel, then so does $\beta$, as desired.
\end{proof}

\end{appendix}

\end{document}